arXiv:math/0501237v4 [math.RA] 16 Jul 2010

# Vector Space Over Division Ring

Aleks Kleyn

Abstract. A system of linear equations over a division ring has properties similar to properties of a system of linear equations over a field. Even noncommutativity of a product creates a new picture the properties of system of linear equations and of vector space over division ring have a close relationship.

## Contents

## 1. Representation of Universal Algebra

**Definition 1.1.** Suppose we defined the structure of $\Omega_2$-algebra on the set $M$ ([1, 7]). We call the endomorphism of $\Omega_2$-algebra

$$t : M \to M$$

**transformation of universal algebra** $M$.[1] □

We denote $\delta$ identical transformation.

Aleks_Kleyn@MailAPS.org.
http://sites.google.com/site/AleksKleyn/.
http://arxiv.org/a/kleyn_a_1.
http://AleksKleyn.blogspot.com/.

[1] If the set of operations of $\Omega_2$-algebra is empty, then $t$ is a map.

**Definition 1.2.** Transformations is **left-side transformation** or **$T\star$-transformation** if it acts from left

$$u' = tu$$

We denote ${}^\star M$ the set of $T\star$-transformations of set $M$. □

**Definition 1.3.** Transformations is **right-side transformations** or **$\star T$-transformation** if it acts from right

$$u' = ut$$

We denote $M^\star$ the set of nonsingular $\star T$-transformations of set $M$. □

**Definition 1.4.** Suppose we defined the structure of $\Omega_1$-algebra on the set ${}^\star M$ ([1]). Let $A$ be $\Omega_1$-algebra. We call homomorphism

$$f : A \rightarrow {}^\star M \tag{1.1}$$

**left-side** or **$T\star$-representation of $\Omega_1$-algebra $A$ in $\Omega_2$-algebra $M$** □

**Definition 1.5.** Suppose we defined the structure of $\Omega_1$-algebra on the set $M^\star$ ([1]). Let $A$ be $\Omega_1$-algebra. We call homomorphism

$$f : A \rightarrow M^\star$$

**right-side** or **$\star T$-representation of $\Omega_1$-algebra $A$ in $\Omega_2$-algebra $M$** □

We extend to representation theory convention described in remark [6]-2.15. We can write duality principle in the following form

**Theorem 1.6** (duality principle)**.** *Any statement which holds for $T\star$-representation of $\Omega_1$-algebra $A$ holds also for $\star T$-representation of $\Omega_1$-algebra $A$.*

*Remark* 1.7. There exist two forms of notation for transformation of $\Omega_2$-algebra $M$. In operational notation, we write the transformation $A$ as either $Aa$ which corresponds to the $T\star$-transformation or $aA$ which corresponds to the $\star T$-transformation. In functional notation, we write the transformation $A$ as $A(a)$ regardless of the fact whether this is $T\star$-transformation or this is $\star T$-transformation. This notation is in agreement with duality principle.

This remark serves as a basis for the following convention. When we use functional notation we do not make a distinction whether this is $T\star$-transformation or this is $\star T$-transformation. We denote ${}^* M$ the set of transformations of $\Omega_2$-algebra $M$. Suppose we defined the structure of $\Omega_1$-algebra on the set ${}^* M$. Let $A$ be $\Omega_1$-algebra. We call homomorphism

$$f : A \rightarrow {}^* M \tag{1.2}$$

**representation of $\Omega_1$-algebra $A$ in $\Omega_2$-algebra $M$.**

Correspondence between operational notation and functional notation is unambiguous. We can select any form of notation which is convenient for presentation of particular subject. □

Diagram

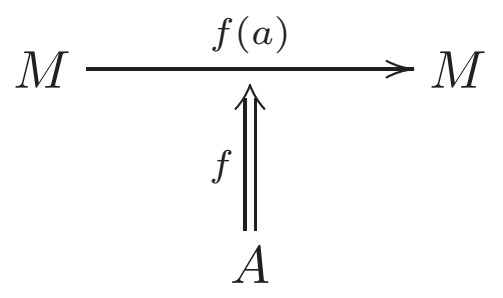

means that we consider the representation of $\Omega_1$-algebra $A$. The map $f(a)$ is image of $a \in A$.

**Definition 1.8.** Suppose map (1.2) is an isomorphism of the $\Omega_1$-algebra $A$ into ${}^*M$. Then the representation of the $\Omega_1$-algebra $A$ is called **effective**. □

*Remark* 1.9. Suppose the $T\star$-representation of $\Omega_1$-algebra is effective. Then we identify an element of $\Omega_1$-algebra and its image and write $T\star$-transormation caused by element $a \in A$ as

$$v' = av$$

Suppose the $\star T$-representation of $\Omega_1$-algebra is effective. Then we identify an element of $\Omega_1$-algebra and its image and write $\star T$-transormation caused by element $a \in A$ as

$$v' = va$$

□

**Definition 1.10.** We call a representation of $\Omega_1$-algebra **transitive** if for any $a, b \in V$ exists such $g$ that

$$a = f(g)(b)$$

We call a representation of $\Omega_1$-algebra **single transitive** if it is transitive and effective. □

**Theorem 1.11.** *$T\star$-representation is single transitive if and only if for any $a, b \in M$ exists one and only one $g \in A$ such that $a = f(g)(b)$*

*Proof.* Corollary of definitions 1.8 and 1.10. □

## 2. Morphism of Representations of Universal Algebra

**Theorem 2.1.** *Let $A$ and $B$ be $\Omega_1$-algebras. Representation of $\Omega_1$-algebra $B$*

$$g : B \to {}^\star M$$

*and homomorphism of $\Omega_1$-algebra*

$$h : A \to B \tag{2.1}$$

*define representation $f$ of $\Omega_1$-algebra $A$*

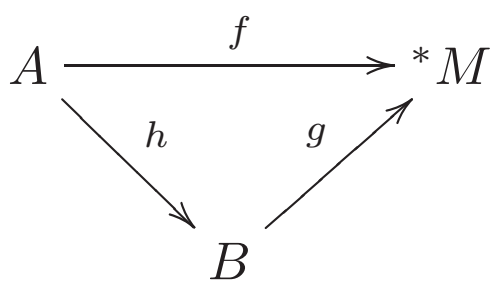

*Proof.* Since mapping $g$ is homomorphism of $\Omega_1$-algebra $B$ into $\Omega_1$-algebra ${}^*M$, the mapping $f$ is homomorphism of $\Omega_1$-algebra $A$ into $\Omega_1$-algebra ${}^*M$. □

Considering representations of $\Omega_1$-algebra in $\Omega_2$-algebras $M$ and $N$, we are interested in a mapping that preserves the structure of representation.

**Definition 2.2.** Let

$$f : A \to {}^*M$$

be representation of $\Omega_1$-algebra $A$ in $\Omega_2$-algebra $M$ and

$$g : B \to {}^*N$$

be representation of $\Omega_1$-algebra $B$ in $\Omega_2$-algebra $N$. Tuple of maps

$$(r : A \to B, R : M \to N) \tag{2.2}$$

such, that

- $r$ is homomorphism of $\Omega_1$-algebra
- $R$ is homomorphism of $\Omega_2$-algebra
-

$$R\circ f(a)=g(r(a))\circ R \tag{2.3}$$

is called **morphism of representations from $f$ into $g$**. We also say that **morphism of representations of $\Omega_1$-algebra in $\Omega_2$-algebra** is defined. □

For any $m\in M$ equation (2.3) has form

$$R(f(a)(m))=g(r(a))(R(m)) \tag{2.4}$$

*Remark* 2.3. We may consider a pair of maps $r$, $R$ as map

$$F:A\cup M\rightarrow B\cup N$$

such that

$$F(A)=B \quad F(M)=N$$

Therefore, hereinafter we will say that we have the map $(r,R)$. □

*Remark* 2.4. Let us consider morphism of representations (2.2). We denote elements of the set $B$ by letter using pattern $b\in B$. However if we want to show that $b$ is image of element $a\in A$, we use notation $r(a)$. Thus equation

$$r(a)=r(a)$$

means that $r(a)$ (in left part of equation) is image $a\in A$ (in right part of equation). Using such considerations, we denote element of set $N$ as $R(m)$. We will follow this convention when we consider correspondences between homomorphisms of $\Omega_1$-algebra and mappings between sets where we defined corresponding representations.

There are two ways to interpret (2.4)

- Let transformation $f(a)$ map $m\in M$ into $f(a)(m)$. Then transformation $g(r(a))$ maps $R(m)\in N$ into $R(f(a)(m))$.
- We represent morphism of representations from $f$ into $g$ using diagram

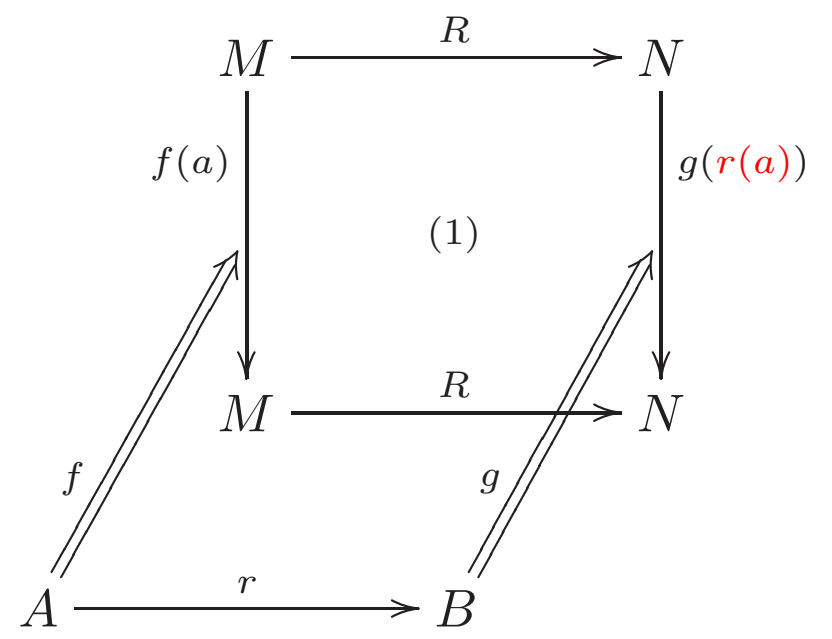

From (2.3) it follows that diagram (1) is commutative.

□

**Theorem 2.5.** *Let us consider representation*

$$f:A\rightarrow {}^*M$$

*of $\Omega_1$-algebra $A$ and representation*

$$g:B\rightarrow {}^*N$$

*of $\Omega_1$-algebra $B$. Morphism*

$$h : A \longrightarrow B \qquad H : M \longrightarrow N$$

*of representations from $f$ into $g$ satisfies equation*

$$H \circ \omega(f(a_1), ..., f(a_n)) = \omega(g(h(a_1)), ..., g(h(a_n))) \circ H \tag{2.5}$$

*for any $n$-ary operation $\omega$ of $\Omega_1$-algebra.*

*Proof.* Since $f$ is homomorphism, we have

$$H \circ \omega(f(a_1), ..., f(a_n)) = H \circ f(\omega(a_1, ..., a_n)) \tag{2.6}$$

From (2.3) and (2.6) it follows that

$$H \circ \omega(f(a_1), ..., f(a_n)) = g(h(\omega(a_1, ..., a_n))) \circ H \tag{2.7}$$

Since $h$ is homomorphism, from (2.7) it follows that

$$H \circ \omega(f(a_1), ..., f(a_n)) = g(\omega(h(a_1), ..., h(a_n))) \circ H \tag{2.8}$$

Since $g$ is homomorphism, (2.5) follows from (2.8). □

**Theorem 2.6.** *Let the map*

$$h : A \longrightarrow B \qquad H : M \longrightarrow N$$

*be morphism from representation*

$$f : A \to {}^*M$$

*of $\Omega_1$-algebra $A$ into representation*

$$g : B \to {}^*N$$

*of $\Omega_1$-algebra $B$. If representation $f$ is effective, then the map*

$${}^*H : {}^*M \to {}^*N$$

*defined by equation*

$${}^*H(f(a)) = g(h(a)) \tag{2.9}$$

*is homomorphism of $\Omega_1$-algebra.*

*Proof.* Because representation $f$ is effective, then for given transformation $f(a)$ element $a$ is determined uniquely. Therefore, transformation $g(h(a))$ is properly defined in equation (2.9).

Since $f$ is homomorphism, we have

$${}^*H(\omega(f(a_1), ..., f(a_n))) = {}^*H(f(\omega(a_1, ..., a_n))) \tag{2.10}$$

From (2.9) and (2.10) it follows that

$${}^*H(\omega(f(a_1), ..., f(a_n))) = g(h(\omega(a_1, ..., a_n))) \tag{2.11}$$

Since $h$ is homomorphism, from (2.11) it follows that

$${}^*H(\omega(f(a_1), ..., f(a_n))) = g(\omega(h(a_1), ..., h(a_n))) \tag{2.12}$$

Since $g$ is homomorphism,

$${}^*H(\omega(f(a_1), ..., f(a_n))) = \omega(g(h(a_1)), ..., g(h(a_n))) = \omega({}^*H(f(a_1)), ..., {}^*H(f(a_n)))$$

follows from (2.12). Therefore, the map ${}^*H$ is homomorphism of $\Omega_1$-algebra. □

**Theorem 2.7.** *Given single transitive representation*

$$f : A \to {}^*M$$

*of $\Omega_1$-algebra $A$ and single transitive $T\star$-epresentation*

$$g : B \to {}^*N$$

*of $\Omega_1$-algebra $B$, there exists morphism*

$$h : A \longrightarrow B \qquad H : M \longrightarrow N$$

*of representations from $f$ into $g$.*

*Proof.* Let us choose homomorphism $h$. Let us choose element $m \in M$ and element $n \in N$. To define map $H$, let us consider following diagram

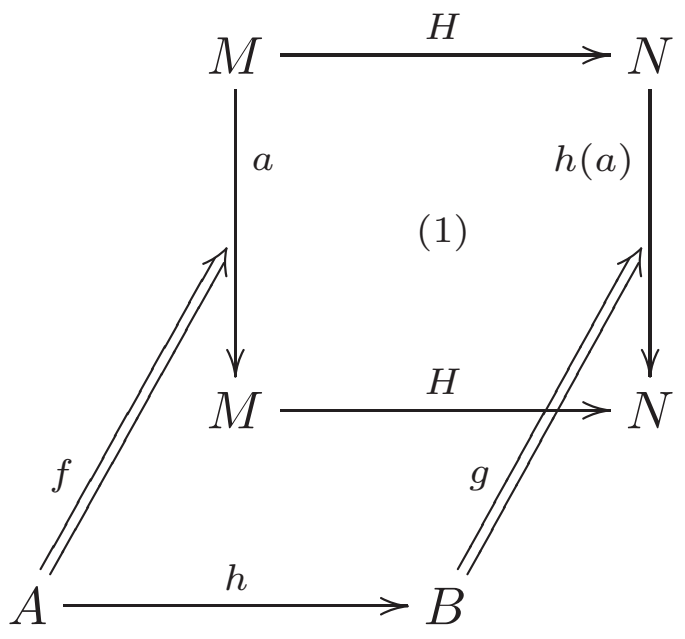

From commutativity of diagram (1), it follows that

$$H(am) = h(a)H(m)$$

For arbitrary $m' \in M$, we defined unambiguously $a \in A$ such that $m' = am$. Therefore, we defined mapping $H$ which satisfies to equation (2.3). □

**Theorem 2.8.** *Let*

$$f : A \to {}^*M$$

*be single transitive representation of $\Omega_1$-algebra $A$ and*

$$g : B \to {}^*N$$

*be single transitive representation of $\Omega_1$-algebra $B$. Given homomorphism of $\Omega_1$-algebra*

$$h : A \longrightarrow B$$

*let us consider a map*

$$H : M \longrightarrow N$$

*such that $(h, H)$ is morphism of representations from $f$ into $g$. This map is unique up to choice of image $n = H(m) \in N$ of given element $m \in M$.*

*Proof.* From proof of theorem 2.7 it follows that choice of homomorphism $h$ and elements $m \in M$, $n \in N$ uniquely defines the map $H$. □

**Theorem 2.9.** *Given single transitive representation*

$$f : A \to {}^*M$$

*of $\Omega_1$-algebra $A$, for any endomorphism of $\Omega_1$-algebra $A$ there exists endomorphism*

$$p : A \longrightarrow A \qquad P : M \longrightarrow M$$

*of representation* $f$.

*Proof.* Let us consider following diagram

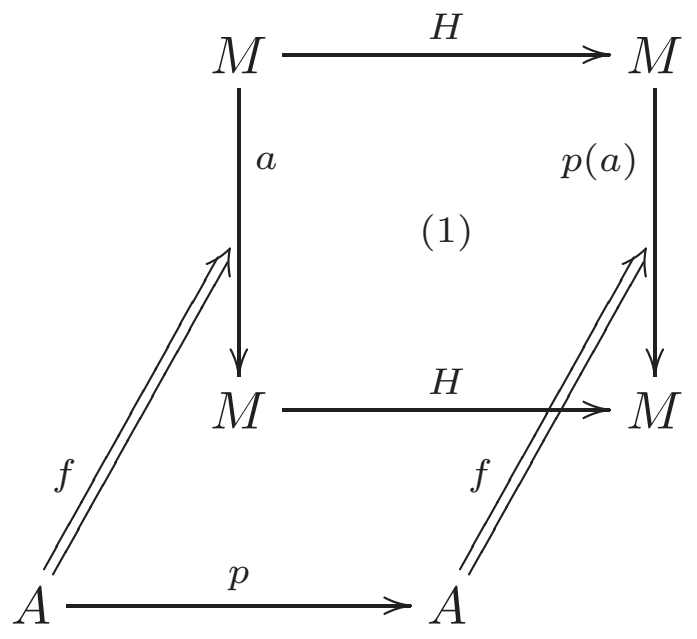

Statement of theorem is corollary of theorem 2.7. □

**Theorem 2.10.** *Let*

$$f : A \to {}^*M$$

*be representation of* $\Omega_1$*-algebra* $A$*,*

$$g : B \to {}^*N$$

*be representation of* $\Omega_1$*-algebra* $B$*,*

$$h : C \to {}^*L$$

*be representation of* $\Omega_1$*-algebra* $C$*. Given morphisms of representations of* $\Omega_1$*-algebra*

$$p : A \longrightarrow B \qquad P : M \longrightarrow N$$

$$q : B \longrightarrow C \qquad Q : N \longrightarrow L$$

*There exists morphism of representations of* $\Omega_1$*-algebra*

$$r : A \longrightarrow C \qquad R : M \longrightarrow L$$

*where* $r = qp$*,* $R = QP$*. We call morphism* $(r, R)$ *of representations from* $f$ *into* $h$ **product of morphisms** $(p, P)$ **and** $(q, Q)$ **of representations of universal algebra**.

*Proof.* We represent statement of theorem using diagram

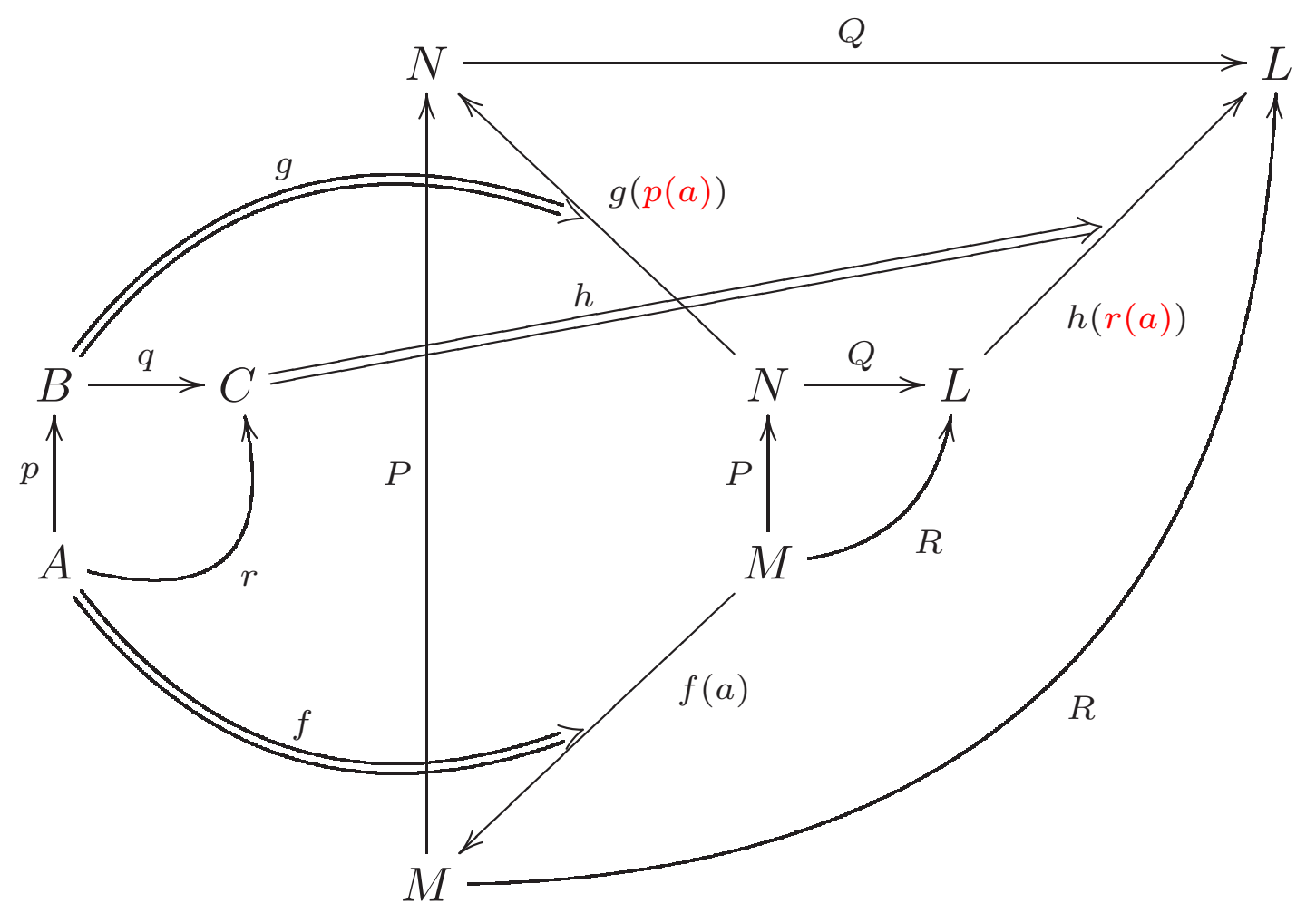

Map $r$ is homomorphism of $\Omega_1$-algebra $A$ into $\Omega_1$-algebra $C$. We need to show that tuple of maps $(r, R)$ satisfies to (2.3):

$$\begin{aligned} R(f(a)m) &= QP(f(a)m) \\ &= Q(g(p(a))P(m)) \\ &= h(qp(a))QP(m)) \\ &= h(r(a))R(m) \end{aligned}$$

□

**Definition 2.11.** Let $\mathcal{A}$ be category of $\Omega_1$-algebras. We define **category** $T\star\mathcal{A}$ **of** $T\star$**-representations of** $\Omega_1$**-algebra from category** $\mathcal{A}$. $T\star$-representations of $\Omega_1$-algebra are objects of this category. Morphisms of $T\star$-representations of $\Omega_1$-algebra are morphisms of this category. □

**Theorem 2.12.** *Endomorphisms of representation $f$ form semigroup.*

*Proof.* From theorem 2.10, it follows that the product of endomorphisms $(p, P)$, $(r, R)$ of the representation $f$ is endomorphism $(pr, PR)$ of the representation $f$. □

**Definition 2.13.** Let us define equivalence $S$ on the set $M$. Transformation $f$ is called **coordinated with equivalence** $S$, when $f(m_1) \equiv f(m_2)(\mathrm{mod} S)$ follows from condition $m_1 \equiv m_2(\mathrm{mod} S)$. □

**Theorem 2.14.** *Let us consider equivalence $S$ on set $M$. Let us consider $\Omega_1$-algebra on set ${}^*M$. Since transformations are coordinated with equivalence $S$, we can define the structure of $\Omega_1$-algebra on the set ${}^*(M/S)$.*

*Proof.* Let $h = \mathrm{nat}\ S$. If $m_1 \equiv m_2(\mathrm{mod} S)$, then $h(m_1) = h(m_2)$. Since $f \in {}^*M$ is coordinated with equivalence $S$, then $h(f(m_1)) = h(f(m_2))$. This allows to define transformation $F$ according to rule

$$F([m]) = h(f(m))$$

Let $\omega$ be n-ary operation of $\Omega_1$-algebra. Suppose $f_1$, ..., $f_n \in {}^*M$ and

$$F_1([m]) = h(f_1(m)) \quad ... \quad F_n([m]) = h(f_n(m))$$

We define operation on the set $^\star(M/S)$ according to rule

$$\omega(F_1,...,F_n)[m]=h(\omega(f_1,...,f_n)m)$$

This definition is proper because $\omega(f_1,...,f_n)\in {}^\star M$ and is coordinated with equivalence $S$. □

**Theorem 2.15.** *Let*

$$f:A\to {}^*M$$

*be representation of $\Omega_1$-algebra $A$,*

$$g:B\to {}^*N$$

*be representation of $\Omega_1$-algebra $B$. Let*

$$r:A\longrightarrow B \qquad R:M\longrightarrow N$$

*be morphism of representations from $f$ into $g$. Suppose*

$$s=rr^{-1} \qquad S=RR^{-1}$$

*Then there exist decompositions of $r$ and $R$, which we describe using diagram*

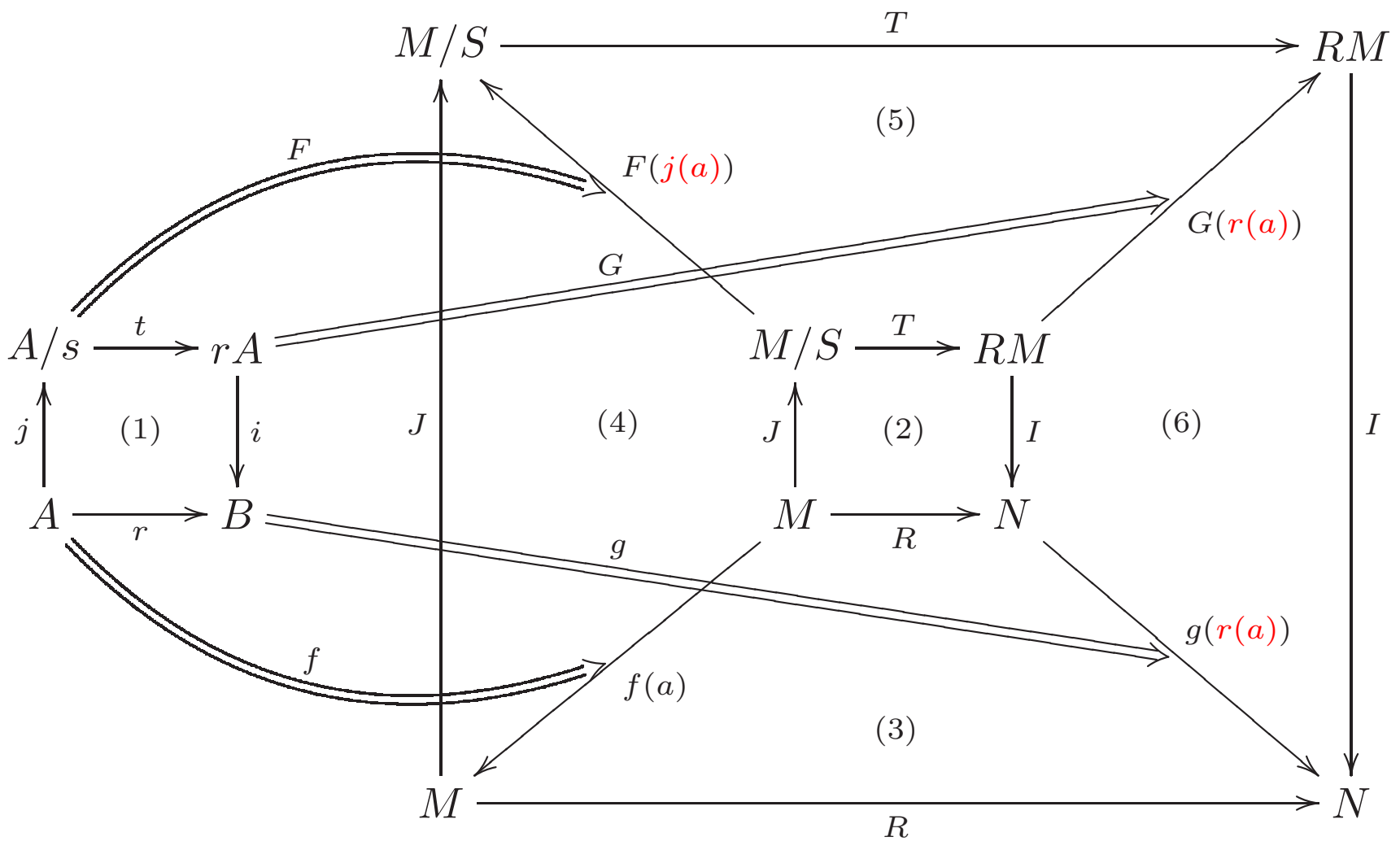

(1) $s=\ker r$ *is a congruence on $A$. There exists decompositions of homomorphism $r$*

$$r=itj \tag{2.13}$$

*$j=\operatorname{nat} s$ is the natural homomorphism*

$$j(a)=j(a) \tag{2.14}$$

*$t$ is isomorphism*

$$r(a)=t(j(a)) \tag{2.15}$$

*$i$ is the inclusion mapping*

$$r(a)=i(r(a)) \tag{2.16}$$

(2) *$S = \ker R$ is an equivalence on $M$. There exists decompositions of homomorphism $R$*

$$R = ITJ \tag{2.17}$$

*$J = \text{nat } S$ is surjection*

$$J(m) = J(m) \tag{2.18}$$

*$T$ is bijection*

$$R(m) = T(J(m)) \tag{2.19}$$

*$I$ is the inclusion mapping*

$$R(m) = I(R(m)) \tag{2.20}$$

(3) *$F$ is $T\star$-representation of $\Omega_1$-algebra $A/s$ in $M/S$*
(4) *$G$ is $T\star$-representation of $\Omega_1$-algebra $rA$ in $RM$*
(5) *$(j, J)$ is morphism of representations $f$ and $F$*
(6) *$(t, T)$ is morphism of representations $F$ and $G$*
(7) *$(t^{-1}, T^{-1})$ is morphism of representations $G$ and $F$*
(8) *$(i, I)$ is morphism of representations $G$ and $g$*
(9) *There exists decompositions of morphism of representations*

$$(r, R) = (i, I)(t, T)(j, J) \tag{2.21}$$

*Proof.* Existence of diagrams (1) and (2) follows from theorem II.3.7 ([7], p. 60).

We start from diagram (4).

Let $m_1 \equiv m_2 (\text{mod } S)$. Then

$$R(m_1) = R(m_2) \tag{2.22}$$

Since $a_1 \equiv a_2 (\text{mod} s)$, then

$$r(a_1) = r(a_2) \tag{2.23}$$

Therefore, $j(a_1) = j(a_2)$. Since $(r, R)$ is morphism of representations, then

$$R(f(a_1)(m_1)) = g(r(a_1))(R(m_1)) \tag{2.24}$$

$$R(f(a_2)(m_2)) = g(r(a_2))(R(m_2)) \tag{2.25}$$

From (2.22), (2.23), (2.24), (2.25), it follows that

$$R(f(a_1)(m_1)) = R(f(a_2)(m_2)) \tag{2.26}$$

From (2.26) it follows

$$f(a_1)(m_1) \equiv f(a_2)(m_2) (\text{mod} S) \tag{2.27}$$

and, therefore,

$$J(f(a_1)(m_1)) = J(f(a_2)(m_2)) \tag{2.28}$$

From (2.28) it follows that we defined map

$$F(j(a))(J(m)) = J(f(a)(m)) \tag{2.29}$$

reasonably and this map is transformation of set $M/S$.

From equation (2.27) (in case $a_1 = a_2$) it follows that for any $a$ transformation is coordinated with equivalence $S$. From theorem 2.14 it follows that we defined structure of $\Omega_1$-algebra on the set ${}^*(M/S)$. Let us consider $n$-ary operation $\omega$ and $n$ transformations

$$F(j(a_i))(J(m)) = J(f(a_i)(m))) \quad i = 1, ..., n$$

of the set $M/S$. We assume

$$\omega(F(j(a_1)), ..., F(j(a_n)))(J(m)) = J(\omega(f(a_1), ..., f(a_n)))(m))$$

Therefore, map $F$ is representations of $\Omega_1$-algebra $A/s$.

From (2.29) it follows that $(j, J)$ is morphism of representations $f$ and $F$ (the statement (5) of the theorem).

Let us consider diagram (5).

Since $T$ is bijection, then we identify elements of the set $M/S$ and the set $MR$, and this identification has form

$$T(J(m)) = R(m) \tag{2.30}$$

We can write transformation $F(j(a))$ of the set $M/S$ as

$$F(j(a)) : J(m) \to F(j(a))(J(m)) \tag{2.31}$$

Since $T$ is bijection, we define transformation

$$T(J(m)) \to T(F(j(a))(J(m))) \tag{2.32}$$

of the set $RM$. Transformation (2.32) depends on $j(a) \in A/s$. Since $t$ is bijection, we identify elements of the set $A/s$ and the set $rA$, and this identification has form

$$t(j(a)) = r(a) \tag{2.33}$$

Therefore, we defined map

$$G : rA \to {}^\star RM$$

according to equation

$$G(t(j(a)))(T(J(m))) = T(F(j(a))(J(m))) \tag{2.34}$$

Let us consider $n$-ary operation $\omega$ and $n$ transformations

$$G(r(a_i))(R(m)) = T(F(j(a_i))(J(m))) \quad i = 1, ..., n$$

of space $RM$. We assume

$$\omega(G(r(a_1)), ..., G(r(a_n)))(R(m)) = T(\omega(F(j(a_1), ..., F(j(a_n)))(J(m))) \tag{2.35}$$

According to (2.34) operation $\omega$ is defined reasonably on the set ${}^\star RM$. Therefore, the map $G$ is representations of $\Omega_1$-algebra.

From (2.34) it follows that $(t, T)$ is morphism of representations $F$ and $G$ (the statement (6) of the theorem).

Since $T$ is bijection, then from equation (2.30) it follows that

$$J(m) = T^{-1}(R(m)) \tag{2.36}$$

We can write transformation $G(r(a))$ of the set $RM$ as

$$G(r(a)) : R(m) \to G(r(a))(R(m)) \tag{2.37}$$

Since $T$ is bijection, we define transformation

$$T^{-1}(R(m)) \to T^{-1}(G(r(a))(R(m))) \tag{2.38}$$

of the set $M/S$. Transformation (2.38) depends on $r(a)\in rA$. Since $t$ is bijection, then from equation (2.33) it follows that

$$j(a)=t^{-1}(r(a)) \tag{2.39}$$

Since, by construction, diagram (5) is commutative, then transformation (2.38) coincides with transformation (2.31). We can write the equation (2.35) as

$$T^{-1}(\omega(G(r(a_1)),...,G(r(a_n)))(R(m)))=\omega(F(j(a_1),...,F(j(a_n)))(J(m)) \tag{2.40}$$

Therefore $(t^{-1},T^{-1})$ is morphism of representations $G$ and $F$ (the statement (7) of the theorem).

Diagram (6) is the most simple case in our prove. Since map $I$ is immersion and diagram (2) is commutative, we identify $n\in N$ and $R(m)$ when $n\in \mathrm{Im}R$. Similarly, we identify corresponding transformations.

$$g'(i(r(a)))(I(R(m)))=I(G(r(a))(R(m))) \tag{2.41}$$

$$\omega(g'(r(a_1)),...,g'(r(a_n)))(R(m))=I(\omega(G(r(a_1),...,G(r(a_n)))(R(m)))$$

Therefore, $(i,I)$ is morphism of representations $G$ and $g$ (the statement (8) of the theorem).

To prove the statement (9) of the theorem we need to show that defined in the proof representation $g'$ is congruent with representation $g$, and operations over transformations are congruent with corresponding operations over ${}^*N$.

$$\begin{aligned}
g'(i(r(a)))(I(R(m)))&=I(G(r(a))(R(m))) && \text{by (2.41)}\\
&=I(G(t(j(a)))(T(J(m)))) && \text{by (2.15), (2.19)}\\
&=IT(F(j(a))(J(m))) && \text{by (2.34)}\\
&=ITJ(f(a)(m)) && \text{by (2.29)}\\
&=R(f(a)(m)) && \text{by (2.17)}\\
&=g(r(a))(R(m)) && \text{by (2.3)}
\end{aligned}$$

$$\begin{aligned}
\omega(G(r(a_1)),...,G(r(a_n)))(R(m))&=T(\omega(F(j(a_1),...,F(j(a_n)))(J(m)))\\
&=T(F(\omega(j(a_1),...,j(a_n)))(J(m)))\\
&=T(F(j(\omega(a_1,...,a_n)))(J(m)))\\
&=T(J(f(\omega(a_1,...,a_n))(m)))
\end{aligned}$$

□

**Definition 2.16.** Let

$$f:A\rightarrow {}^*M$$

be representation of $\Omega_1$-algebra $A$,

$$g:B\rightarrow {}^*N$$

be representation of $\Omega_1$-algebra $B$. Let

$$r:A\longrightarrow B\qquad R:M\longrightarrow N$$

be morphism of representations from $f$ into $g$ such that $f$ is isomorphism of $\Omega_1$-algebra and $g$ is isomorphism of $\Omega_2$-algebra. Then map $(r,R)$ is called **isomorphism of repesentations**. □

**Theorem 2.17.** *In the decomposition* (2.21)*, the map* $(t, T)$ *is isomorphism of representations* $F$ *and* $G$*.*

*Proof.* The statement of the theorem is corollary of definition 2.16 and statements (6) and (7) of the theorem 2.15. □

From theorem 2.15 it follows that we can reduce the problem of studying of morphism of representations of $\Omega_1$-algebra to the case described by diagram

(2.42)

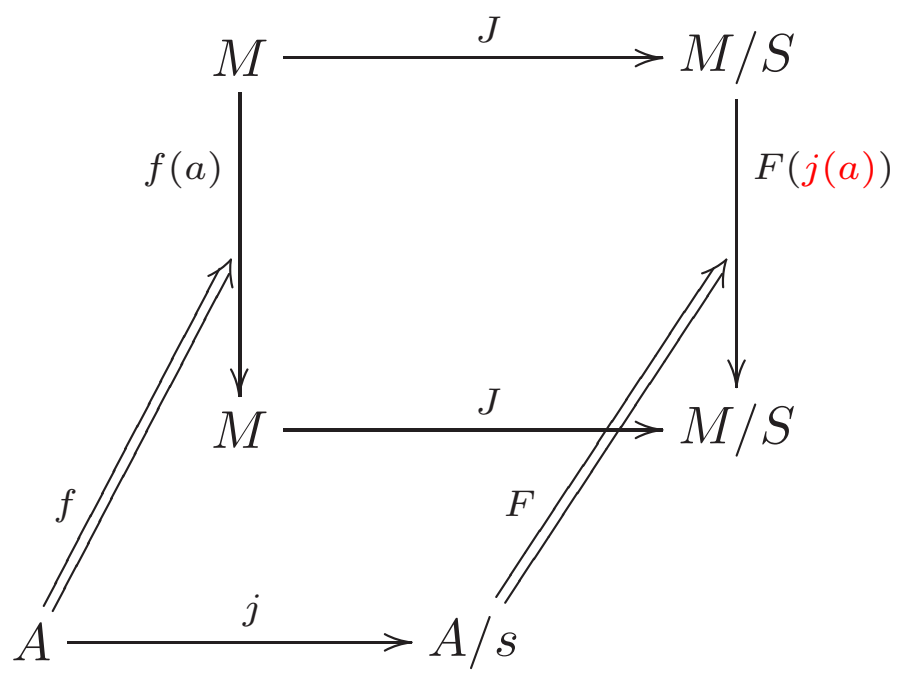

**Theorem 2.18.** *We can supplement diagram* (2.42) *with representation* $F_1$ *of* $\Omega_1$*-algebra* $A$ *into set* $M/S$ *such that diagram*

(2.43)

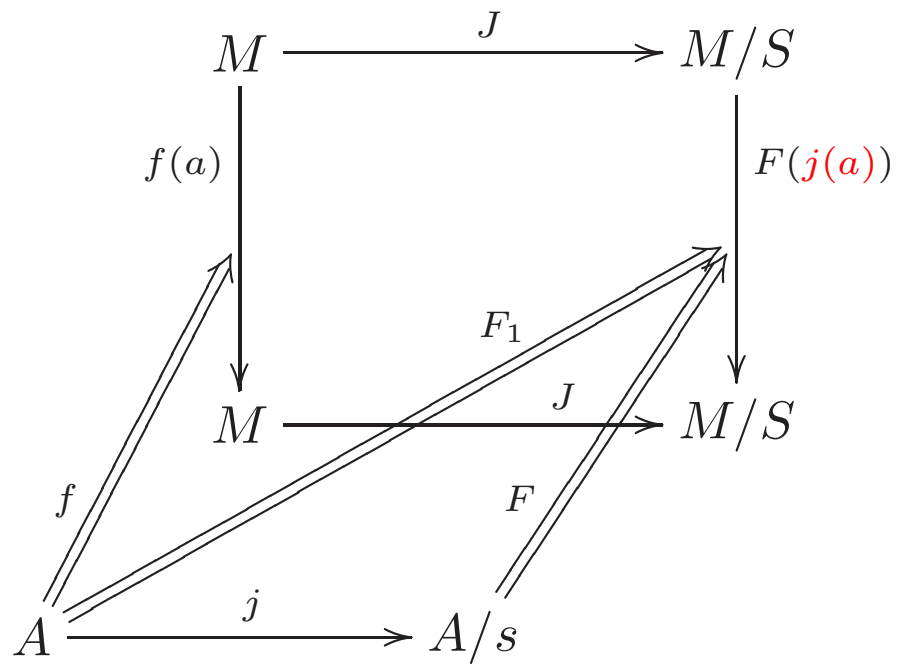

*is commutative. The set of transformations of representation* $F$ *and the set of transformations of representation* $F_1$ *coincide.*

*Proof.* To prove theorem it is enough to assume

$$F_1(a) = F(j(a))$$

Since map $j$ is surjection, then $\mathrm{Im}F_1 = \mathrm{Im}F$. Since $j$ and $F$ are homomorphisms of $\Omega_1$-algebra, then $F_1$ is also homomorphism of $\Omega_1$-algebra. □

Theorem 2.18 completes the series of theorems dedicated to the structure of morphism of representations $\Omega_1$-algebra. From these theorems it follows that we can simplify task of studying of morphism of representations $\Omega_1$-algebra and not go beyond morphism of representations of form

$$id : A \longrightarrow A \qquad R : M \longrightarrow N$$

In this case we identify morphism of $(id, R)$ representations of $\Omega_1$-algebra and map $R$. We will use diagram

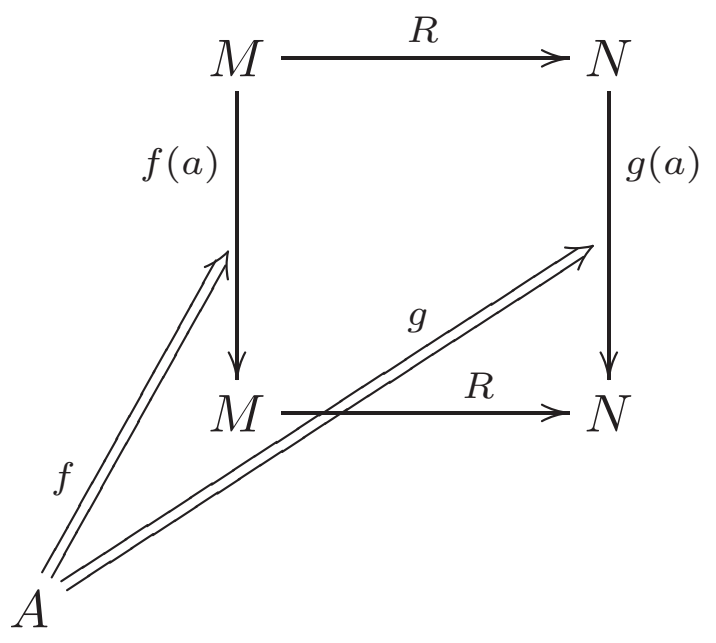

to represent morphism $(id, R)$ of representations of $\Omega_1$-algebra. From diagram it follows

$$R \circ f(a) = g(a) \circ R \tag{2.44}$$

By analogy with definition 2.11. we give following definition.

**Definition 2.19.** We define **category** $T \star A$ **$T\star$-representations of $\Omega_1$-algebra** $A$. $T\star$-representations of $\Omega_1$-algebra $A$ are objects of this category. Morphisms $(id, R)$ of $T\star$-representations of $\Omega_1$-algebra $A$ are morphisms of this category. □

## 3. Automorphism of Representation of Universal Algebra

**Definition 3.1.** Let

$$f : A \to {}^*M$$

be representation of $\Omega_1$-algebra $A$ in $\Omega_2$-algebra $M$. The morphism of representations of $\Omega_1$-algebra

$$(r : A \to A, R : M \to M)$$

such, that $r$ is endomorphism of $\Omega_1$-algebra and $R$ is endomorphism of $\Omega_2$-algebra is called **endomorphism of representation** $f$. □

**Definition 3.2.** Let

$$f : A \to {}^*M$$

be representation of $\Omega_1$-algebra $A$ in $\Omega_2$-algebra $M$. The morphism of representations of $\Omega_1$-algebra

$$(r : A \to A, R : M \to M)$$

such, that $r$ is automorphism of $\Omega_1$-algebra and $R$ is automorphism of $\Omega_2$-algebra is called **automorphism of representation** $f$. □

**Theorem 3.3.** *Let*

$$f : A \to {}^*M$$

*be representation of $\Omega_1$-algebra $A$ in $\Omega_2$-algebra $M$. The set of automorphisms of the representation $f$ forms* **loop** $\mathfrak{A}(f)$ .[2]

[2]Look [3], p. 24, [2] for definition of loop.

*Proof.* Let $(r,R)$, $(p,P)$ be automorphisms of the representation $f$. According to definition 3.2 maps $r$, $p$ are automorphisms of $\Omega_1$-algebra $A$ and maps $R$, $P$ are automorphisms of $\Omega_2$-algebra $M$. According to theorem II.3.2 ([7], p. 57), the map $rp$ is automorphism of $\Omega_1$-algebra $A$ and the map $RP$ is automorphism of $\Omega_2$-algebra $M$. From the theorem 2.10 and the definition 3.2, it follows that product of automorphisms $(rp,RP)$ of representation $f$ is automorphism of the representation $f$.

Let $(r,R)$ be an automorphism of the representation $f$. According to definition 3.2 the map $r$ is automorphism of $\Omega_1$-algebra $A$ and the map $R$ is automorphism of $\Omega_2$-algebra $M$. Therefore, the map $r^{-1}$ is automorphism of $\Omega_1$-algebra $A$ and the map $R^{-1}$ is automorphism of $\Omega_2$-algebra $M$. The equation (2.4) is true for automorphism $(r,R)$. Assume $a'=r(a)$, $m'=R(m)$. Since $r$ and $R$ are automorphisms then $a=r^{-1}(a')$, $m=R^{-1}(m')$ and we can write (2.4) in the form

$$R(f(r^{-1}(a'))(R^{-1}(m')))=g(a')(m') \tag{3.1}$$

Since the map $R$ is automorphism of $\Omega_2$-algebra $M$, then from the equation (3.1) it follows that

$$f(r^{-1}(a'))(R^{-1}(m'))=R^{-1}(g(a')(m')) \tag{3.2}$$

The equation (3.2) corresponds to the equation (2.4) for the map $(r^{-1},R^{-1})$. Therefore, map $(r^{-1},R^{-1})$ of the representation $f$. □

*Remark* 3.4. It is evident that the set of automorphisms of $\Omega_1$-algebra $A$ also forms loop. Of course, it is attractive to assume that the set of automorphisms forms a group. Since the product of automorphisms $f$ and $g$ is automorphism $fg$, then automorphisms $(fg)h$ and $f(gh)$ are defined. However, it does not follow from this statement that

$$(fg)h=f(gh)$$

□

## 4. Vector Space

To define $T\star$-representation

$$f:R\to {}^\star\overline{M}$$

of ring $R$ on the set $\overline{M}$ we need to define the structure of the ring on the set ${}^\star\overline{M}$.[3]

[3]Is it possible to define an addition on the set ${}^\star\overline{M}$, if this operation is not defined on the set $\overline{M}$. The answer on this question is positive.

Let $\overline{M}=B\cup C$ and let $F:B\to C$ be one to one map. We define the set ${}^\star\overline{M}$ of $T\star$-transformations of the set $\overline{M}$ according to the following rule. Let $V\subseteq B$. Let the $T\star$-transformation $F_V$ be given by

$$F_Vx=\begin{cases}x & x\in B\backslash V\\ Fx & x\in V\\ x & x\in C\backslash F(V)\\ F^{-1}x & x\in F(V)\end{cases}$$

We define sum of $T\star$-transformations according rule

$$F_V+F_W=F_{V\triangle W}$$
$$V\triangle W=(V\cup W)\backslash(V\cap W)$$

It is evident that

$$F_\emptyset+F_V=F_V$$
$$F_V+F_V=F_\emptyset$$

**Theorem 4.1.** *$T\star$-representation $f$ of the ring $R$ on the set $\overline{M}$ is defined iff $T\star$-representations of multiplicative and additive groups of the ring $R$ are defined and these $T\star$-representations hold relationship*

$$f(a(b+c)) = f(a)f(b) + f(a)f(c)$$

*Proof.* Theorem follows from definition 1.4. □

**Definition 4.2.** $\overline{M}$ is a $R\star$**-module** over a ring $R$ if $\overline{M}$ is an Abelian group and there exists $T\star$-representation of ring $R$. □

According to our notation $R\star$-module is **left module** and $\star R$**-module** is **right module**.

Since a field is a special case of a ring, vector space over the field has more properties then module over the ring. It is very hard, if possible at all, to extend definitions, which work in a vector space, to a module over an arbitrary ring. A definition of a basis and dimension of vector space are closely linked with the possibility of finding a solution of a linear equation in a ring. Properties of the linear equation in division ring are close to properties of the linear equation in field. This is why we hope that properties of vector space over division ring are close to properties of the vector space over the field.

**Theorem 4.3.** *$T\star$-representation of the division ring $D$ is* **effective** *iff $T\star$-representation of its multiplicative group is effective.*

*Proof.* Suppose

$$f : D \to {}^\star M \tag{4.1}$$

is $T\star$-representation of the division ring $D$. Suppose elements $a$, $b$ of of the multiplicative group cause the same $T\star$-transformation. Then

$$f(a)m = f(b)m \tag{4.2}$$

for any $m \in M$. Performing transformation $f(a^{-1})$ on both sides of the equation (4.2), we obtain

$$m = f(a^{-1})(f(b)m) = f(a^{-1}b)m$$

□

According to the remark 1.9, since the representation of the division ring is effective, we identify an element of the division ring and $T\star$-transformation corresponding to this element.

**Definition 4.4.** $\overline{V}$ is a $D\star$**-vector space** over a division ring $D$ if $\overline{V}$ is an Abelian group and there exists effective $T\star$-representation of division ring $D$. □

**Theorem 4.5.** *Following conditions hold for $D\star$-vector space:*

- **associative law**

$$(ab)\overline{m} = a(b\overline{m}) \tag{4.3}$$

- **distributive law**

$$a(\overline{m}+\overline{n}) = a\overline{m} + a\overline{n} \tag{4.4}$$

$$(a+b)\overline{m} = a\overline{m} + b\overline{m} \tag{4.5}$$

---

Therefore, the map $F_\emptyset$ is zero of the addition, and the set ${}^\star\overline{M}$ is the Abelian group.

- **unitarity law**

$$1\overline{m} = \overline{m} \tag{4.6}$$

*for any* $a, b \in D$*,* $\overline{m}, \overline{n} \in \overline{V}$*. We call* $T\star$*-representation* $D\star$**-product of vector over scalar**.

*Proof.* Since $T\star$-transformation $a$ is automorphism of the Abelian group, we obtain the equation (4.4). Since $T\star$-representation is homomorphism of the aditive group of division ring $D$, we obtain the equation (4.5). Since $T\star$-representation is $T\star$-representation of the multiplicative group of division ring $D$, we obtain the equations (4.3) and (4.6). □

According to our notation $D\star$-vector space is **left vector space** and $\star D$**-vector space** is **right vector space**.

**Definition 4.6.** Let $\overline{V}$ be a $D\star$-vector space over a division ring $D$. Set of vectors $\overline{N}$ is a **subspace of** $D\star$**-vector space** $\overline{V}$ if

$$\overline{a} + \overline{b} \in \overline{N}$$
$$k\overline{a} \in \overline{N}$$
$$\overline{a}, \overline{b} \in \overline{N} \quad k \in D$$

□

**Example 4.7.** Let $D^m_n$ be set of $m \times n$ matrices over division ring $D$. We define addition

$$a + b = \left(a^j_i\right) + \left(b^j_i\right) = \left(a^j_i + b^j_i\right)$$

and product over scalar

$$da = d\left(a^j_i\right) = \left(da^j_i\right)$$

$a = 0$ iff $a^j_i = 0$ for any $i$, $j$. We can verify directly that $D^m_n$ is a $D\star$-vector space. when product is defined from left. Otherwise $D^m_n$ is $\star D$-vector space. Vector space $D^m_n$ is called $D\star$**-matrices vector space**. □

## 5. Vector Space Type

The product of vector over scalar is asymmetric because the product is defined for objects of different sets. However we see difference between $D\star$- and $\star D$-vector space only when we work with coordinate representation. When we speak vector space is $D\star$- or $\star D$- we point out how we multiply coordinates of vector over elements of division ring: from left or right.

**Definition 5.1.** Suppose $u$, $v$ are vectors of $D\star$-vector space $\overline{V}$. We call vector $w$ $D\star$**-linear composition of vectors** $u$ and $v$ when we can write $w = au + bv$ where $a$ and $b$ are scalars. □

We can extend definition of the linear composition on any finite set of vectors. Using generalized indexes to enumerate vectors we can represent set of vectors as one dimensional matrix. We use the convention that we represent any set of vectors of the vector space or as $_*$-row either as $^*$-row. This representation defines type of notation of linear composition. Getting this representation in $D\star$- or $\star D$-vector space we get four different models of vector space.

For an opportunity to show without change of the notation what kind of vector space ($D\star$- or $\star D$-) we study we introduce new notation. We call symbol $D\star$ **vector space type** and this symbol means that we study $D\star$-vector space over division ring $D$. The symbol of product in the type of vector space points to matrix operation used in the linear composition.

**Example 5.2.** Let ${}^*$-row

$$\overline{a}=\begin{pmatrix} {}^1\overline{a} \\ ... \\ {}^n\overline{a} \end{pmatrix}$$

represent the set of vectors ${}^i\overline{a}$, $i\in I$, of $D_*{}^*$**-vector space** $\overline{V}$ and ${}_*$-row

$$c=\begin{pmatrix} c_1 & ... & c_n \end{pmatrix}$$

represent the set of scalars $c_i$, $i\in I$. Then we can write the linear composition of vectors ${}^i\overline{a}$ as

$$c_i\ {}^i\overline{a}=c_*{}^*\overline{a}$$

We use notation ${}_{D_*{}^*}\overline{V}$ when we want to tell that $\overline{V}$ is $D_*{}^*$-vector space. □

**Example 5.3.** Let ${}_*$-row

$$\overline{a}=\begin{pmatrix} {}_1\overline{a} & ... & {}_n\overline{a} \end{pmatrix}$$

represent the set of vectors ${}_i\overline{a}$, $i\in I$, of $D^*{}_*$**-vector space** $\overline{V}$ and ${}^*$-row

$$c=\begin{pmatrix} c^1 \\ ... \\ c^n \end{pmatrix}$$

represent the set of scalars $c^i$, $i\in I$. Then we can write the linear composition of vectors ${}_i\overline{a}$ as

$$c^i\ {}_i\overline{a}=c^*{}_*\overline{a}$$

We use notation ${}_{D^*{}_*}\overline{V}$ when we want to tell that $\overline{V}$ is $D^*{}_*$-vector space. □

**Example 5.4.** Let ${}^*$-row

$$\overline{a}=\begin{pmatrix} \overline{a}^1 \\ ... \\ \overline{a}^n \end{pmatrix}$$

represent the set of vectors $\overline{a}^i$, $i\in I$, of ${}^*{}_*D$**-vector space** $\overline{V}$ and ${}_*$-row

$$c=\begin{pmatrix} {}_1c & ... & {}_nc \end{pmatrix}$$

represent the set of scalars ${}_ic$, $i\in I$. Then we can write the linear composition of vectors $\overline{a}^i$ as

$$\overline{a}^i\ {}_ic=\overline{a}^*{}_*c$$

We use notation $\overline{V}_{{}^*{}_*D}$ when we want to tell that $\overline{V}$ is ${}^*{}_*D$-vector space. □

**Example 5.5.** Let ${}_*$-row

$$\overline{a} = \begin{pmatrix} \overline{a}_1 & ... & \overline{a}_n \end{pmatrix}$$

represent the set of vectors $\overline{A}_i,\ i \in I,$ of ${}_*{}^*D$**-vector space** $\overline{V}$ and ${}^*$-row

$$c = \begin{pmatrix} {}^1c \\ ... \\ {}^nc \end{pmatrix}$$

represent the set of scalars ${}_ic,\ i \in I$. Then we can write the linear composition of vectors $\overline{a}_i$ as

$$\overline{a}_i\ {}^ic = \overline{a}_*{}^*c$$

We use notation $\overline{V}_{{}_*{}^*D}$ when we want to tell that $\overline{V}$ is ${}_*{}^*D$-vector space. □

*Remark* 5.6. We extend to vector space and its type convention described in remark [6]-2.15. For instance, we execute operations in expression

$$A_*{}^*B_*{}^*v\lambda$$

from left to right. This corresponds to the ${}_*{}^*D$-vector space. However we can execute product from right to left. In custom notation this expression is

$$\lambda v^*{}_*B^*{}_*A$$

and corresponds to $D^*{}_*$-vector space. Similarly, reading this expression from down up we get expression

$$A^*{}_*B^*{}_*v\lambda$$

corresponding to ${}^*{}_*D$-vector space. □

## 6. ${}_*{}^*D$-Basis of Vector Space

**Definition 6.1.** Vectors $\overline{A}_i,\ i \in I,$ of ${}_*{}^*D$-vector space $\overline{V}$ are ${}_*{}^*D$**-linearly independent** if $c = 0$ follows from the equation

$$\overline{A}_*{}^*c = 0$$

Otherwise vectors $\overline{A}_i$ are ${}_*{}^*D$**-linearly dependent**. □

**Definition 6.2.** We call set of vectors $\overline{\overline{e}} = (\overline{e}_i, i \in I)$ a ${}_*{}^*D$**-basis for vector space** if vectors $\overline{e}_i$ are ${}_*{}^*D$-linearly independent and adding to this system any other vector we get a new system which is ${}_*{}^*D$-linearly dependent. □

**Theorem 6.3.** *If $\overline{\overline{e}}$ is a ${}_*{}^*D$-basis of vector space $\overline{V}$ then any vector $\overline{v} \in \overline{V}$ has one and only one expansion*

$$\overline{v} = \overline{e}_*{}^*v \tag{6.1}$$

*relative to this ${}_*{}^*D$-basis.*

*Proof.* Because system of vectors $\overline{e}_i$ is a maximal set of ${}_*{}^*D$-linearly independent vectors the system of vectors $\overline{v}, \overline{e}_i$ is $D_*{}^*$-linearly dependent and in equation

$$\overline{v}b + \overline{e}_*{}^*c = 0 \tag{6.2}$$

at least $b$ is different from 0. Then equation

$$\overline{v} = \overline{e}_*{}^*(-cb^{-1}) \tag{6.3}$$

follows from (6.2). (6.1) follows from (6.3).

Assume we get another expansion

$$\overline{v}=\overline{e}_*{}^*v' \tag{6.4}$$

We subtract (6.1) from (6.4) and get

$$0=\overline{e}_*{}^*(v'-v)$$

Because vectors $\overline{e}_i$ are ${}_*{}^*D$-linearly independent we get

$$v'-v=0$$

□

**Definition 6.4.** We call the matrix $v$ in expansion (6.1) **coordinate matrix of vector $\overline{v}$ in ${}_*{}^*D$-basis $\overline{\overline{e}}$** and we call its elements **coordinates of vector $\overline{v}$ in ${}_*{}^*D$-basis $\overline{\overline{e}}$**. □

**Theorem 6.5.** *Set of coordinates $a$ of vector $\overline{a}$ relative ${}_*{}^*D$-basis $\overline{\overline{e}}$ forms ${}_*{}^*D$-vector space $D^n$ isomorphic ${}_*{}^*D$-vector space $\overline{V}$. This ${}_*{}^*D$-vector space is called* **coordinate ${}_*{}^*D$-vector space**. *This isomorphism is called* **coordinate ${}_*{}^*D$-isomorphism**.

*Proof.* Suppose vectors $\overline{a}$ and $\overline{b}\in\overline{V}$ have expansion

$$\overline{a}=\overline{e}_*{}^*a$$
$$\overline{b}=\overline{e}_*{}^*b$$

relative basis $\overline{\overline{e}}$. Then

$$\overline{a}+\overline{b}=\overline{e}_*{}^*a+\overline{e}_*{}^*b=\overline{e}_*{}^*(a+b)$$
$$\overline{a}m=(\overline{e}_*{}^*a)m=\overline{e}_*{}^*(ma)$$

for any $m\in D$. Thus, operations in a vector space are defined by coordinates

$${}^i(a+b)={}^ia+{}^ib$$
$${}^i(am)={}^iam$$

This completes the proof. □

**Example 6.6.** Let $\overline{\overline{e}}=(\overline{e}_i,i\in I,|I|=n)$ be a ${}_*{}^*D$-basis of vector space $D^n$. The coordinate matrix

$$a=\begin{pmatrix}{}^1a\\ \dots\\ {}^na\end{pmatrix}=({}^ia,i\in I)$$

of vector $\overline{a}$ in ${}_*{}^*D$-basis $\overline{\overline{e}}$ is called **${}_*{}^*D$-vector**.[4] We call vector space $D^n$ **${}^*$-rows ${}_*{}^*D$-vector space**.

Let ${}_*$-row

$$\overline{A}=\begin{pmatrix}\overline{A}_1 & \dots & \overline{A}_m\end{pmatrix}=(\overline{A}_j,j\in J) \tag{6.5}$$

be set of vectors. Vectors $\overline{A}_j$ have expansion

$$\overline{A}_j=\overline{e}_*{}^*A_j$$

[4] ${}_*{}^*D$-vector is an analogue of column vector. We can also call it ${}^*$-row $D\star$-vector.

If we substitute coordinate matrices of vectors $\overline{A}_j$ in the matrix (6.5) we get matrix

$$A = \begin{pmatrix} \begin{pmatrix} {}^1A_1 \\ ... \\ {}^nA_1 \end{pmatrix} & ... & \begin{pmatrix} {}^1A_m \\ ... \\ {}^nA_m \end{pmatrix} \end{pmatrix} = \begin{pmatrix} {}^1A_1 & ... & {}^1A_m \\ ... & ... & ... \\ {}^nA_1 & ... & {}^nA_m \end{pmatrix} = ({}^iA_j)$$

We call the matrix $A$ **coordinate matrix of set of vectors** $(\overline{A}_j, j \in J)$ in basis $\overline{\overline{e}}$ and we call its elements **coordinates of set of vectors** $(\overline{A}_j, j \in J)$ in basis $\overline{\overline{e}}$.

${}_*$-Row

$$\overline{f} = \begin{pmatrix} \overline{f}_1 & ... & \overline{f}_n \end{pmatrix} = (\overline{f}_j, j \in I)$$

represents the ${}_*{}^*D$**-basis** $\overline{\overline{f}}$ **of** ${}^*$**-rows vector space** $D^n$. We tell that coordinate matrix $f$ of set of vectors $(\overline{f}_j, j \in I)$ defines coordinates ${}^if_j$ of basis $\overline{\overline{f}}$ relative basis $\overline{\overline{e}}$. □

**Example 6.7.** Let $\overline{\overline{e}} = ({}^j\overline{e}, j \in J, |J| = m)$ be a $D_*{}^*$-basis of vector space $D_n$. The coordinate matrix

$$a = \begin{pmatrix} a_1 & ... & a_m \end{pmatrix} = (a_j, j \in J)$$

of vector $\overline{a}$ in $D_*{}^*$-basis $\overline{\overline{e}}$ is called $D_*{}^*$**-vector**.[5] We call vector space $D_n$ ${}_*$**-rows** $D_*{}^*$**-vector space**.

Let ${}^*$-row

$$\overline{A} = \begin{pmatrix} {}^1\overline{A} \\ ... \\ {}^n\overline{A} \end{pmatrix} = ({}^i\overline{A}, i \in I) \tag{6.6}$$

be set of vectors. Vectors ${}^i\overline{A}$ have expansion

$${}^i\overline{A} = {}^iA_*{}^*\overline{e}$$

If we substitute coordinate matrices of vectors ${}^i\overline{A}$ in the matrix (6.6) we get matrix

$$A = \begin{pmatrix} \begin{pmatrix} {}^1A_1 & ... & {}^1A_m \end{pmatrix} \\ ... \\ \begin{pmatrix} {}^nA_1 & ... & {}^nA_m \end{pmatrix} \end{pmatrix} = \begin{pmatrix} {}^1A_1 & ... & {}^1A_m \\ ... & ... & ... \\ {}^nA_1 & ... & {}^nA_m \end{pmatrix} = ({}^iA_j)$$

We call the matrix $A$ **coordinate matrix of set of vectors** $({}^i\overline{A}, i \in M)$ in basis $\overline{\overline{e}}$ and we call its elements **coordinates of set of vectors** $({}^i\overline{A}, i \in M)$ in basis $\overline{\overline{e}}$.

${}^*$-Row

$$\overline{f} = \begin{pmatrix} {}^1\overline{f} & ... & {}^n\overline{f} \end{pmatrix} = ({}^i\overline{f}, i \in J)$$

represents the $D_*{}^*$**-basis** $\overline{\overline{f}}$ **of** ${}_*$**-rows vector space** $D_n$. We tell that coordinate matrix $f$ of set of vectors $({}^i\overline{f}, i \in J)$ defines coordinates ${}^if_j$ of basis $\overline{\overline{f}}$ relative basis $\overline{\overline{e}}$. □

Since we express linear composition using matrices we can extend the duality principle to the vector space theory. We can write duality principle in one of the following forms

[5] $D_*{}^*$-vector is an analogue of row vector. We can also call it ${}_*$-row $D\star$-vector

**Theorem 6.8** (duality principle)**.** *Let $\mathfrak{A}$ be true statement about vector spaces. If we exchange the same time*

- $D_*{}^*$*-vector and* $D^*{}_*$*-vector*
- ${}_*{}^*D$*-vector and* ${}^*{}_*D$*-vector*
- ${}_*{}^*$*-product and* ${}^*{}_*$*-product*

*then we soon get true statement.*

**Theorem 6.9** (duality principle)**.** *Let $\mathfrak{A}$ be true statement about vector spaces. If we exchange the same time*

- $D_*{}^*$*-vector and* ${}_*{}^*D$*-vector or* $D^*{}_*$*-vector and* ${}^*{}_*D$*-vector*
- ${}_*{}^*$*-quasideterminant and* ${}^*{}_*$*-quasideterminant*

*then we soon get true statement.*

## 7. ${}_*{}^*D$-Linear Map of Vector Spaces

**Definition 7.1.** Suppose $\overline{V}$ is ${}_*{}^*S$-vector space. Suppose $\overline{U}$ is ${}_*{}^*T$-vector space. Morphism

$$f : S \longrightarrow T \qquad \overline{A} : \overline{V} \longrightarrow \overline{U}$$

of $T\star$-representations of division ring in Abelian group is called $({}_*{}^*S, {}_*{}^*T)$**-linear map of vector spaces**. □

By theorem 2.15 studying $({}_*{}^*S, {}_*{}^*T)$-linear map we can consider case $S = T$.

**Definition 7.2.** Suppose $\overline{V}$ and $\overline{W}$ are ${}_*{}^*D$-vector spaces. We call map

$$\overline{A} : \overline{V} \to \overline{W} \tag{7.1}$$

${}_*{}^*D$**-linear map of vector spaces** if[6]

$$\overline{A}(\overline{m}_*{}^*a) = \overline{A}(\overline{m})_*{}^*a \tag{7.2}$$

for any ${}^ia \in D$, $\overline{m}_i \in \overline{V}$. □

**Theorem 7.3.** *Let*

$$\overline{\overline{f}} = (\overline{f}_i, i \in I)$$

*be a* ${}_*{}^*D$*-basis of vector space* $\overline{V}$ *and*

$$\overline{\overline{e}} = (\overline{e}_j, j \in J)$$

*be a* ${}_*{}^*D$*-basis of vector space* $\overline{U}$. *Then* ${}_*{}^*D$*-linear map* (7.1) *of vector spaces has presentation*

$$b = A_*{}^*a \tag{7.3}$$

*relative to selected bases. Here*

- $a$ *is coordinate matrix of vector* $\overline{a}$ *relative the* ${}_*{}^*D$*-basis* $\overline{\overline{f}}$
- $b$ *is coordinate matrix of vector*
$$\overline{b} = \overline{A}(\overline{a})$$
*relative the* ${}_*{}^*D$*-basis* $\overline{\overline{e}}$
- $A$ *is coordinate matrix of set of vectors* $(\overline{A}(\overline{f}_i))$ *in* ${}_*{}^*D$*-basis* $\overline{\overline{e}}$ *called* **matrix of** ${}_*{}^*D$**-linear mapping** *relative bases* $\overline{\overline{f}}$ *and* $\overline{\overline{e}}$

[6]Expression $\overline{A}(\overline{m})_*{}^*a$ means expression $\overline{A}(\overline{m}_i)\ {}^ia$

*Proof.* Vector $\overline{a}\in\overline{V}$ has expansion

$$\overline{a}=\overline{f}{}_*{}^*a$$

relative to ${}_*{}^*D$-basis $\overline{\overline{f}}$. Vector $\overline{b}=f(\overline{a})\in\overline{U}$ has expansion

$$\overline{b}=\overline{e}{}_*{}^*b \tag{7.4}$$

relative to ${}_*{}^*D$-basis $\overline{\overline{e}}$.

Since $\overline{A}$ is a ${}_*{}^*D$-linear map, from (7.2) it follows that

$$\overline{b}=\overline{A}(\overline{a})=\overline{A}(\overline{f}{}_*{}^*a)=\overline{A}(\overline{f}){}_*{}^*a \tag{7.5}$$

$\overline{A}(\overline{f}_i)$ is also a vector of $\overline{U}$ and has expansion

$$\overline{A}(\overline{f}_i)=\overline{e}{}_*{}^*A_i=\overline{e}_j{}^jA_i \tag{7.6}$$

relative to basis $\overline{\overline{e}}$. Combining (7.5) and (7.6) we get

$$\overline{b}=\overline{e}{}_*{}^*A{}_*{}^*a \tag{7.7}$$

(7.3) follows from comparison of (7.4) and (7.7) and theorem 6.3. □

On the basis of theorem 7.3 we identify the ${}_*{}^*D$-linear map (7.1) of vector spaces and the matrix of its presentation (7.3).

**Theorem 7.4.** *Let*

$$\overline{\overline{f}}=(\overline{f}_i,i\in I)$$

*be a ${}_*{}^*D$-basis of vector space $\overline{V}$,*

$$\overline{\overline{e}}=(\overline{e}_j,j\in J)$$

*be a $D_*{}^*$-basis of vector space $\overline{U}$, and*

$$\overline{\overline{g}}=(\overline{g}_l,l\in L)$$

*be a ${}_*{}^*D$-basis of vector space $\overline{W}$. Suppose diagram of mappings*

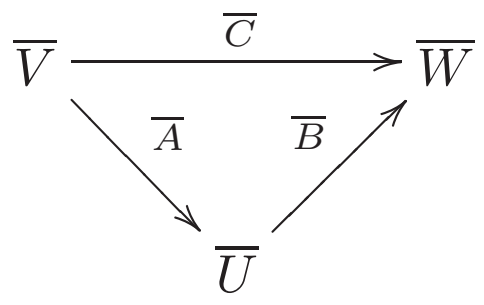

*is commutative diagram where ${}_*{}^*D$-linear mapping $\overline{A}$ has presentation*

$$b=A{}_*{}^*a \tag{7.8}$$

*relative to selected bases and ${}_*{}^*D$-linear mapping $\overline{B}$ has presentation*

$$c=B{}_*{}^*b \tag{7.9}$$

*relative to selected bases. Then mapping $\overline{C}$ is ${}_*{}^*D$-linear mapping and has presentation*

$$c=B{}_*{}^*A{}_*{}^*a \tag{7.10}$$

*relative to selected bases.*

*Proof.* The mapping $\overline{C}$ is ${}_*{}^*D$-linear mapping because

$$\begin{aligned}\overline{C}_*{}^*(\overline{f}_*{}^*a) &= (\overline{A}_*{}^*\overline{B})_*{}^*(\overline{f}_*{}^*a)\\ &= \overline{B}_*{}^*(\overline{A}_*{}^*(\overline{f}_*{}^*a))\\ &= \overline{B}_*{}^*(\overline{e}_*{}^*(A_*{}^*a))\\ &= \overline{g}_*{}^*(B_*{}^*(A_*{}^*a))\\ &= \overline{g}_*{}^*((B_*{}^*A)_*{}^*a)\\ &= \overline{g}_*{}^*(C_*{}^*a)\end{aligned}$$

Equation (7.10) follows from substituting (7.8) into (7.9). □

Presenting ${}_*{}^*D$-linear map as ${}_*{}^*$-product we can rewrite (7.2) as

$$\overline{A}_*{}^*(\overline{a}k) = (\overline{A}_*{}^*\overline{a})k \tag{7.11}$$

We can express the statement of the theorem 7.4 in the next form

$$\overline{B}_*{}^*(\overline{A}_*{}^*\overline{a}) = (\overline{B}_*{}^*\overline{A})_*{}^*\overline{a} \tag{7.12}$$

Equations (7.11) and (7.12) represent the **associative law for ${}_*{}^*D$-linear maps of vector spaces**. This allows us writing of such expressions without using of brackets.

Equation (7.3) is coordinate notation for ${}_*{}^*D$-linear map. Based theorem 7.3 non coordinate notation also can be expressed using ${}_*{}^*$-product

$$\overline{b} = \overline{A}_*{}^*\overline{a} = \overline{A}_*{}^*\overline{f}_*{}^*a = \overline{e}_*{}^*A_*{}^*a \tag{7.13}$$

If we substitute equation (7.13) into theorem 7.4, then we get chain of equations

$$\begin{aligned}\overline{c} &= \overline{B}_*{}^*\overline{b} = \overline{B}_*{}^*\overline{e}_*{}^*b = \overline{g}_*{}^*B_*{}^*b\\ \overline{c} &= \overline{B}_*{}^*\overline{A}_*{}^*\overline{a} = \overline{B}_*{}^*\overline{A}_*{}^*\overline{f}_*{}^*a = \overline{g}_*{}^*B_*{}^*A_*{}^*a\end{aligned}$$

*Remark* 7.5. One can easily see from the an example of ${}_*{}^*D$-linear map how theorem 2.15 makes our reasoning simpler in study of the morphism of $T\star$-representations of $\Omega$-algebra. In the framework of this remark, we agree to call the theory of ${}_*{}^*D$-linear mappings reduced theory, and theory stated in this remark is called enhanced theory.

Suppose $\overline{V}$ is ${}_*{}^*S$-vector space. Suppose $\overline{U}$ is ${}_*{}^*T$-vector space. Suppose

$$r : S \longrightarrow T \qquad \overline{A} : \overline{V} \longrightarrow \overline{U}$$

is $({}_*{}^*S, {}_*{}^*T)$-linear map of vector spaces. Let

$$\overline{\overline{f}} = (\overline{f}_i, i \in I)$$

be a ${}_*{}^*S$-basis of vector space $\overline{V}$ and

$$\overline{\overline{e}} = (\overline{e}_j, j \in J)$$

be a ${}_*{}^*T$-basis of vector space $\overline{U}$.

From definitions 7.1 and 2.2 it follows

$$\overline{b} = \overline{A}(\overline{a}) = \overline{A}(\overline{f}_*{}^*a) = \overline{A}(\overline{f})_*{}^*r(a) \tag{7.14}$$

$\overline{A}(\overline{f}_i)$ is also a vector of $\overline{U}$ and has expansion

$$\overline{A}(\overline{f}_i) = \overline{e}_*{}^*A_i = \overline{e}_j{}^jA_i \tag{7.15}$$

relative to basis $\overline{\overline{e}}$. Combining (7.14) and (7.15), we get

$$\overline{b} = \overline{e}_*{}^*A_*{}^*r(a) \tag{7.16}$$

Suppose $\overline{W}$ is ${}_*{}^*D$-vector space. Suppose

$$p : T \longrightarrow D \qquad \overline{B} : \overline{U} \longrightarrow \overline{W}$$

is $({}_*{}^*T, {}_*{}^*D)$-linear map of vector spaces. Let

$$\overline{\overline{g}} = (\overline{g}_l, l \in L)$$

be ${}_*{}^*D$-basis of vector space $\overline{W}$. Then, according to (7.16), the product of $({}_*{}^*S, {}_*{}^*T)$-linear map $(r, \overline{A})$ and $({}_*{}^*T, {}_*{}^*D)$-linear map $(p, \overline{B})$ has form

$$\overline{c} = \overline{h}_*{}^*B_*{}^*p(A)_*{}^*pr(a) \tag{7.17}$$

Comparison of equations (7.10) and (7.17) that extended theory of linear maps is more complicated then reduced theory.

If we need we can use extended theory, however we will not get new results comparing with reduced theory. At the same time plenty of details makes picture less clear and demands permanent attention. □

## 8. System of ${}_*{}^*D$-Linear Equations

**Definition 8.1.** Let $\overline{V}$ be a ${}_*{}^*D$-vector space and $\{\overline{A}_i \in \overline{V}, i \in I\}$ be set of vectors. ${}_*{}^*D$**-linear span in vector space** is set $\mathrm{span}(\overline{A}_i, i \in I)$ of vectors ${}_*{}^*D$-linearly dependent on vectors $\overline{A}_i$. □

**Theorem 8.2.** *Let span*$(\overline{A}_i, i \in I)$ *be* ${}_*{}^*D$*-linear span in vector space* $\overline{V}$*. Then span*$(\overline{A}_i, i \in I)$ *is subspace of* $\overline{V}$*.*

*Proof.* Suppose

$$\overline{b} \in \mathrm{span}(\overline{A}_i, i \in I)$$
$$\overline{c} \in \mathrm{span}(\overline{A}_i, i \in I)$$

According to definition 8.1

$$\overline{b} = \overline{A}_*{}^*b$$
$$\overline{c} = \overline{A}_*{}^*c$$

Then

$$\overline{b} + \overline{c} = \overline{A}_*{}^*b + \overline{A}_*{}^*c = \overline{A}_*{}^*(b + c) \in \mathrm{span}(\overline{A}_i, i \in I)$$
$$\overline{b}k = (\overline{A}_*{}^*b)k = \overline{A}_*{}^*(bk) \in \mathrm{span}(\overline{A}_i, i \in I)$$

This proves the statement. □

**Example 8.3.** Let $\overline{V}$ be a ${}_*{}^*D$-vector space and ${}_*$-row

$$\overline{\overline{A}} = \begin{pmatrix} \overline{A}_1 & \dots & \overline{A}_n \end{pmatrix} = (\overline{A}_i, i \in I)$$

be set of vectors. To answer the question of whether vector $\overline{b} \in \mathrm{span}(\overline{A}_i, i \in I)$ we write linear equation

$$\overline{b} = \overline{A}_*{}^*x \tag{8.1}$$

where

$$x=\begin{pmatrix}{}^1x\\ ...\\ {}^nx\end{pmatrix}$$

is ${}^*$-row of unknown coefficients of expansion. $\overline{b}\in\mathrm{span}(\overline{A}_i,i\in I)$ if equation (8.1) has a solution. Suppose $\overline{\overline{f}}=(\overline{f}_j,j\in J)$ is a ${}_*{}^*D$-basis. Then vectors $\overline{b}$, $\overline{A}_i$ have expansion

$$\overline{b}=b_*{}^*\overline{f} \tag{8.2}$$

$$\overline{A}_i=\overline{f}_*{}^*A_i \tag{8.3}$$

If we substitute (8.2) and (8.3) into (8.1) we get

$$\overline{f}_*{}^*b=\overline{f}_*{}^*A_*{}^*x \tag{8.4}$$

Applying theorem 6.3 to (8.4) we get **system of ${}_*{}^*D$-linear equations**

$$A_*{}^*x=b \tag{8.5}$$

We can write system of ${}_*{}^*D$-linear equations (8.5) in one of the next forms

$$\begin{pmatrix}{}^1A_1 & ... & {}^1A_n\\ ... & ... & ...\\ {}^mA_1 & ... & {}^mA_n\end{pmatrix}{}_*{}^*\begin{pmatrix}{}^1x\\ ...\\ {}^nx\end{pmatrix}=\begin{pmatrix}{}^1b\\ ...\\ {}^mb\end{pmatrix}$$

$${}^jA_i\ {}^ix={}^jb \tag{8.6}$$

$$\begin{matrix}{}^1A_1\ {}^1x & +... & +{}^1A_n\ {}^nx & ={}^1b\\ ... & ... & ... & ...\\ {}^1A_m\ {}^1x & +... & +{}^mA_n\ {}^nx & ={}^mb\end{matrix}$$

□

**Example 8.4.** Let $\overline{V}$ be a $D_*{}^*$-vector space and ${}^*$-row

$$\overline{\overline{A}}=\begin{pmatrix}{}^1\overline{A}\\ ...\\ {}^m\overline{A}\end{pmatrix}=({}^j\overline{A},j\in J)$$

be set of vectors. To answer the question of whether vector $\overline{b}\in\mathrm{span}({}^j\overline{A},j\in J)$ we write linear equation

$$\overline{b}=x_*{}^*\overline{A} \tag{8.7}$$

where

$$x=\begin{pmatrix}x_1 & ... & x_m\end{pmatrix}$$

is ${}_*$-row of unknown coefficients of expansion. $\overline{b}\in\mathrm{span}({}^j\overline{A},j\in J)$ if equation (8.7) has a solution. Suppose $\overline{\overline{f}}=({}^i\overline{f},i\in I)$ is a $D_*{}^*$-basis. Then vectors $\overline{b}$, ${}^j\overline{A}$ have expansion

$$\overline{b}=b_*{}^*\overline{f} \tag{8.8}$$

$${}^j\overline{A}={}^jA_*{}^*\overline{f} \tag{8.9}$$

If we substitute (8.8) and (8.9) into (8.7) we get

$$b_*{}^*\overline{f} = x_*{}^*A_*{}^*\overline{f} \tag{8.10}$$

Applying theorem 6.3 to (8.10) we get **system of $D_*{}^*$-linear equations**[7]

$$x_*{}^*A = b \tag{8.11}$$

We can write system of $D_*{}^*$-linear equations (8.11) in one of the next forms

$$\begin{pmatrix} x_1 & ... & x_m \end{pmatrix} {}_*{}^* \begin{pmatrix} {}^1A_1 & ... & {}^1A_n \\ ... & ... & ... \\ {}^mA_1 & ... & {}^mA_n \end{pmatrix} = \begin{pmatrix} b_1 & ... & b_n \end{pmatrix}$$

$$x_j\ {}^jA_i = b_i \tag{8.12}$$

$$\begin{matrix} x_1\ {}^1A_1 & ... & x_1\ {}^1A_n \\ + & ... & + \\ ... & ... & ... \\ + & ... & + \\ x_m\ {}^mA_1 & ... & x_m\ {}^mA_n \\ = & ... & = \\ b_1 & ... & b_n \end{matrix}$$

□

**Definition 8.5.** If $n\times n$ matrix $A$ has ${}_*{}^*$-inverse matrix we call such matrix **${}_*{}^*$-nonsingular matrix**. Otherwise, we call such matrix **${}_*{}^*$-singular matrix**. □

**Definition 8.6.** Suppose $A$ is ${}_*{}^*$-nonsingular matrix. We call appropriate system of ${}_*{}^*D$-linear equations

$$A_*{}^*x = b \tag{8.13}$$

**nonsingular system of ${}_*{}^*D$-linear equations**. □

**Theorem 8.7.** *Solution of nonsingular system of ${}_*{}^*D$-linear equations (8.13) is determined uniquely and can be presented in either form*[8]

$$x = A^{-1_*{}^*}{}_*{}^*b \tag{8.14}$$

$$x = \overline{H}\det(A, {}_*{}^*)\, {}_*{}^*b \tag{8.15}$$

*Proof.* Multiplying both sides of equation (8.13) from left by $A^{-1_*{}^*}$ we get (8.14). Using definition [6]-(3.12) we get (8.15). Since theorem [6]-(2.16) the solution is unique. □

---

[7] Reading system of ${}_*{}^*D$-linear equations (8.5) from bottom up and from left to right we get system of $D^*{}_*$-linear equations (8.11).

[8] We can see a solution of system (8.13) in theorem [4]-1.6.1. I repeat this statement because I slightly changed the notation.

## 9. Rank of Matrix

**Definition 9.1.** Matrix[9] ${}^SA_T$ is a minor of an order $k$. □

**Definition 9.2.** If minor ${}^SA_T$ is ${}_*{}^*$-nonsingular matrix then we say that ${}_*{}^*$-rank of matrix $A$ is not less then $k$. **${}_*{}^*$-rank of matrix** $A$

$$\mathrm{rank}_{{}_*{}^*} A$$

is the maximal value of $k$. We call an appropriate minor the **${}_*{}^*$-major minor**. □

**Theorem 9.3.** *Let matrix $A$ be ${}_*{}^*$-singular matrix and minor ${}^SA_T$ be major minor. Then*

$$
{}^p\det\left({}^{S\cup\{p\}}A_{T\cup\{r\}}, {}_*{}^*\right)_r = 0 \tag{9.1}
$$

*Proof.* To understand why minor ${}^{S\cup\{p\}}A_{T\cup\{r\}}$ does not have ${}_*{}^*$-inverse matrix,[10] we assume that it has the ${}_*{}^*$-inverse matrix and write down the respective system [6]-(3.2), [6]-(3.3). We assume $i = r$, $j = p$ and will try to solve this system. Assume

$$
B = {}^{S\cup\{p\}}A_{T\cup\{r\}} \tag{9.2}
$$

Than we get system

$$
{}^SB_T{}_*{}^*{}^TB^{-1_*{}^*}{}_p + {}^SB_r\ {}^rB^{-1_*{}^*}{}_p = 0 \tag{9.3}
$$

$$
{}^pB_T{}_*{}^*{}^TB^{-1_*{}^*}{}_p + {}^pB_r\ {}^rB^{-1_*{}^*}{}_p = 1 \tag{9.4}
$$

We multiply (9.3) by $({}^SB_T)^{-1_*{}^*}$

$$
{}^TB^{-1_*{}^*}{}_p + ({}^SB_T)^{-1_*{}^*}{}_*{}^*{}^SB_r\ {}^rB^{-1_*{}^*}{}_p = 0 \tag{9.5}
$$

Now we can substitute (9.5) into (9.4)

$$
-{}^pB_T{}_*{}^*({}^SB_T)^{-1_*{}^*}{}_*{}^*{}^SB_r\ {}^rB^{-1_*{}^*}{}_p + {}^pB_r\ {}^rB^{-1_*{}^*}{}_p = 1 \tag{9.6}
$$

From (9.6) it follows that

$$
({}^pB_r - {}^pB_T{}_*{}^*({}^SB_T)^{-1_*{}^*}{}_*{}^*{}^SB_r)\ {}^rB^{-1_*{}^*}{}_p = 1 \tag{9.7}
$$

Expression in brackets is quasideterminant ${}^p\det(B, {}_*{}^*)_r$. Substituting this expression into (9.7), we get

$$
{}^p\det(B, {}_*{}^*)_r\ {}^rB^{-1_*{}^*}{}_p = 1 \tag{9.8}
$$

Thus we proved that quasideterminant ${}^p\det(B, {}_*{}^*)_r$ is defined and its equation to 0 is necessary and sufficient condition that the matrix $B$ is singular. Since (9.2) the statement of theorem is proved. □

[9] In this section, we will make the following assumption.

- $i \in M$, $|M| = m$, $j \in N$, $|N| = n$.
- $A = ({}^iA_j)$ is an arbitrary matrix.
- $k$, $s \in S \supseteq M$, $l$, $t \in T \supseteq N$, $k = |S| = |T|$.
- $p \in M \setminus S$, $r \in N \setminus T$.

[10] It is natural to expect relationship between ${}_*{}^*$-singularity of the matrix and its ${}_*{}^*$-quasideterminant similar to relationship which is known in commutative case. However ${}_*{}^*$-quasideterminant is defined not always. For instance, it is not defined when ${}_*{}^*$-inverse matrix has too much elements equal 0. As it follows from this theorem, the ${}_*{}^*$-quasideterminant is undefined also in case when ${}_*{}^*$-rank of the matrix is less then $n - 1$.

**Theorem 9.4.** *Suppose $A$ is a matrix,*

$$\mathrm{rank}_{*}{}^{*} A = k < m$$

*and ${}^SA_T$ is ${}_*{}^*$-major minor. Then ${}_*$-row ${}^pA$ is a $D_*{}^*$-linear composition of ${}_*$-rows ${}^SA$.*

$$ {}^{M\setminus S}A = R_*{}^*{}^SA \tag{9.9}$$

$$ {}^pA = {}^pR_*{}^*{}^SA \tag{9.10}$$

$$ {}^pA_b = {}^pR_s\ {}^sA_b \tag{9.11}$$

*Proof.* If a number of ${}^*$-rows is $k$ then assuming that ${}_*$-row ${}^pA$ is a $D_*{}^*$-linear combination (9.10) of ${}_*$-rows ${}_sA$ with coefficients ${}^pR_s$ we get system (9.11). According to theorem 8.7 this system has a unique solution[11] and it is nontrivial because all ${}_*{}^*$-quasideterminants are different from 0.

It remains to prove this statement in case when a number of ${}^*$-rows is more then $k$. I get ${}_*$-row ${}^pA$ and ${}^*$-row $A_r$. According to assumption minor ${}^{S\cup\{p\}}A_{T\cup\{r\}}$ is a ${}_*{}^*$-singular matrix and its ${}_*{}^*$-quasideterminant

$$ {}^p\det\left({}^{S\cup\{p\}}A_{T\cup\{r\}},{}_*{}^*\right)_r = 0 \tag{9.12}$$

According to [6]-(3.14) the equation (9.12) has form

$$ {}^pA_r - {}^pA_T{}_*{}^*(({}^SA_T)^{-1_*{}^*})_*{}^*{}^SA_r = 0$$

Matrix

$$ {}^pR = {}^pA_T{}_*{}^*(({}^SA_T)^{-1_*{}^*}) \tag{9.13}$$

do not depend on $r$, Therefore, for any $r \in N\setminus T$

$$ {}^pA_r = {}^pR_*{}^*{}^SA_r \tag{9.14}$$

From equation

$$(({}^SA_T)^{-1_*{}^*})_*{}^*{}^SA_l = {}^T\delta_l$$

it follows that

$$ {}^pA_l = {}^pA_T{}_*{}^*{}^T\delta_l = {}^pA_T{}_*{}^*(({}^SA_T)^{-1_*{}^*})_*{}^*{}^SA_l \tag{9.15}$$

Substituting (9.13) into (9.15) we get

$$ {}^pA_l = {}^pR_*{}^*{}^SA_l \tag{9.16}$$

(9.14) and (9.16) finish the proof. □

**Corollary 9.5.** *Suppose $A$ is a matrix, $\mathrm{rank}_{*}{}^{*} A = k < m$. Then ${}_*$-rows of the matrix are $D_*{}^*$-linearly dependent.*

$$\lambda_*{}^*A = 0 \tag{9.17}$$

*Proof.* Suppose ${}_*$-row ${}^pA$ is a $D_*{}^*$-linear composition (9.10). We assume $\lambda_p = -1$, $\lambda_s = {}^pR_s$ and the rest $\lambda_c = 0$. □

**Theorem 9.6.** *Let $({}^i\overline{A}, i \in M, |M| = m)$ be set of $D_*{}^*$-linear independent vectors. Then ${}_*{}^*$-rank of their coordinate matrix equal $m$.*

[11]We assume that unknown variables are $x_s = {}^pR_s$

*Proof.* According to model built in example 6.7, the coordinate matrix of set of vectors $({}^i\overline{A})$ relative basis $\overline{\overline{e}}$ consists from ${}_*$-rows which are coordinate matrices of vectors ${}^i\overline{A}$ relative the basis $\overline{\overline{e}}$. Therefore ${}_*{}^*$-rank of this matrix cannot be more then $m$.

Let ${}_*{}^*$-rank of the coordinate matrix be less then $m$. According to corollary 9.5 ${}_*$-rows of matrix are $D_*{}^*$-linear dependent. Let us multiply both parts of equation (9.17) over ${}^*$-row $\overline{e}$. Suppose $c = \lambda_*{}^*A$. Then we get that $D_*{}^*$-linear composition

$$c_*{}^*\overline{e} = 0$$

of vectors of basis equal 0. This contradicts to statement that vectors $\overline{e}$ form basis. We proved statement of theorem. □

**Theorem 9.7.** *Suppose $A$ is a matrix,*

$$\operatorname{rank}_{{}_*{}^*} A = k < n$$

*and ${}^SA_T$ is ${}_*{}^*$-major minor. Then ${}^*$-row $A_r$ is a ${}_*{}^*D$-linear composition of ${}^*$-rows $A_t$*

$$A_{N\setminus T} = A_{T*}{}^*R \tag{9.18}$$

$$A_r = A_{T*}{}^*R_r \tag{9.19}$$

$${}^aA_r = {}^aA_t\ {}^tR_r \tag{9.20}$$

*Proof.* If a number of ${}_*$-rows is $k$ then assuming that ${}^*$-row $A_r$ is a ${}_*{}^*D$-linear combination (9.19) of ${}^*$-rows $A_t$ with coefficients ${}^tR_r$ we get system (9.20). According to theorem 8.7 this system has a unique solution[12] and it is nontrivial because all ${}_*{}^*$-quasideterminants are different from 0.

It remains to prove this statement in case when a number of ${}_*$-rows is more then $k$. I get ${}^*$-row $A_r$ and ${}_*$-row ${}^pA$. According to assumption minor ${}^{S\cup\{p\}}A_{T\cup\{r\}}$ is a ${}_*{}^*$-singular matrix and its $\overline{RC}$-quasideterminant

$${}^p\det\left({}^{S\cup\{p\}}A_{T\cup\{r\}}, {}_*{}^*\right)_r = 0 \tag{9.21}$$

According to [6]-(3.14) (9.21) has form

$${}^pA_r - {}^pA_{T*}{}^*(({}^SA_T)^{-1_*{}^*})_*{}^{*S}A_r = 0$$

Matrix

$$R_r = (({}^SA_T)^{-1_*{}^*})_*{}^{*S}A_r \tag{9.22}$$

do not depend on $p$, Therefore, for any $p \in M\setminus S$

$${}^pA_r = {}^pA_{T*}{}^*R_r \tag{9.23}$$

From equation

$${}^kA_{T*}{}^*(({}^SA_T)^{-1_*{}^*})_s = {}^k\delta_s$$

it follows that

$${}^kA_r = {}^k\delta_{S*}{}^{*S}A_r = {}^kA_{T*}{}^*(({}^SA_T)^{-1_*{}^*})^s{}_*{}^{*S}A_r \tag{9.24}$$

Substituting (9.22) into (9.24) we get

$${}^kA_r = {}^kA_{T*}{}^*R_r \tag{9.25}$$

(9.23) and (9.25) finish the proof. □

[12] We assume that unknown variables are ${}_tx = {}^tR_r$

**Corollary 9.8.** *Suppose $A$ is a matrix,* $\mathrm{rank}_{*}{}^{*} A = k < m$*. Then ${}^{*}$-rows of the matrix are ${}_{*}{}^{*}D$-linearly dependent.*

$$A_{*}{}^{*}\lambda = 0$$

*Proof.* Suppose ${}^{*}$-row $A_r$ is a right linear composition (9.19). We assume ${}^{r}\lambda = -1$, ${}^{t}\lambda{=}{}^{t}R_r$ and the rest ${}^{c}\lambda = 0$. □

Base on theorem [6]-3.8 we can write similar statements for ${}^{*}{}_{*}$-rank of matrix.

**Theorem 9.9.** *Suppose $A$ is a matrix,*

$$\mathrm{rank}_{*}{}_{*}\ A = k < m$$

*and ${}_{T}A^{S}$ is ${}^{*}{}_{*}$-major minor. Then ${}_{*}$-row $A^{p}$ is a ${}^{*}{}_{*}D$-linear composition of ${}_{*}$-rows $A^{s}$.*

$$A^{M\setminus S} = A^{S*}{}_{*}R \tag{9.26}$$

$$A^{p} = A^{S*}{}_{*}R^{p} \tag{9.27}$$

$$_{b}A^{p} = {}_{b}A^{s}\ {}_{s}R^{p} \tag{9.28}$$

**Corollary 9.10.** *Suppose $A$ is a matrix,* $\mathrm{rank}_{*}{}_{*}\ A = k < m$*. Then ${}_{*}$-rows of matrix are ${}^{*}{}_{*}D$-linearly dependent.*

$$A^{*}{}_{*}\lambda = 0$$

**Theorem 9.11.** *Suppose $A$ is a matrix,*

$$\mathrm{rank}_{*}{}_{*}\ A = k < n$$

*and ${}_{T}A^{S}$ is ${}^{*}{}_{*}$-major minor. Then ${}^{*}$-row ${}_{r}A$ is a $D^{*}{}_{*}$-linear composition of ${}^{*}$-rows ${}_{t}A$*

$$_{N\setminus T}A = R^{*}{}_{*T}A \tag{9.29}$$

$$_{r}A = {}_{r}R^{*}{}_{*T}A \tag{9.30}$$

$$_{r}A^{a} = {}_{r}R^{t}\ {}_{t}A^{a} \tag{9.31}$$

**Corollary 9.12.** *Suppose $A$ is a matrix,* $\mathrm{rank}_{*}{}_{*}\ a = k < m$*. Then ${}^{*}$-rows of matrix are $D^{*}{}_{*}$-linearly dependent.*

$$\lambda^{*}{}_{*}A = 0$$

## 10. System of ${}_{*}{}^{*}D$-Linear Equations

**Definition 10.1.** Suppose[9] $A$ is a matrix of system of $D_{*}{}^{*}$-linear equations (8.12). We call matrix

$$\begin{pmatrix} {}^{j}A_i \\ b_i \end{pmatrix} = \begin{pmatrix} {}^{1}A_1 & ... & {}^{1}A_n \\ ... & ... & ... \\ {}^{m}A_1 & ... & {}^{m}A_n \\ b_1 & ... & b_n \end{pmatrix} \tag{10.1}$$

an **extended matrix** of this system. □

**Definition 10.2.** Suppose[9] $A$ is a matrix of system of ${}_*{}^*D$-linear equations (8.6). We call matrix

$$\begin{pmatrix} {}_aA^b & {}_ab \end{pmatrix} = \begin{pmatrix} {}_1A^1 & ... & {}_1A^n & {}_1b \\ ... & ... & ... & ... \\ {}_mA^1 & ... & {}_mA^n & {}_mb \end{pmatrix} \tag{10.2}$$

an **extended matrix** of this system. □

**Theorem 10.3.** *System of ${}_*{}^*D$-linear equations (8.6) has a solution iff*

$$\mathrm{rank}_{{}_*{}^*}({}^jA_i) = \mathrm{rank}_{{}_*{}^*}\begin{pmatrix} {}^jA_i & {}^jb \end{pmatrix} \tag{10.3}$$

*Proof.* Let ${}^SA_T$ be ${}_*{}^*$-major minor of matrix $A$.

Let a system of ${}_*{}^*D$-linear equations (8.6) have solution ${}^ix = {}^id$. Then

$$A_*{}^*d = b \tag{10.4}$$

Equation (10.4) can be rewritten in form

$$A_T{}_*{}^*{}^Td + A_{N\setminus T}{}_*{}^*{}^{N\setminus T}d = b \tag{10.5}$$

Substituting (9.18) into (10.5) we get

$$A_T{}_*{}^*{}^Td + A_T{}_*{}^*R_*{}^*{}^{N\setminus T}d = b \tag{10.6}$$

From (10.6) it follows that ${}_*$-row $b$ is a ${}_*{}^*D$-linear combination of ${}_*$-rows $A_T$

$$A_T{}_*{}^*({}^Td + R_*{}^*{}^{N\setminus T}d) = b$$

This holds equation (10.3).

It remains to prove that an existence of solution of system of ${}_*{}^*D$-linear equations (8.6) follows from (10.3). Holding (10.3) means that ${}^SA_T$ is ${}_*{}^*$-major minor of extended matrix as well. From theorem 9.7 it follows that ${}_*$-row $b$ is a ${}_*{}^*D$-linear composition of ${}_*$-rows $A_T$

$$b = A_T{}_*{}^*{}^TR$$

Assigning ${}^rR = 0$ we get

$$b = A_*{}^*R$$

Therefore, we found at least one solution of system of ${}_*{}^*D$-linear equations (8.6). □

**Theorem 10.4.** *Suppose (8.6) is a system of ${}_*{}^*D$-linear equations satisfying (10.3). If $\mathrm{rank}_{{}_*{}^*} A = k \le m$ then solution of the system depends on arbitrary values of $m-k$ variables not included into ${}_*{}^*$-major minor.*

*Proof.* Let ${}^SA_T$ be ${}_*{}^*$-major minor of matrix $a$. Suppose

$${}^pA_*{}^*x = {}^pb \tag{10.7}$$

is an equation with number $p$. Applying theorem 9.4 to extended matrix (10.1) we get

$${}^pA = {}^pR_*{}^*{}^SA \tag{10.8}$$

$${}^pb = {}^pR_*{}^*{}^Sb \tag{10.9}$$

Substituting (10.8) and (10.9) into (10.7) we get

$${}^pR_*{}^*{}^SA_*{}^*x = {}^pR_*{}^*{}^Sb \tag{10.10}$$

(10.10) means that we can exclude equation (10.7) from system (8.6) and the new system is equivalent to the old one. Therefore, a number of equations can be reduced to $k$.

At this point, we have two choices. If the number of variables is also $k$ then according to theorem 8.7 the system has unique solution (8.15). If the number of variables $m > k$ then we can move $m-k$ variables that are not included into ${}_*{}^*$-major minor into right side. Giving arbitrary values to these variables, we determine value of the right side and for this value we get a unique solution according to theorem 8.7. □

**Corollary 10.5.** *System of ${}_*{}^*D$-linear equations (8.6) has a unique solution iff its matrix is nonsingular.* □

**Theorem 10.6.** *Solutions of a homogenous system of ${}_*{}^*D$-linear equations*

$$A{}_*{}^*x = 0 \tag{10.11}$$

*form a ${}_*{}^*D$-vector space.*

*Proof.* Let $\overline{X}$ be set of solutions of system of ${}_*{}^*D$-linear equations (10.11). Suppose $x = ({}^ax) \in \overline{X}$ and $y = ({}^ay) \in \overline{X}$. Then

$$x^a\ {}_aa^b = 0$$
$$y^a\ {}_aa^b = 0$$

Therefore

$${}^iA_j\ ({}^jx + {}^jy) = {}^iA_j\ {}^jx + {}^iA_j\ {}^jy = 0$$
$$x + y = ({}^jx + {}^jy) \in \overline{X}$$

The same way we see

$${}^iA_j\ ({}^jxb) = ({}^iA_j\ {}^jx)b = 0$$
$$xb = ({}^jxb) \in \overline{X}$$

According to definition 4.4 $\overline{X}$ is a ${}_*{}^*D$-vector space. □

## 11. Nonsingular Matrix

Suppose $A$ is $n \times n$ matrix. Corollaries 9.5 and 9.8 tell us that if $\operatorname{rank}_{{}_*{}^*} A < n$ then ${}_*$-rows are $D_*{}^*$-linearly dependent and ${}^*$-rows are ${}_*{}^*D$-linearly dependent.[13]

**Theorem 11.1.** *Let $A$ be $n \times n$ matrix and ${}^*$-row $A_r$ be a ${}_*{}^*D$-linear combination of other ${}^*$-rows. Then $\operatorname{rank}_{{}_*{}^*} A < n$.*

*Proof.* The statement that ${}^*$-row $A_r$ is a ${}_*{}^*D$-linear combination of other ${}^*$-rows means that system of ${}_*{}^*D$-linear equations

$$A_r = A_{[r]}{}_*{}^*\lambda$$

has at least one solution. According theorem 10.3

$$\operatorname{rank}_{{}_*{}^*} A = \operatorname{rank}_{{}_*{}^*} A_{[r]}$$

Since a number of ${}_*$-rows is less then $n$ we get $\operatorname{rank}_{{}_*{}^*} A[r] < n$. □

**Theorem 11.2.** *Let $A$ be $n \times n$ matrix and ${}_*$-row ${}^pA$ be a $D_*{}^*$-linear combination of other ${}_*$-rows. Then $\operatorname{rank}_{{}_*{}^*} A < n$.*

---

[13]This statement is similar to proposition [5]-1.2.5.

*Proof.* Proof of statement is similar to proof of theorem 11.1 □

**Theorem 11.3.** *Suppose $A$ and $B$ are $n\times n$ matrices and*

$$C = A_*{}^*B \tag{11.1}$$

*$C$ is ${}_*{}^*$-singular matrix iff either matrix $A$ or matrix $B$ is ${}_*{}^*$-singular matrix.*

*Proof.* Suppose matrix $B$ is ${}_*{}^*$-singular. According theorem 9.7 ${}^*$-rows of matrix $B$ are ${}_*{}^*D$-linearly dependent. Therefore

$$0 = B_*{}^*\lambda \tag{11.2}$$

where $\lambda \neq 0$. From (11.1) and (11.2) it follows that

$$C_*{}^*\lambda = A_*{}^*B_*{}^*\lambda = 0$$

According theorem 11.1 matrix $C$ is ${}_*{}^*$-singular.

Suppose matrix $B$ is not ${}_*{}^*$-singular, but matrix $A$ is ${}_*{}^*$-singular According to theorem 9.4 ${}^*$-rows of matrix $A$ are ${}_*{}^*D$-linearly dependent. Therefore

$$0 = A_*{}^*\mu \tag{11.3}$$

where $\mu \neq 0$. According to theorem 8.7 the system

$$B_*{}^*\lambda = \mu$$

has only solution where $\lambda \neq = 0$. Therefore

$$C_*{}^*\lambda = A_*{}^*B_*{}^*\lambda = A_*{}^*\mu = 0$$

According to theorem 11.1 matrix $C$ is ${}_*{}^*$-singular.

Suppose matrix $C$ is ${}_*{}^*$-singular matrix. According to the theorem 9.4 ${}_*$-rows of matrix $C$ are ${}_*{}^*D$-linearly dependent. Therefore

$$0 = C_*{}^*\lambda \tag{11.4}$$

where $\lambda \neq 0$. From (11.1) and (11.4) it follows that

$$0 = A_*{}^*B_*{}^*\lambda$$

If

$$0 = B_*{}^*\lambda$$

satisfied then matrix $B$ is ${}_*{}^*$-singular. Suppose that matrix $B$ is not ${}_*{}^*$-singular. Let us introduce

$$\mu = B_*{}^*\lambda$$

where $\mu \neq 0$. Then

$$0 = A_*{}^*\mu \tag{11.5}$$

From (11.5) it follows that matrix $A$ is ${}_*{}^*$-singular. □

Basing theorem [6]-3.8 we can write similar statements for $D^*{}_*$-linear combination of ${}_*$-rows or ${}^*{}_*D$-linear combination of ${}^*$-rows and ${}^*{}_*$-quasideterminant.

**Theorem 11.4.** *Let $A$ be $n\times n$ matrix and ${}^*$-row ${}_rA$ be a ${}^*{}_*D$-linear combination of other ${}^*$-rows. Then* $\mathrm{rank}_{{}^*{}_*} A < n$.

**Theorem 11.5.** *Let $A$ be $n\times n$ matrix and ${}_*$-row $A^p$ be a $D^*{}_*$-linear combination of other ${}_*$-rows. Then* $\mathrm{rank}_{{}^*{}_*} A < n$.

**Theorem 11.6.** *Suppose $A$ and $B$ are $n\times n$ matrices and $C = A^*{}_*B$. $C$ is ${}^*{}_*$-singular matrix iff either matrix $A$ or matrix $B$ is ${}^*{}_*$-singular matrix.*

**Definition 11.7.** ${}_*{}^*$**-matrix group** $GL_n{}_*{}^*{}^D$ is a group of ${}_*{}^*$-nonsingular matrices where we define ${}_*{}^*$-product of matrices [6]-(2.1) and ${}_*{}^*$-inverse matrix $A^{-1_*{}^*}$. □

**Definition 11.8.** ${}^*{}_*$**-matrix group** $GL_n{}^*{}_*{}^D$ is a group of ${}^*{}_*$-nonsingular matrices where we define ${}^*{}_*$-product of matrices [6]-(2.2) and ${}^*{}_*$-inverse matrix $A^{-1^*{}_*}$. □

**Theorem 11.9.**

$$GL_n{}_*{}^*{}^D \neq GL_n{}^*{}_*{}^D$$

*Remark* 11.10. From theorem [6]-(3.10) it follows that there are matrices which are ${}^*{}_*$-nonsingular and ${}_*{}^*$-nonsingular. Theorem 11.9 implies that sets of ${}^*{}_*$-nonsingular matrices and ${}_*{}^*$-nonsingular matrices are not identical. For instance, there exists such ${}_*{}^*$-nonsingular matrix which is ${}^*{}_*$-singular matrix. □

*Proof.* It is enough to prove this statement for $n = 2$. Assume every ${}^*{}_*$-singular matrix

$$A = \begin{pmatrix} {}^1A_1 & {}^2A_1 \\ {}^1A_2 & {}^2A_2 \end{pmatrix} \tag{11.6}$$

is ${}_*{}^*$-singular matrix. It follows from theorem 9.4 and theorem 9.7 that ${}_*{}^*$-singular matrix satisfies to condition

$$ {}^1A_2 = b\ {}^1A_1 \tag{11.7}$$

$$ {}^2A_2 = b\ {}^2A_1 \tag{11.8}$$

$$ {}^2A_1 = {}^1A_1 c \tag{11.9}$$

$$ {}^2A_2 = {}^1A_2 c \tag{11.10}$$

If we substitute (11.9) into (11.8) we get

$$ {}^2A_2 = b\ {}^1A_1\ c$$

$b$ and $c$ are arbitrary elements of division ring $D$ and ${}_*{}^*$-singular matrix matrix (11.6) has form ($d = {}^1A_1$)

$$A = \begin{pmatrix} d & dc \\ bd & bdc \end{pmatrix} \tag{11.11}$$

The similar way we can show that ${}^*{}_*$-singular matrix has form

$$A = \begin{pmatrix} d & c'd \\ db' & c'db' \end{pmatrix} \tag{11.12}$$

From assumption it follows that (11.12) and (11.11) represent the same matrix. Comparing (11.12) and (11.11) we get that for every $d, c \in D$ exists such $c' \in D$ which does not depend on $d$ and satisfies equation

$$dc = c'd$$

This contradicts the fact that $D$ is division ring. □

**Example 11.11.** Since we get division ring of quaternions we assume $b = 1 + k$, $c = j$, $d = k$. Then we get

$$A = \begin{pmatrix} k & kj \\ (1+k)k & (1+k)kj \end{pmatrix} = \begin{pmatrix} k & -i \\ k-1 & -i-j \end{pmatrix}$$

$$\begin{aligned}{}^2\det(A,{}_*{}^*)_2&={}^2A_2-{}^2A_1({}^1A_1)^{-1}\ {}^1A_2\\&=-i-j-(k-1)(k)^{-1}(-i)=-i-j-(k-1)(-k)(-i)\\&=-i-j-kki+ki=-i-j+i+j\\&=0\end{aligned}$$

$$\begin{aligned}{}_1\det(A,{}_*{}^*)^1&={}_1A^1-{}_1A^2({}_2A^2)^{-1}\ {}_2A^1\\&=k-(k-1)(-i-j)^{-1}(-i)=k-(k-1)\frac12(i+j)(-i)\\&=k+\frac12((k-1)i+(k-1)j)i=k+\frac12(ki-i+kj-j)i\\&=k+\frac12(j-i-i-j)i=k-ii\\&=k+1\end{aligned}$$

$$\begin{aligned}{}_1\det(A,{}^*{}_*)^2&={}_1A^2-{}_1A^1({}_2A^1)^{-1}\ {}_2A^2\\&=k-1-k(-i)^{-1}(-i-j)=k-1+ki(i+j)\\&=k-1+j(i+j)=k-1+ji+jj\\&=k-1-k-1\\&=-2\end{aligned}$$

$$\begin{aligned}{}_2\det(A,{}^*{}_*)^1&={}_2A^1-{}_2A^2({}_1A^2)^{-1}\ {}_1A^1\\&=(-i)-(-i-j)(k-1)^{-1}k=-i+(i+j)\frac12(-k-1)k\\&=-i-\frac12(i+j)(k+1)k=-i-\frac12(ik+i+jk+j)k\\&=-i-\frac12(-j+i+i+j)k=-i-ik\\&=-i+j\end{aligned}$$

$$\begin{aligned}{}_2\det(A,{}^*{}_*)^2&={}_2A^2-{}_2A^1({}_1A^1)^{-1}\ {}_1A^2\\&=-i-j-(-i)(k)^{-1}(k-1)=-i-j+i(-k)(k-1)\\&=-i-j+j(k-1)=-i-j+jk-j=-i-j+i-j\\&=-2j\end{aligned}$$

The system of ${}_*{}^*D$-linear equations

$$\begin{pmatrix}k&-i\\k-1&-i-j\end{pmatrix}{}_*{}^*\begin{pmatrix}{}^1x\\{}^2x\end{pmatrix}=\begin{pmatrix}{}^1b\\{}^2b\end{pmatrix}\tag{11.13}$$

has ${}_*{}^*$-singular matrix. We can write the system of ${}_*{}^*D$-linear equations (11.13) in the form

$$\left\{\begin{array}{rl}k\ {}^1x-i\ {}^2x&={}^1b\\(k-1)\ {}^1x-(i+j)\ {}^2x&={}^2b\end{array}\right.$$

The system of ${}^*{}_*D$-linear equations

$$\begin{pmatrix} k & -i \\ k-1 & -i-j \end{pmatrix} {}^*{}_* \begin{pmatrix} {}_1x & {}_2x \end{pmatrix} = \begin{pmatrix} {}_1b & {}_2b \end{pmatrix} \tag{11.14}$$

has ${}^*{}_*$-nonsingular matrix. We can write the system of ${}^*{}_*D$-linear equations (11.14) in the form

$$\left\{ \begin{matrix} k\ {}_1x & - & i\ {}_1x \\ +\ (k-1)\ {}_2x & -\ (i+j)\ {}_2x \\ = & {}_1b\ = & {}_2b \end{matrix} \right. \qquad \left\{ \begin{matrix} k\ {}_1x + (k-1)\ {}_2x = {}_1b \\ -i\ {}_1x - \ (i+j)\ {}_2x = {}_2b \end{matrix} \right.$$

The system of $D_*{}^*$-linear equations

$$\begin{pmatrix} x_1 & x_2 \end{pmatrix} {}_*{}^* \begin{pmatrix} k & -i \\ k-1 & -i-j \end{pmatrix} = \begin{pmatrix} b_1 & b_2 \end{pmatrix} \tag{11.15}$$

has ${}_*{}^*$-singular matrix. We can write the system of $D_*{}^*$-linear equations (11.15) in the form

$$\left\{ \begin{matrix} x_1k & -x_1i \\ +\ x_2(k-1) & -x_2(i+j) \\ = b_1 & = b_2 \end{matrix} \right. \qquad \left\{ \begin{matrix} x_1k & +x_2(k-1) = b_1 \\ -x_1i & -x_2(i+j) \ = b_2 \end{matrix} \right.$$

The system of $D^*{}_*$-linear equations

$$\begin{pmatrix} x^1 \\ x^2 \end{pmatrix} {}^*{}_* \begin{pmatrix} k & -i \\ k-1 & -i-j \end{pmatrix} = \begin{pmatrix} {}^1b \\ {}^2b \end{pmatrix} \tag{11.16}$$

has ${}^*{}_*$-nonsingular matrix. We can write the system of $D^*{}_*$-linear equations (11.16) in the form

$$\left\{ \begin{matrix} k\ {}^1x - & i\ {}^2x = {}^1b \\ (k-1)\ {}^1x - & (i+j)\ {}^2x = {}^2b \end{matrix} \right.$$

□

## 12. Dimension of ${}_*{}^*D$-Vector Space

**Theorem 12.1.** *Let $\overline{V}$ be a ${}_*{}^*D$-vector space. Suppose $\overline{V}$ has ${}_*{}^*D$-bases $\overline{\overline{e}} = (\overline{e}_i, i \in I)$ and $\overline{\overline{g}} = (\overline{g}_j, j \in J)$. If $|I|$ and $|J|$ are finite numbers then $|I| = |J|$.*

*Proof.* Suppose $|I| = m$ and $|J| = n$. Suppose

$$m < n \tag{12.1}$$

Because $\overline{\overline{e}}$ is a ${}_*{}^*D$-basis any vector $\overline{g}_j,\ j \in J$ has expansion

$$_b\overline{g} = {}_ba_*{}^*\overline{e}$$

Because $\overline{\overline{g}}$ is a ${}_*{}^*D$-basis

$$\lambda = 0 \tag{12.2}$$

should follow from

$$\overline{g}_*{}^*\lambda = \overline{e}_*{}^*A_*{}^*\lambda = 0$$

Because $\overline{\overline{e}}$ is a ${}_*{}^*D$-basis we get

$$A_*{}^*\lambda = 0 \tag{12.3}$$

According to (12.1) $\operatorname{rank}_{*}{}^{*} A \le m$ and system (12.3) has more variables then equations. According to theorem 10.4 $\lambda \ne 0$. This contradicts statement (12.2). Therefore, statement $m < n$ is not valid.

In the same manner we can prove that the statement $n < m$ is not valid. This completes the proof of the theorem. □

**Definition 12.2.** We call **dimension of ${}_*{}^*D$-vector space** the number of vectors in a basis □

**Theorem 12.3.** *The coordinate matrix of ${}_*{}^*D$-basis $\overline{\overline{g}}$ relative ${}_*{}^*D$-basis $\overline{\overline{e}}$ of vector space $\overline{V}$ is ${}_*{}^*$-nonsingular matrix.*

*Proof.* According to theorem 9.6 ${}_*{}^*D$-rank of the coordinate matrix of basis $\overline{\overline{g}}$ relative basis $\overline{\overline{e}}$ equal to the dimension of vector space. This proves the statement of the theorem. □

**Definition 12.4.** We call one-to-one map $\overline{A} : \overline{V} \to \overline{W}$ **${}_*{}^*D$-isomorphism of vector spaces** if this map is a ${}_*{}^*D$-linear map of vector spaces. □

**Definition 12.5.** **${}_*{}^*D$-automorphism of vector space $\overline{V}$** is ${}_*{}^*D$-isomorphism $\overline{A} : \overline{V} \to \overline{V}$. □

**Theorem 12.6.** *Suppose that $\overline{\overline{f}}$ is a ${}_*{}^*D$-basis of vector space $\overline{V}$. Then any ${}_*{}^*D$-automorphism $\overline{A}$ of vector space $\overline{V}$ has form*

$$v' = A_*{}^*v \tag{12.4}$$

*where $A$ is a ${}_*{}^*$-nonsingular matrix.*

*Proof.* (12.4) follows from theorem 7.3. Because $\overline{A}$ is an isomorphism for each vector $\overline{v}'$ exist one and only one vector $\overline{v}$ such that $\overline{v}' = \overline{v}_*{}^*\overline{A}$. Therefore, system of ${}_*{}^*D$-linear equations (12.4) has a unique solution. According to corollary 10.5 matrix $A$ is a nonsingular matrix. □

**Theorem 12.7.** *Automorphisms of ${}_*{}^*D$-vector space form a group $GL_n^{{}_*{}^*D}$.*

*Proof.* If we have two automorphisms $\overline{A}$ and $\overline{B}$ then we can write

$$v' = A_*{}^*v$$

$$v'' = B_*{}^*v' = B_*{}^*A_*{}^*v$$

Therefore, the resulting automorphism has matrix $A_*{}^*B$. □

## 13. References

[1] S. Burris, H.P. Sankappanavar, A Course in Universal Algebra, Springer-Verlag (March, 1982),
eprint http://www.math.uwaterloo.ca/ snburris/htdocs/ualg.html
(The Millennium Edition)

[2] A. G. Kurosh, Lectures on General Algebra, Chelsea Pub Co, 1965

[3] Lev V. Sabinin, Smooth Quasigroups and Loops, Kluwer Academic Publisher, 1999

[4] I. Gelfand, S. Gelfand, V. Retakh, R. Wilson, Quasideterminants,
eprint arXiv:math.QA/0208146 (2002)

[5] I.Gelfand, V.Retakh, Quasideterminants, I,
eprint arXiv:q-alg/9705026 (1997)

[6] Aleks Kleyn, Biring of Matrices,
eprint arXiv:math.OA/0612111 (2007)
[7] Paul M. Cohn, Universal Algebra, Springer, 1981

# 14. Index

# 15. Special Symbols and Notations

# Векторное пространство над телом

Александр Клейн

Аннотация. Система линейных уравнений над телом имеет свойства, похожие на свойства систем линейных уравнений над полем. Хотя некоммутативность произведения порождает новую картину, существует тесная связь между свойствами системы линейных уравнений и векторного пространства над телом.

## Содержание

## 1. Представление универсальной алгебры

**Определение 1.1.** Пусть на множестве $M$ определена структура $\Omega_2$-алгебры ([1, 7]). Мы будем называть эндоморфизм $\Omega_2$-алгебры

$$t : M \to M$$

**преобразование универсальной алгебры** $M$.[1] □

Aleks_Kleyn@MailAPS.org.
http://sites.google.com/site/AleksKleyn/.
http://arxiv.org/a/kleyn_a_1.
http://AleksKleyn.blogspot.com/.

[1]Если множество операций $\Omega_2$-алгебры пусто, то $t$ является отображением.

Мы будем обозначать $\delta$ тождественное преобразование.

**Определение 1.2.** Преобразование называется **левосторонним преобразованием** или **$T\star$-преобразованием**, если оно действует слева

$$u' = tu$$

Мы будем обозначать ${}^\star M$ множество $T\star$-преобразований множества $M$. □

**Определение 1.3.** Преобразование называется **правосторонним преобразованием** или **$\star T$-преобразованием**, если оно действует справа

$$u' = ut$$

Мы будем обозначать $M^\star$ множество $\star T$-преобразований множества $M$. □

**Определение 1.4.** Пусть на множестве ${}^\star M$ определена структура $\Omega_1$-алгебры ([1]). Пусть $A$ является $\Omega_1$-алгеброй. Мы будем называть гомоморфизм

$$f : A \to {}^\star M \tag{1.1}$$

**левосторонним** или **$T\star$-представлением $\Omega_1$-алгебры $A$ в $\Omega_2$-алгебре $M$** □

**Определение 1.5.** Пусть на множестве $M^\star$ определена структура $\Omega_1$-алгебры ([1]). Пусть $A$ является $\Omega_1$-алгеброй. Мы будем называть гомоморфизм

$$f : A \to M^\star$$

**правосторонним** или **$\star T$-представлением $\Omega_1$-алгебры $A$ в $\Omega_2$-алгебре $M$** □

Мы распространим на теорию представлений соглашение, описанное в замечании [6]-2.15. Мы можем записать принцип двойственности в следующей форме

**Теорема 1.6** (принцип двойственности)**.** *Любое утверждение, справедливое для $T\star$-представления $\Omega_1$-алгебры $A$, будет справедливо для $\star T$-представления $\Omega_1$-алгебры $A$.*

*Замечание* 1.7. Существует две формы записи преобразования $\Omega_2$-алгебры $M$. Если мы пользуемся операторной записью, то преобразование $A$ записывается в виде $Aa$ или $aA$, что соответствует $T\star$-преобразованию или $\star T$-преобразованию. Если мы пользуемся функциональной записью, то преобразование $A$ записывается в виде $A(a)$ независимо от того, это $T\star$-преобразование или $\star T$-преобразование. Эта запись согласована с принципом двойственности.

Это замечание является основой следующего соглашения. Когда мы пользуемся функциональной записью, мы не различаем $T\star$-преобразование и $\star T$-преобразование. Мы будем обозначать ${}^* M$ множество преобразований $\Omega_2$-алгебры $M$. Пусть на множестве ${}^* M$ определена структура $\Omega_1$-алгебры. Пусть $A$ является $\Omega_1$-алгеброй. Мы будем называть гомоморфизм

$$f : A \to {}^* M \tag{1.2}$$

**представлением $\Omega_1$-алгебры $A$ в $\Omega_2$-алгебре $M$.**

Соответствие между операторной записью и функциональной записью однозначно. Мы можем выбирать любую форму записи, которая удобна для изложения конкретной темы. □

Диаграмма

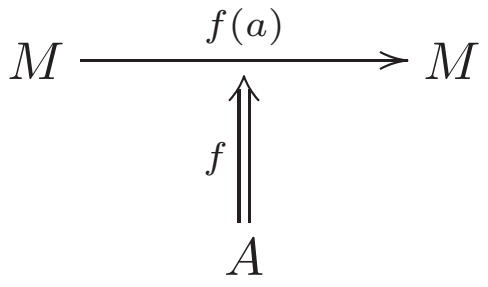

означает, что мы рассматриваем представление $\Omega_1$-алгебры $A$. Отображение $f(a)$ является образом $a \in A$.

**Определение 1.8.** Мы будем называть представление $\Omega_1$-алгебры $A$ **эффективным**, если отображение (1.2) - изоморфизм $\Omega_1$-алгебры $A$ в ${}^*M$. □

*Замечание* 1.9. Если $T\star$-представление $\Omega_1$-алгебры эффективно, мы можем отождествлять элемент $\Omega_1$-алгебры с его образом и записывать $T\star$-преобразование, порождённое элементом $a \in A$, в форме

$$v' = av$$

Если $\star T$-представление $\Omega_1$-алгебры эффективно, мы можем отождествлять элемент $\Omega_1$-алгебры с его образом и записывать $\star T$-преобразование, порождённое элементом $a \in A$, в форме

$$v' = va$$

□

**Определение 1.10.** Мы будем называть представление $\Omega_1$-алгебры **транзитивным**, если для любых $a, b \in V$ существует такое $g$, что

$$a = f(g)(b)$$

Мы будем называть представление $\Omega_1$-алгебры **однотранзитивным**, если оно транзитивно и эффективно. □

**Теорема 1.11.** *$T\star$-представление однотранзитивно тогда и только тогда, когда для любых $a, b \in M$ существует одно и только одно $g \in A$ такое, что $a = f(g)(b)$*

*Доказательство.* Следствие определений 1.8 и 1.10. □

## 2. Морфизм представлений универсальной алгебры

**Теорема 2.1.** *Пусть $A$ и $B$ - $\Omega_1$-алгебры. Представление $\Omega_1$-алгебры $B$*

$$g : B \rightarrow {}^\star M$$

*и гомоморфизм $\Omega_1$-алгебры*

$$h : A \rightarrow B \tag{2.1}$$

*определяют представление $f$ $\Omega_1$-алгебры $A$*

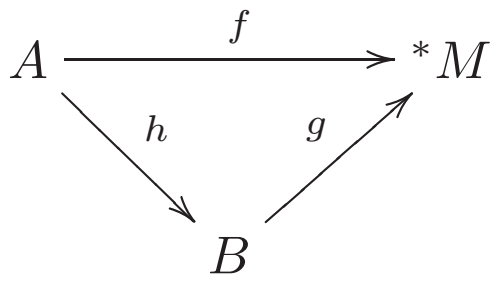

*Доказательство.* Отображение $f$ является гомоморфизмом $\Omega_1$-алгебры $A$ в $\Omega_1$-алгебру ${}^*M$, так как отображение $g$ является гомоморфизмом $\Omega_1$-алгебры $B$ в $\Omega_1$-алгебру ${}^*M$. □

Если мы изучаем представление $\Omega_1$-алгебры в $\Omega_2$-алгебрах $M$ и $N$, то нас интересуют отображения из $M$ в $N$, сохраняющие структуру представления.

**Определение 2.2.** Пусть

$$f : A \to {}^*M$$

представление $\Omega_1$-алгебры $A$ в $\Omega_2$-алгебре $M$ и

$$g : B \to {}^*N$$

представление $\Omega_1$-алгебры $B$ в $\Omega_2$-алгебре $N$. Пара отображений

$$(r : A \to B, R : M \to N) \tag{2.2}$$

таких, что

- $r$ - гомоморфизм $\Omega_1$-алгебры
- $R$ - гомоморфизм $\Omega_2$-алгебры
- $$R \circ f(a) = g(r(a)) \circ R \tag{2.3}$$

называется **морфизмом представлений из $f$ в $g$**. Мы также будем говорить, что определён **морфизм представлений $\Omega_1$-алгебры в $\Omega_2$-алгебре**. □

Для произвольного $m \in M$ равенство (2.3) имеет вид

$$R(f(a)(m)) = g(r(a))(R(m)) \tag{2.4}$$

*Замечание* 2.3. Мы можем рассматривать пару отображений $r$, $R$ как отображение

$$F : A \cup M \to B \cup N$$

такое, что

$$F(A) = B \quad F(M) = N$$

Поэтому в дальнейшем мы будем говорить, что дано отображение $(r, R)$. □

*Замечание* 2.4. Рассмотрим морфизм представлений (2.2). Мы можем обозначать элементы множества $B$, пользуясь буквой по образцу $b \in B$. Но если мы хотим показать, что $b$ является образом элемента $a \in A$, мы будем пользоваться обозначением $r(a)$. Таким образом, равенство

$$r(a) = r(a)$$

означает, что $r(a)$ (в левой части равенства) является образом $a \in A$ (в правой части равенства). Пользуясь подобными соображениями, мы будем обозначать элемент множества $N$ в виде $R(m)$. Мы будем следовать этому соглашению, изучая соотношения между гомоморфизмами $\Omega_1$-алгебр и отображениями между множествами, где определены соответствующие представления.

Мы можем интерпретировать (2.4) двумя способами

- Пусть преобразование $f(a)$ отображает $m \in M$ в $f(a)(m)$. Тогда преобразование $g(r(a))$ отображает $R(m) \in N$ в $R(f(a)(m))$.

- Мы можем представить морфизм представлений из $f$ в $g$, пользуясь диаграммой

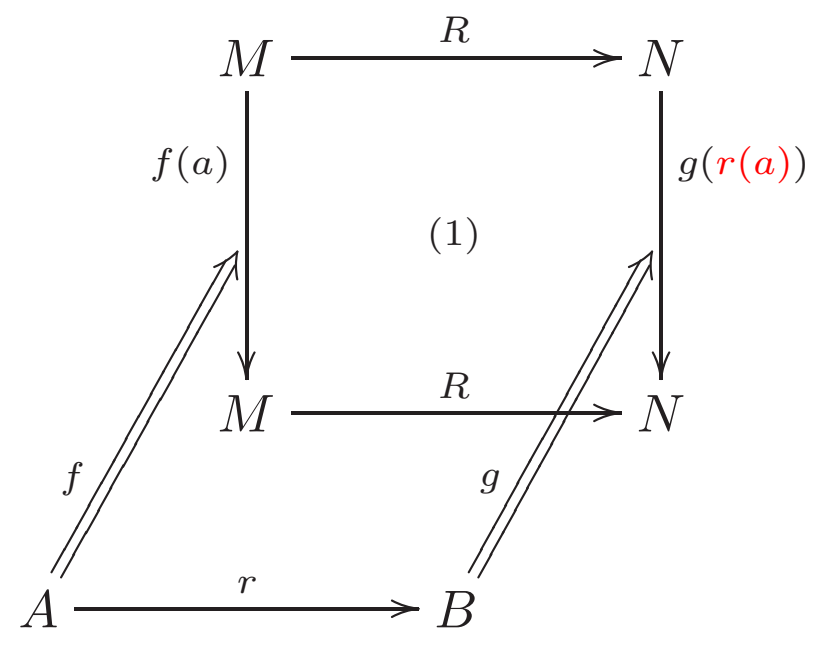

Из (2.3) следует, что диаграмма (1) коммутативна.

□

**Теорема 2.5.** *Рассмотрим представление*

$$f : A \to {}^*M$$

*$\Omega_1$-алгебры $A$ и представление*

$$g : B \to {}^*N$$

*$\Omega_1$-алгебры $B$. Морфизм*

$$h : A \longrightarrow B \qquad H : M \longrightarrow N$$

*представлений из $f$ в $g$ удовлетворяет соотношению*

$$H \circ \omega(f(a_1), ..., f(a_n)) = \omega(g(h(a_1)), ..., g(h(a_n))) \circ H \tag{2.5}$$

*для произвольной $n$-арной операции $\omega$ $\Omega_1$-алгебры.*

*Доказательство.* Так как $f$ - гомоморфизм, мы имеем

$$H \circ \omega(f(a_1), ..., f(a_n)) = H \circ f(\omega(a_1, ..., a_n)) \tag{2.6}$$

Из (2.3) и (2.6) следует

$$H \circ \omega(f(a_1), ..., f(a_n)) = g(h(\omega(a_1, ..., a_n))) \circ H \tag{2.7}$$

Так как $h$ - гомоморфизм, из (2.7) следует

$$H \circ \omega(f(a_1), ..., f(a_n)) = g(\omega(h(a_1), ..., h(a_n))) \circ H \tag{2.8}$$

Так как $g$ - гомоморфизм, из (2.8) следует (2.5). □

**Теорема 2.6.** *Пусть отображение*

$$h : A \longrightarrow B \qquad H : M \longrightarrow N$$

*является морфизмом из представления*

$$f : A \to {}^*M$$

*$\Omega_1$-алгебры $A$ в представление*

$$g : B \to {}^*N$$

*$\Omega_1$-алгебры $B$. Если представление $f$ эффективно, то отображение*

$${}^*H : {}^*M \to {}^*N$$

*определённое равенством*

$$^*H(f(a)) = g(h(a)) \tag{2.9}$$

*является гомоморфизмом* $\Omega_1$*-алгебры.*

*Доказательство.* Так как представление $f$ эффективно, то для выбранного преобразования $f(a)$ выбор элемента $a$ определён однозначно. Следовательно, преобразование $g(h(a))$ в равенстве (2.9) определено корректно.

Так как $f$ - гомоморфизм, мы имеем

$$^*H(\omega(f(a_1), ..., f(a_n))) = {}^*H(f(\omega(a_1, ..., a_n))) \tag{2.10}$$

Из (2.9) и (2.10) следует

$$^*H(\omega(f(a_1), ..., f(a_n))) = g(h(\omega(a_1, ..., a_n))) \tag{2.11}$$

Так как $h$ - гомоморфизм, из (2.11) следует

$$^*H(\omega(f(a_1), ..., f(a_n))) = g(\omega(h(a_1), ..., h(a_n))) \tag{2.12}$$

Так как $g$ - гомоморфизм,

$$^*H(\omega(f(a_1), ..., f(a_n))) = \omega(g(h(a_1)), ..., g(h(a_n))) = \omega({}^*H(f(a_1)), ..., {}^*H(f(a_n)))$$

следует из (2.12). Следовательно, отображение $^*H$ является гомоморфизмом $\Omega_1$-алгебры. □

**Теорема 2.7.** *Если представление*

$$f : A \to {}^*M$$

$\Omega_1$*-алгебры* $A$ *однотранзитивно и представление*

$$g : B \to {}^*N$$

$\Omega_1$*-алгебры* $B$ *однотранзитивно, то существует морфизм*

$$h : A \longrightarrow B \qquad H : M \longrightarrow N$$

*представлений из* $f$ *в* $g$*.*

*Доказательство.* Выберем гомоморфизм $h$. Выберем элемент $m \in M$ и элемент $n \in N$. Чтобы построить отображение $H$, рассмотрим следующую диаграмму

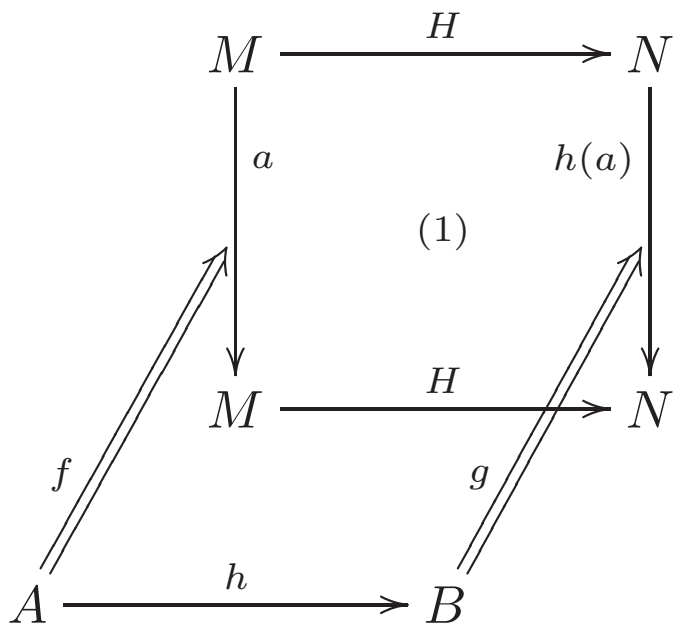

Из коммутативности диаграммы (1) следует

$$H(am) = h(a)H(m)$$

Для произвольного $m' \in M$ однозначно определён $a \in A$ такой, что $m' = am$. Следовательно, мы построили отображении $H$, которое удовлетворяет равенству (2.3). □

**Теорема 2.8.** *Если представление*

$$f : A \to {}^*M$$

*$\Omega_1$-алгебры $A$ однотранзитивно и представление*

$$g : B \to {}^*N$$

*$\Omega_1$-алгебры $B$ однотранзитивно, то для заданного гомоморфизма $\Omega_1$-алгебры*

$$h : A \longrightarrow B$$

*отображение*

$$H : M \longrightarrow N$$

*такое, что $(h, H)$ является морфизмом представлений из $f$ в $g$, единственно с точностью до выбора образа $n = H(m) \in N$ заданного элемента $m \in M$.*

*Доказательство.* Из доказательства теоремы 2.7 следует, что выбор гомоморфизма $h$ и элементов $m \in M$, $n \in N$ однозначно определяет отображение $H$. □

**Теорема 2.9.** *Если представление*

$$f : A \to {}^*M$$

*$\Omega_1$-алгебры $A$ однотранзитивно, то для любого эндоморфизма $\Omega_1$-алгебры $A$ существует эндоморфизм*

$$p : A \longrightarrow A \qquad P : M \longrightarrow M$$

*представления $f$.*

*Доказательство.* Рассмотрим следующую диаграмму

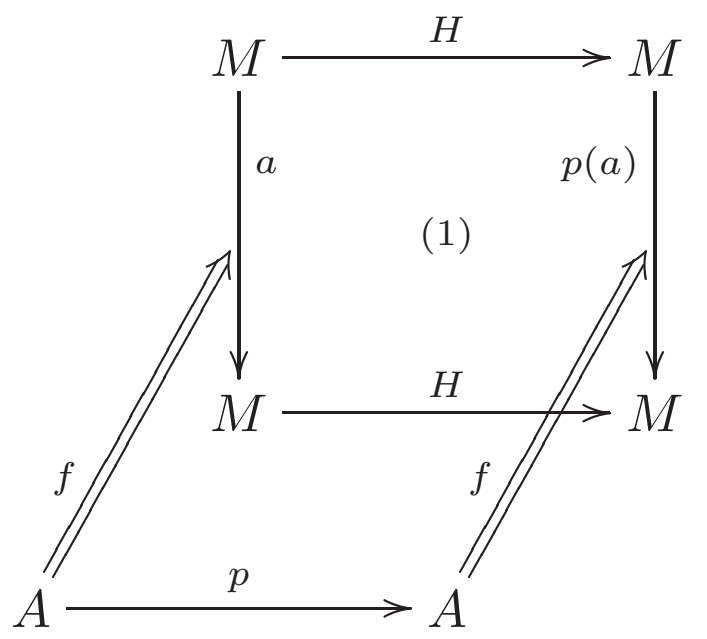

Утверждение теоремы является следствием теоремы 2.7. □

**Теорема 2.10.** *Пусть*

$$f : A \to {}^*M$$

*представление $\Omega_1$-алгебры $A$,*

$$g : B \to {}^*N$$

*представление $\Omega_1$-алгебры $B$,*

$$h : C \to {}^*L$$

*представление $\Omega_1$-алгебры $C$. Пусть определены морфизмы представлений $\Omega_1$-алгебры*

$$p : A \longrightarrow B \qquad P : M \longrightarrow N$$

$$q : B \longrightarrow C \qquad Q : N \longrightarrow L$$

*Тогда определён морфизм представлений $\Omega_1$-алгебры*

$$r : A \longrightarrow C \qquad R : M \longrightarrow L$$

*где $r = qp$, $R = QP$. Мы будем называть морфизм $(r, R)$ представлений из $f$ в $h$* **произведением морфизмов** $(p, P)$ **и** $(q, Q)$ **представлений универсальной алгебры**.

*Доказательство.* Мы можем представить утверждение теоремы, пользуясь диаграммой

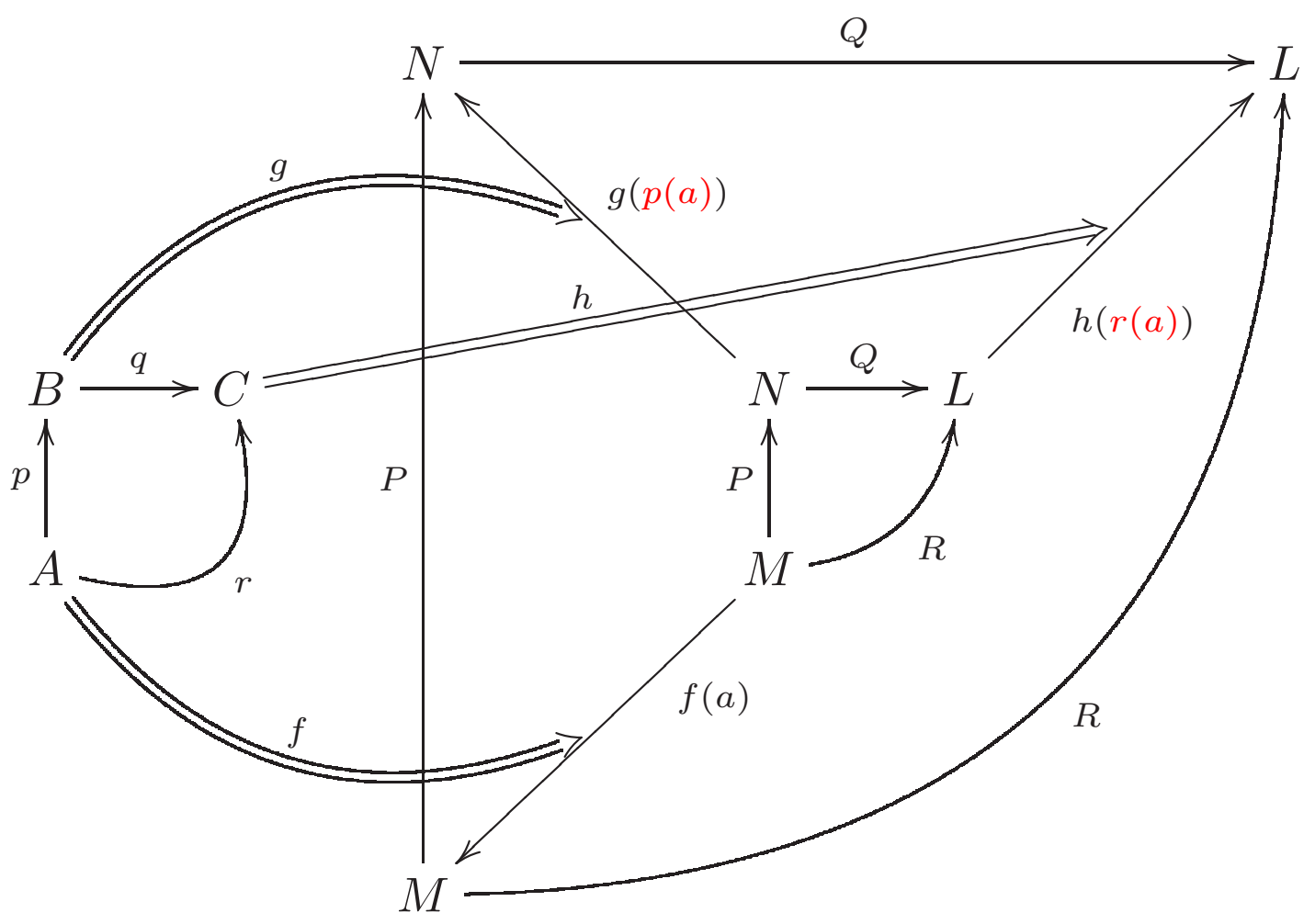

Отображение $r$ является гомоморфизмом $\Omega_1$-алгебры $A$ в $\Omega_1$-алгебру $C$. Нам надо показать, что пара отображений $(r, R)$ удовлетворяет (2.3):

$$\begin{aligned} R(f(a)m) &= QP(f(a)m) \\ &= Q(g(p(a))P(m)) \\ &= h(qp(a))QP(m)) \\ &= h(r(a))R(m) \end{aligned}$$

□

**Определение 2.11.** Допустим $\mathcal{A}$ категория $\Omega_1$-алгебр. Мы определим **категорию** $T\star\mathcal{A}$ **$T\star$-представлений универсальной алгебры из категории** $\mathcal{A}$. Объектами этой категории являются $T\star$-представлениями $\Omega_1$-алгебры. Морфизмами этой категории являются морфизмы $T\star$-представлений $\Omega_1$-алгебры. □

**Теорема 2.12.** *Эндоморфизмы представления $f$ порождают полугруппу.*

*Доказательство.* Из теоремы 2.10 следует, что произведение эндоморфизмов $(p, P)$, $(r, R)$ представления $f$ является эндоморфизмом $(pr, PR)$ представления $f$. □

**Определение 2.13.** Пусть на множестве $M$ определена эквивалентность $S$. Преобразование $f$ называется **согласованным с эквивалентностью** $S$, если из условия $m_1 \equiv m_2(\mathrm{mod} S)$ следует $f(m_1) \equiv f(m_2)(\mathrm{mod} S)$. □

**Теорема 2.14.** *Пусть на множестве $M$ определена эквивалентность $S$. Пусть на множестве ${}^*M$ определена $\Omega_1$-алгебра. Если преобразования согласованны с эквивалентностью $S$, то мы можем определить структуру $\Omega_1$-алгебры на множестве ${}^*(M/S)$.*

*Доказательство.* Пусть $h = \mathrm{nat}\ S$. Если $m_1 \equiv m_2(\mathrm{mod} S)$, то $h(m_1) = h(m_2)$. Поскольку $f \in {}^*M$ согласованно с эквивалентностью $S$, то $h(f(m_1)) = h(f(m_2))$. Это позволяет определить преобразование $F$ согласно правилу

$$F([m]) = h(f(m))$$

Пусть $\omega$ - n-арная операция $\Omega_1$-алгебры. Пусть $f_1$, ..., $f_n \in {}^\star M$ и

$$F_1([m]) = h(f_1(m)) \quad ... \quad F_n([m]) = h(f_n(m))$$

Мы определим операцию на множестве ${}^\star(M/S)$ по правилу

$$\omega(F_1, ..., F_n)[m] = h(\omega(f_1, ..., f_n)m)$$

Это определение корректно, так как $\omega(f_1, ..., f_n) \in {}^\star M$ и согласованно с эквивалентностью $S$. □

**Теорема 2.15.** *Пусть*

$$f : A \to {}^*M$$

*представление $\Omega_1$-алгебры $A$,*

$$g : B \to {}^*N$$

*представление $\Omega_1$-алгебры $B$. Пусть*

$$r : A \longrightarrow B \qquad R : M \longrightarrow N$$

*морфизм представлений из $f$ в $g$. Положим*

$$s = rr^{-1} \qquad S = RR^{-1}$$

*Тогда для отображений $r$, $R$ существуют разложения, которые можно описать диаграммой*

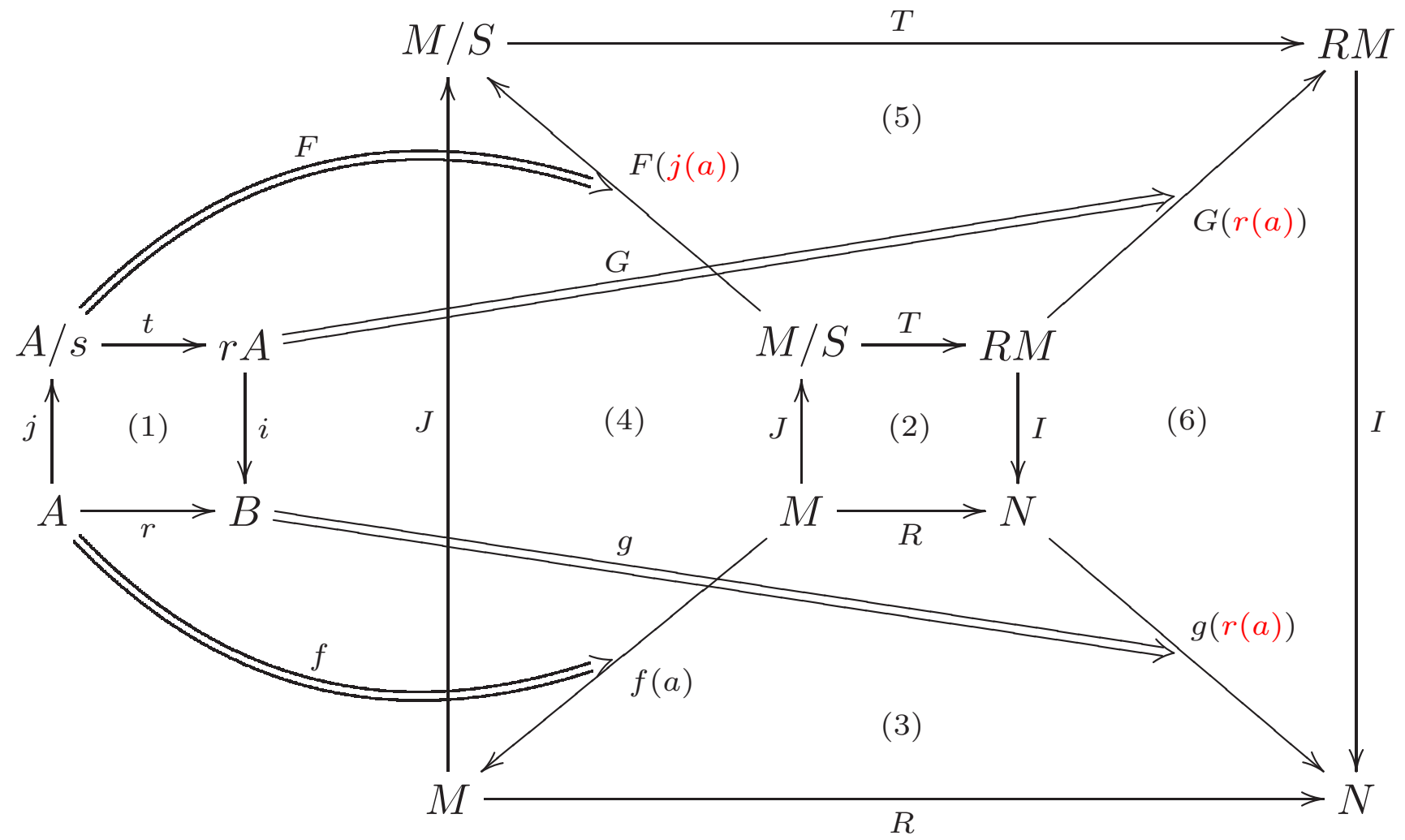

(1) *$s = \ker r$ является конгруэнцией на $A$. Существует разложение гомоморфизма $r$*

$$r = itj \tag{2.13}$$

*$j = \operatorname{nat} s$ - естественный гомоморфизм*

$$j(a) = j(a) \tag{2.14}$$

*$t$ - изоморфизм*

$$r(a) = t(j(a)) \tag{2.15}$$

*$i$ - вложение*

$$r(a) = i(r(a)) \tag{2.16}$$

(2) *$S = \ker R$ является эквивалентностью на $M$. Существует разложение отображения $R$*

$$R = ITJ \tag{2.17}$$

*$J = \operatorname{nat} S$ - сюръекция*

$$J(m) = J(m) \tag{2.18}$$

*$T$ - биекция*

$$R(m) = T(J(m)) \tag{2.19}$$

*$I$ - вложение*

$$R(m) = I(R(m)) \tag{2.20}$$

(3) *$F$ - $T\star$-представление $\Omega_1$-алгебры $A/s$ в $M/S$*
(4) *$G$ - $T\star$-представление $\Omega_1$-алгебры $rA$ в $RM$*
(5) *$(j, J)$ - морфизм представлений $f$ и $F$*
(6) *$(t, T)$ - морфизм представлений $F$ и $G$*
(7) *$(t^{-1}, T^{-1})$ - морфизм представлений $G$ и $F$*
(8) *$(i, I)$ - морфизм представлений $G$ и $g$*
(9) *Существует разложение морфизма представлений*

$$(r, R) = (i, I)(t, T)(j, J) \tag{2.21}$$

*Доказательство.* Существование диаграмм (1) и (2) следует из теоремы II.3.7 ([7], с. 74).

Мы начнём с диаграммы (4).

Пусть $m_1 \equiv m_2 (\operatorname{mod} S)$. Следовательно,

$$R(m_1) = R(m_2) \tag{2.22}$$

Если $a_1 \equiv a_2 (\operatorname{mod} s)$, то

$$r(a_1) = r(a_2) \tag{2.23}$$

Следовательно, $j(a_1) = j(a_2)$. Так как $(r, R)$ - морфизм представлений, то

$$R(f(a_1)(m_1)) = g(r(a_1))(R(m_1)) \tag{2.24}$$

$$R(f(a_2)(m_2)) = g(r(a_2))(R(m_2)) \tag{2.25}$$

Из (2.22), (2.23), (2.24), (2.25) следует

(2.26) $$R(f(a_1)(m_1)) = R(f(a_2)(m_2))$$

Из (2.26) следует

(2.27) $$f(a_1)(m_1) \equiv f(a_2)(m_2)(\mathrm{mod} S)$$

и, следовательно,

(2.28) $$J(f(a_1)(m_1)) = J(f(a_2)(m_2))$$

Из (2.28) следует, что отображение

(2.29) $$F(j(a))(J(m)) = J(f(a)(m)))$$

определено корректно и является преобразованием множества $M/S$.

Из равенства (2.27) (в случае $a_1 = a_2$) следует, что для любого $a$ преобразование согласованно с эквивалентностью $S$. Из теоремы 2.14 следует, что на множестве ${}^*(M/S)$ определена структура $\Omega_1$-алгебры. Рассмотрим $n$-арную операцию $\omega$ и $n$ преобразований

$$F(j(a_i))(J(m)) = J(f(a_i)(m))) \quad i = 1, ..., n$$

пространства $M/S$. Мы положим

$$\omega(F(j(a_1)), ..., F(j(a_n)))(J(m)) = J(\omega(f(a_1), ..., f(a_n)))(m))$$

Следовательно, отображение $F$ является представлением $\Omega_1$-алгебры $A/s$.

Из (2.29) следует, что $(j, J)$ является морфизмом представлений $f$ и $F$ (утверждение (5) теоремы).

Рассмотрим диаграмму (5).

Так как $T$ - биекция, то мы можем отождествить элементы множества $M/S$ и множества $MR$, причём это отождествление имеет вид

(2.30) $$T(J(m)) = R(m)$$

Мы можем записать преобразование $F(j(a))$ множества $M/S$ в виде

(2.31) $$F(j(a)) : J(m) \to F(j(a))(J(m))$$

Так как $T$ - биекция, то мы можем определить преобразование

(2.32) $$T(J(m)) \to T(F(j(a))(J(m)))$$

множества $RM$. Преобразование (2.32) зависит от $j(a) \in A/s$. Так как $t$ - биекция, то мы можем отождествить элементы множества $A/s$ и множества $rA$, причём это отождествление имеет вид

(2.33) $$t(j(a)) = r(a)$$

Следовательно, мы определили отображение

$$G : rA \to {}^{\star}RM$$

согласно равенству

(2.34) $$G(t(j(a)))(T(J(m))) = T(F(j(a))(J(m)))$$

Рассмотрим $n$-арную операцию $\omega$ и $n$ преобразований

$$G(r(a_i))(R(m)) = T(F(j(a_i))(J(m))) \quad i = 1, ..., n$$

пространства $RM$. Мы положим

(2.35) $$\omega(G(r(a_1)), ..., G(r(a_n)))(R(m)) = T(\omega(F(j(a_1), ..., F(j(a_n)))(J(m)))$$

Согласно (2.34) операция $\omega$ корректно определена на множестве ${}^\star RM$. Следовательно, отображение $G$ является представлением $\Omega_1$-алгебры.

Из (2.34) следует, что $(t,T)$ является морфизмом представлений $F$ и $G$ (утверждение (6) теоремы).

Так как $T$ - биекция, то из равенства (2.30) следует

$$J(m)=T^{-1}(R(m)) \tag{2.36}$$

Мы можем записать преобразование $G(r(a))$ множества $RM$ в виде

$$G(r(a)):R(m)\to G(r(a))(R(m)) \tag{2.37}$$

Так как $T$ - биекция, то мы можем определить преобразование

$$T^{-1}(R(m))\to T^{-1}(G(r(a))(R(m))) \tag{2.38}$$

множества $M/S$. Преобразование (2.38) зависит от $r(a)\in rA$. Так как $t$ - биекция, то из равенства (2.33) следует

$$j(a)=t^{-1}(r(a)) \tag{2.39}$$

Так как по построению диаграмма (5) коммутативна, то преобразование (2.38) совпадает с преобразованием (2.31). Равенство (2.35) можно записать в виде

$$T^{-1}(\omega(G(r(a_1)),...,G(r(a_n)))(R(m)))=\omega(F(j(a_1),...,F(j(a_n)))(J(m)) \tag{2.40}$$

Следовательно, $(t^{-1},T^{-1})$ является морфизмом представлений $G$ и $F$ (утверждение (7) теоремы).

Диаграмма (6) является самым простым случаем в нашем доказательстве. Поскольку отображение $I$ является вложением и диаграмма (2) коммутативна, мы можем отождествить $n\in N$ и $R(m)$, если $n\in \mathrm{Im}R$. Аналогично, мы можем отождествить соответствующие преобразования.

$$\begin{aligned} &g'(i(r(a)))(I(R(m)))=I(G(r(a))(R(m))) \\ &\omega(g'(r(a_1)),...,g'(r(a_n)))(R(m))=I(\omega(G(r(a_1),...,G(r(a_n)))(R(m))) \end{aligned} \tag{2.41}$$

Следовательно, $(i,I)$ является морфизмом представлений $G$ и $g$ (утверждение (8) теоремы).

Для доказательства утверждения (9) теоремы осталось показать, что определённое в процессе доказательства представление $g'$ совпадает с представлением $g$, а операции над преобразованиями совпадают с соответствующими операциями на ${}^*N$.

$$\begin{aligned} g'(i(r(a)))(I(R(m)))&=I(G(r(a))(R(m))) && \text{согласно (2.41)}\\ &=I(G(t(j(a)))(T(J(m)))) && \text{согласно (2.15), (2.19)}\\ &=IT(F(j(a))(J(m))) && \text{согласно (2.34)}\\ &=ITJ(f(a)(m)) && \text{согласно (2.29)}\\ &=R(f(a)(m)) && \text{согласно (2.17)}\\ &=g(r(a))(R(m)) && \text{согласно (2.3)} \end{aligned}$$

$$\begin{aligned}\omega(G(r(a_1)),...,G(r(a_n)))(R(m)) &= T(\omega(F(j(a_1)),...,F(j(a_n)))(J(m)))\\ &= T(F(\omega(j(a_1),...,j(a_n)))(J(m)))\\ &= T(F(j(\omega(a_1,...,a_n)))(J(m)))\\ &= T(J(f(\omega(a_1,...,a_n))(m)))\end{aligned}$$

□

**Определение 2.16.** Пусть

$$f:A\to {}^*M$$

представление $\Omega_1$-алгебры $A$,

$$g:B\to {}^*N$$

представление $\Omega_1$-алгебры $B$. Пусть

$$r:A\longrightarrow B \qquad R:M\longrightarrow N$$

морфизм представлений из $r$ в $R$ такой, что $f$ - изоморфизм $\Omega_1$-алгебры и $g$ - изоморфизм $\Omega_2$-алгебры. Тогда отображение $(r,R)$ называется **изоморфизмом представлений**. □

**Теорема 2.17.** *В разложении* (2.21) *отображение* $(t,T)$ *является изоморфизмом представлений* $F$ *и* $G$.

*Доказательство.* Следствие определения 2.16 и утверждений (6) и (7) теоремы 2.15. □

Из теоремы 2.15 следует, что мы можем свести задачу изучения морфизма представлений $\Omega_1$-алгебры к случаю, описываемому диаграммой

(2.42)

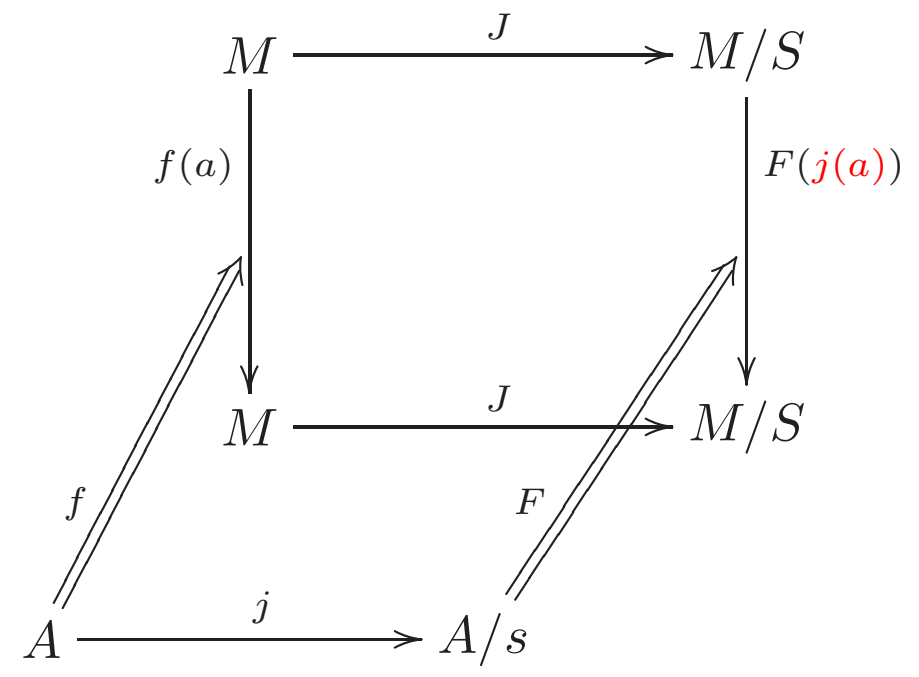

**Теорема 2.18.** *Диаграмма* (2.42) *может быть дополнена представлением* $F_1$ $\Omega_1$*-алгебры* $A$ *в множестве* $M/S$ *так, что диаграмма*

(2.43)

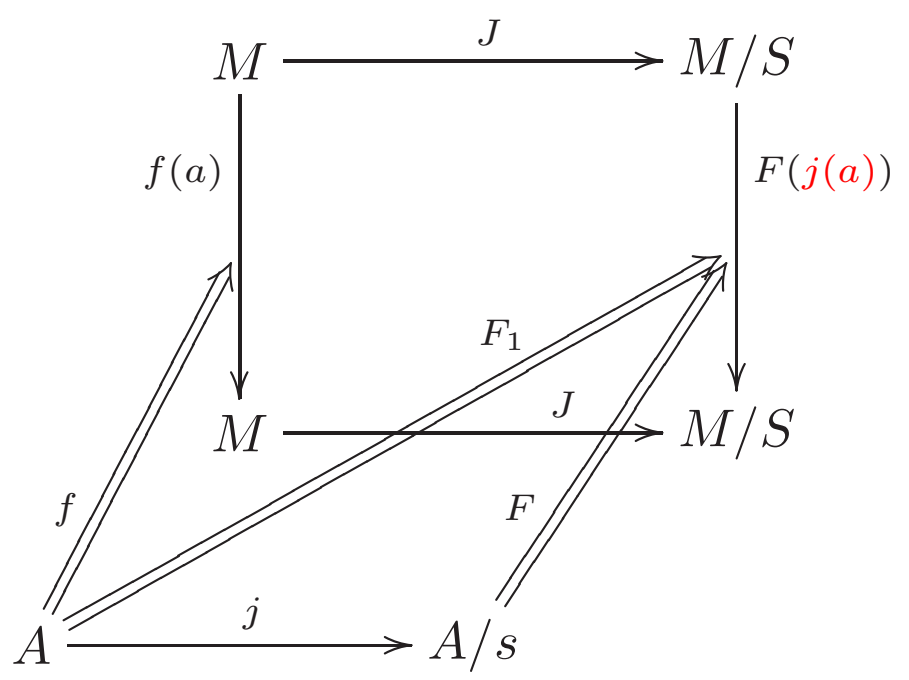

*коммутативна. При этом множество преобразований представления* $F$ *и множество преобразований представления* $F_1$ *совпадают.*

*Доказательство.* Для доказательства теоремы достаточно положить

$$F_1(a) = F(j(a))$$

Так как отображение $j$ - сюрьекция, то $\mathrm{Im}F_1 = \mathrm{Im}F$. Так как $j$ и $F$ - гомоморфизмы $\Omega_1$-алгебры, то $F_1$ - также гомоморфизм $\Omega_1$-алгебры. □

Теорема 2.18 завершает цикл теорем, посвящённых структуре морфизма представлений $\Omega_1$-алгебры. Из этих теорем следует, что мы можем упростить задачу изучения морфизма представлений $\Omega_1$-алгебры и ограничиться морфизмом представлений вида

$$id : A \longrightarrow A \qquad R : M \longrightarrow N$$

В этом случае мы можем отождествить морфизм $(id, R)$ представлений $\Omega_1$-алгебры и отображение $R$. Мы будем пользоваться диаграммой

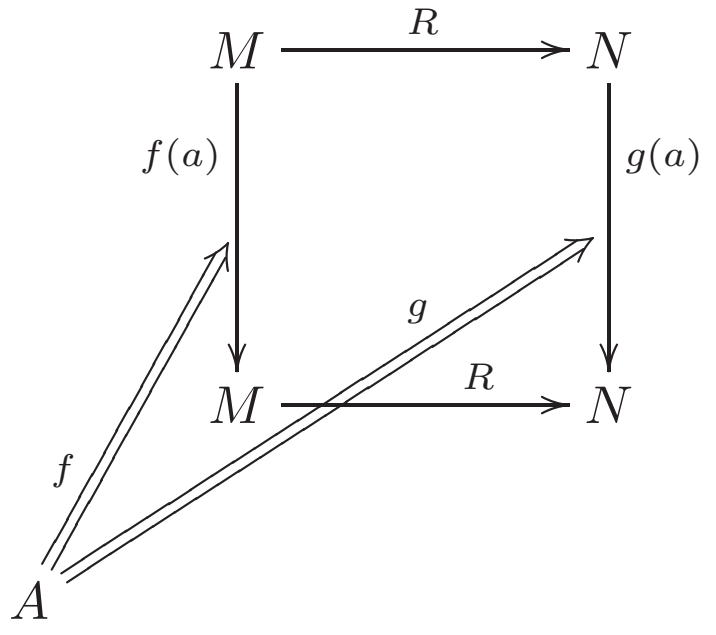

для представления морфизма $(id, R)$ представлений $\Omega_1$-алгебры. Из диаграммы следует

(2.44) $$R \circ f(a) = g(a) \circ R$$

Мы дадим следующее определение по аналогии с определением 2.11.

**Определение 2.19.** Мы определим **категорию** $T\star A$ $\quad$ **$T\star$-представлений $\Omega_1$-алгебры** $A$. Объектами этой категории являются $T\star$-представлениями $\Omega_1$-алгебры $A$. Морфизмами этой категории являются морфизмы $(id, R)$ $T\star$-представлений $\Omega_1$-алгебры $A$. □

## 3. Автоморфизм представления универсальной алгебры

**Определение 3.1.** Пусть
$$f : A \to {}^*M$$
представление $\Omega_1$-алгебры $A$ в $\Omega_2$-алгебре $M$. Морфизм представлений $\Omega_1$-алгебры
$$(r : A \to A, R : M \to M)$$
такой, что $r$ - эндоморфизм $\Omega_1$-алгебры и $R$ - эндоморфизм $\Omega_2$-алгебры называется **эндоморфизмом представления** $f$. □

**Определение 3.2.** Пусть
$$f : A \to {}^*M$$
представление $\Omega_1$-алгебры $A$ в $\Omega_2$-алгебре $M$. Морфизм представлений $\Omega_1$-алгебры
$$(r : A \to A, R : M \to M)$$
такой, что $r$ - автоморфизм $\Omega_1$-алгебры и $R$ - автоморфизм $\Omega_2$-алгебры называется **автоморфизмом представления** $f$. □

**Теорема 3.3.** *Пусть*
$$f : A \to {}^*M$$
*представление $\Omega_1$-алгебры $A$ в $\Omega_2$-алгебре $M$. Множество автоморфизмов представления $f$ порождает* **лупу** $\mathfrak{A}(f)$ .[2]

*Доказательство.* Пусть $(r, R)$, $(p, P)$ - автоморфизмы представления $f$. Согласно определению 3.2 отображения $r$, $p$ являются автоморфизмами $\Omega_1$-алгебры $A$ и отображения $R$, $P$ являются автоморфизмами $\Omega_2$-алгебры $M$. Согласно теореме II.3.2 ([7], с. 60) отображение $rp$ является автоморфизмом $\Omega_1$-алгебры $A$ и отображение $RP$ является автоморфизмом $\Omega_2$-алгебры $M$. Из теоремы 2.10 и определения 3.2 следует, что произведение автоморфизмов $(rp, RP)$ представления $f$ является автоморфизмом представления $f$.

Пусть $(r, R)$ - автоморфизм представления $f$. Согласно определению 3.2 отображение $r$ является автоморфизмом $\Omega_1$-алгебры $A$ и отображение $R$ является автоморфизмом $\Omega_2$-алгебры $M$. Следовательно, отображение $r^{-1}$ является автоморфизмом $\Omega_1$-алгебры $A$ и отображение $R^{-1}$ является автоморфизмом $\Omega_2$-алгебры $M$. Для автоморфизма $(r, R)$ справедливо равенство (2.4). Положим $a' = r(a)$, $m' = R(m)$. Так как $r$ и $R$ - автоморфизмы, то $a = r^{-1}(a')$, $m = R^{-1}(m')$ и равенство (2.4) можно записать в виде

$$R(f(r^{-1}(a'))(R^{-1}(m'))) = g(a')(m') \tag{3.1}$$

аа Так как отображение $R$ является автоморфизмом $\Omega_2$-алгебры $M$, то из равенства (3.1) следует

$$f(r^{-1}(a'))(R^{-1}(m')) = R^{-1}(g(a')(m')) \tag{3.2}$$

Равенство (3.2) соответствует равенству (2.4) для отображения $(r^{-1}, R^{-1})$. Следовательно, отображение $(r^{-1}, R^{-1})$ является автоморфизмом представления $f$. □

[2]Определение лупы приведенно в [3], с. 24, [2], с. 39.

*Замечание* 3.4. Очевидно, что множество автоморфизмов $\Omega_1$-алгебры $A$ также порождает лупу. Конечно, заманчиво предположить, что множество автоморфизмов порождает группу. Так как произведение автоморфизмов $f$ и $g$ является автоморфизмом $fg$, то определены автоморфизмы $(fg)h$ и $f(gh)$. Но из этого утверждения не следует, что

$$(fg)h = f(gh)$$

□

## 4. Векторное пространство

Чтобы определить $T\star$-представление

$$f : R \to {}^\star\overline{M}$$

кольца $R$ на множестве $\overline{M}$ мы должны определить структуру кольца на множестве ${}^\star\overline{M}$.[3]

**Теорема 4.1.** *$T\star$-представление $f$ кольца $R$ на множестве $\overline{M}$ определенно тогда и только тогда, когда определены $T\star$-представления мультипликативной и аддитивной групп кольца $R$ и эти $T\star$-представления удовлетворяют соотношению*

$$f(a(b+c)) = f(a)f(b) + f(a)f(c)$$

*Доказательство.* Теорема следует из определения 1.4. □

**Определение 4.2.** $\overline{M}$ - **$R\star$-модуль** над кольцом $R$, если $\overline{M}$ - абелева группа и определено $T\star$-представление кольца $R$. □

Согласно нашим обозначениям $R\star$-модуль - это **левый модуль** и **$\star R$-модуль** - это **правый модуль**.

Поле является частным случаем кольца. Поэтому векторное пространство над полем имеет больше свойств чем модуль над кольцом. Очень трудно, если вообще возможно, распространить определения, работающие в векторном пространстве, на модуль над произвольным кольцом. Определение базиса и размерности векторного пространства тесно связаны с возможностью найти

[3] Можно ли определить операцию сложения на множестве ${}^\star\overline{M}$, если эта операция не определена на множестве $\overline{M}$. Ответ на этот вопрос положительный.

Допустим $\overline{M} = B \cup C$ и $F : B \to C$ - взаимно однозначное отображение. Мы определим множество ${}^\star\overline{M}$ $T\star$-преобразований множества $\overline{M}$ согласно следующему правилу. Пусть $V \subseteq B$. $T\star$-отображение $F_V$ имеет вид

$$F_V x = \begin{cases} x & x \in B\backslash V \\ Fx & x \in V \\ x & x \in C\backslash F(V) \\ F^{-1}x & x \in F(V) \end{cases}$$

Мы определим сложение $T\star$-преобразований по правилу

$$F_V + F_W = F_{V \triangle W}$$
$$V \triangle W = (V \cup W)\backslash(V \cap W)$$

Очевидно, что

$$F_\emptyset + F_V = F_V$$
$$F_V + F_V = F_\emptyset$$

Следовательно, отображение $F_\emptyset$ является нулём относительно сложения, а множество ${}^\star\overline{M}$ является абелевой группой.

решение линейного уравнения в кольце. Свойства линейного уравнения в теле близки к свойствам линейного уравнения в поле. Поэтому мы надеемся что свойства векторного пространства над телом близки к свойствам векторного пространства над полем.

**Теорема 4.3.** *$T\star$-представление тела $D$* **эффективно**, *если эффективно $T\star$-представление мультипликативной группы тела $D$.*

*Доказательство.* Пусть

$$f : D \to {}^{\star}M \tag{4.1}$$

$T\star$-представление тела $D$. Если элементы $a$, $b$ мультипликативной группы порождают одно и то $T\star$-преобразование, то

$$f(a)m = f(b)m \tag{4.2}$$

для любого $m \in M$. Выполняя преобразование $f(a^{-1})$ над обеими частями равенства (4.2), мы получим

$$m = f(a^{-1})(f(b)m) = f(a^{-1}b)m$$

□

Согласно замечанию 1.9, если представление тела эффективно, мы отождествляем элемент тела и соответствующее ему $T\star$-преобразование.

**Определение 4.4.** $\overline{V}$ - **$D\star$-векторное пространство** над телом $D$, если $\overline{V}$ - абелева группа и определено эффективное $T\star$-представление тела $D$. □

**Теорема 4.5.** *Элементы $D\star$-векторного пространства $\overline{V}$ удовлетворяют соотношениям*

- **закону ассоциативности**

$$(ab)\overline{m} = a(b\overline{m}) \tag{4.3}$$

- **закону дистрибутивности**

$$a(\overline{m} + \overline{n}) = a\overline{m} + a\overline{n} \tag{4.4}$$

$$(a + b)\overline{m} = a\overline{m} + b\overline{m} \tag{4.5}$$

- **закону унитарности**

$$1\overline{m} = \overline{m} \tag{4.6}$$

*для любых $a, b \in D$, $\overline{m}, \overline{n} \in \overline{V}$. Мы будем называть $T\star$-представление $D\star$-* **произведением вектора на скаляр**.

*Доказательство.* Равенство (4.4) следует из утверждения, что $T\star$-преобразование $a$ является автоморфизмом абелевой группы. Равенство (4.5) следует из утверждения, что $T\star$-представление является гомоморфизмом аддитивной группы тела $D$. Равенства (4.3) и (4.6) следуют из утверждения, что $T\star$-представление является ковариантным $T\star$-представлением мультипликативной группы тела $D$. □

Согласно нашим обозначениям $D\star$-векторное пространство - это **левое векторное пространство** и **$\star D$-векторное пространство** - это **правое векторное пространство**.

**Определение 4.6.** Пусть $\overline{V}$ - $D\star$-векторное пространство над телом $D$. Множество векторов $\overline{N}$ - **подпространство $D\star$-векторного пространства** $\overline{V}$, если

$$\overline{a}+\overline{b}\in\overline{N}$$
$$k\overline{a}\in\overline{N}$$
$$\overline{a},\overline{b}\in\overline{N}\quad k\in D$$

□

**Пример 4.7.** Определим на множестве $D_n^m$ $m\times n$ матриц над телом $D$ операцию сложения

$$a+b=\left(a_i^j\right)+\left(b_i^j\right)=\left(a_i^j+b_i^j\right)$$

и умножения на скаляр

$$da=d\left(a_i^j\right)=\left(da_i^j\right)$$

$a=0$ тогда и только тогда, когда $a_i^j=0$ для любых $i$, $j$. Непосредственная проверка показывает, что $D_n^m$ является $D\star$-векторным пространством, если произведение действует слева. В противном случае $D_n^m$ - $\star D$-векторное пространство. Мы будем называть векторное пространство $D_n^m$ **$D\star$-векторным пространством матриц**. □

## 5. Тип векторного пространства

Произведение вектора на скаляр асимметрично, так как произведение определено для объектов разных множеств. Однако различие между $D\star$- и $\star D$ -векторным пространством появляется только тогда, когда мы переходим к координатному представлению. Говоря, векторное пространство является $D\star$ - или $\star D$-, мы указываем, с какой стороны, слева или справа, мы умножаем координаты вектора на элементы тела.

**Определение 5.1.** Допустим $u$, $v$ - векторы $D\star$-векторного пространства $\overline{V}$. Мы будем говорить, что вектор $w$ является **$D\star$-линейной комбинацией векторов** $u$ и $v$, если мы можем записать $w=au+bv$, где $a$ и $b$ - скаляры. □

Мы можем распространить понятие линейной комбинации на любое конечное семейство векторов. Пользуясь обобщённым индексом для нумерации векторов, мы можем представить семейство векторов в виде одномерной матрицы. Мы пользуемся соглашением, что в заданном векторном пространстве мы представляем любое семейство векторов либо в виде ${}_*$-строки, либо ${}^*$-строки. Это представление определяет характер записи линейной комбинации. Рассматривая это представление в $D\star$- или $\star D$-векторном пространстве, мы получаем четыре разных модели векторного пространства.

Чтобы иметь возможность, не меняя записи указать, является векторное пространство $D\star$- или $\star D$-векторным пространством, мы вводим новое обозначение. Символ $D\star$ называется **типом векторного пространства** и означает, что мы изучаем $D\star$-векторное пространство над телом $D$. Операция умножения в типе векторного пространства указывает на матричную операцию, используемую в линейной комбинации.

**Пример 5.2.** Представим множество векторов ${}^i\overline{a}$, $i \in I$, $D_*{}^*$**-векторного пространства** $\overline{V}$ в виде ${}^*$-строки

$$\overline{a} = \begin{pmatrix} {}^1\overline{a} \\ ... \\ {}^n\overline{a} \end{pmatrix}$$

и множество скаляров $c_i$, $i \in I$, в виде ${}_*$-строки

$$c = \begin{pmatrix} c_1 & ... & c_n \end{pmatrix}$$

Тогда мы можем записать линейную комбинацию векторов ${}^i\overline{a}$ в виде

$$c_i\ {}^i\overline{a} = c_*{}^*\overline{a}$$

Мы пользуемся записью ${}_{D_*{}^*}\overline{V}$, когда мы хотим сказать, что $\overline{V}$ - $D_*{}^*$-векторное пространство. □

**Пример 5.3.** Представим множество векторов ${}_i\overline{a}$, $i \in I$, $D^*{}_*$**-векторного пространства** $\overline{V}$ в виде ${}_*$-строки

$$\overline{a} = \begin{pmatrix} {}_1\overline{a} & ... & {}_n\overline{a} \end{pmatrix}$$

и множество скаляров $c^i$, $i \in I$, в виде ${}^*$-строки

$$c = \begin{pmatrix} c^1 \\ ... \\ c^n \end{pmatrix}$$

Тогда мы можем записать линейную комбинацию векторов ${}_i\overline{a}$ в виде

$$c^i\ {}_i\overline{a} = c^*{}_*\overline{a}$$

Мы пользуемся записью ${}_{D^*{}_*}\overline{V}$, когда мы хотим сказать, что $\overline{V}$ - $D^*{}_*$-векторное пространство. □

**Пример 5.4.** Представим множество векторов $\overline{a}^i$, $i \in I$, ${}^*{}_*D$**-векторного пространства** $\overline{V}$ в виде ${}^*$-строки

$$\overline{a} = \begin{pmatrix} \overline{a}^1 \\ ... \\ \overline{a}^n \end{pmatrix}$$

и множество скаляров ${}_ic$, $i \in I$, в виде ${}_*$-строки

$$c = \begin{pmatrix} {}_1c & ... & {}_nc \end{pmatrix}$$

Тогда мы можем записать линейную комбинацию векторов $\overline{a}^i$ в виде

$$\overline{a}^i\ {}_ic = \overline{a}^*{}_*c$$

Мы пользуемся записью $\overline{V}_{{}^*{}_*D}$, когда мы хотим сказать, что $\overline{V}$ - ${}^*{}_*D$-векторное пространство. □

**Пример 5.5.** Представим множество векторов $\overline{A}_i$, $i \in I$, ${}_*{}^*D$**-векторного пространства** $\overline{V}$ в виде ${}_*$-строки

$$\overline{a} = \begin{pmatrix} \overline{a}_1 & \dots & \overline{a}_n \end{pmatrix}$$

и множество скаляров ${}_ic$, $i \in I$, в виде ${}^*$-строки

$$c = \begin{pmatrix} {}^1c \\ \dots \\ {}^nc \end{pmatrix}$$

Тогда мы можем записать линейную комбинацию векторов $\overline{a}_i$ в виде

$$\overline{a}_i\ {}^ic = \overline{a}_*{}^*c$$

Мы пользуемся записью $\overline{V}_{{}_*{}^*D}$, когда мы хотим сказать, что $\overline{V}$ - ${}_*{}^*D$-векторное пространство. □

*Замечание* 5.6. Мы распространим на векторное пространство и его тип соглашение, описанное в замечании [6]-2.15. Например, в выражении

$$A_*{}^*B_*{}^*v\lambda$$

мы выполняем операцию умножения слева направо. Это соответствует ${}_*{}^*D$-векторному пространству. Однако мы можем выполнять операцию умножения справа налево. В традиционной записи это выражение примет вид

$$\lambda v^*{}_*B^*{}_*A$$

и будет соответствовать $D^*{}_*$-векторному пространству. Аналогично, читая выражение снизу вверх мы получим выражение

$$A^*{}_*B^*{}_*v\lambda$$

соответствующего ${}^*{}_*D$-векторному пространству. □

## 6. ${}_*{}^*D$-базис векторного пространства

**Определение 6.1.** Векторы $\overline{A}_i$, $i \in I$, ${}_*{}^*D$-векторного пространства $\overline{V}$ ${}_*{}^*D$**-линейно независимы**, если $c = 0$ следует из уравнения

$$\overline{A}_*{}^*c = 0$$

В противном случае, векторы $\overline{A}_i$ ${}_*{}^*D$**-линейно зависимы**. □

**Определение 6.2.** Множество векторов $\overline{\overline{e}} = (\overline{e}_i, i \in I)$ - ${}_*{}^*D$**-базис векторного пространства**, если векторы $\overline{e}_i$ ${}_*{}^*D$-линейно независимы и добавление любого вектора к этой системе делает эту систему ${}_*{}^*D$-линейно зависимой. □

**Теорема 6.3.** *Если $\overline{\overline{e}}$ - ${}_*{}^*D$-базис векторного пространства $\overline{V}$, то любой вектор $\overline{v} \in \overline{V}$ имеет одно и только одно разложение*

$$\overline{v} = \overline{e}_*{}^*v \tag{6.1}$$

*относительно этого ${}_*{}^*D$-базиса.*

*Доказательство.* Так как система векторов $\overline{e}_i$ является максимальным множеством ${}_*{}^*D$-линейно независимых векторов, система векторов $\overline{v}, \overline{e}_i$ - $D_*{}^*$-линейно зависима и в уравнении

$$\overline{v}b + \overline{e}_*{}^*c = 0 \tag{6.2}$$

по крайней мере $b$ отлично от 0. Тогда равенство

$$\overline{v} = \overline{e}_*{}^*(-cb^{-1}) \tag{6.3}$$

следует из (6.2). (6.1) следует из (6.3).

Допустим мы имеем другое разложение

$$\overline{v} = \overline{e}_*{}^*v' \tag{6.4}$$

Вычтя (6.1) из (6.4), мы получим

$$0 = \overline{e}_*{}^*(v' - v)$$

Так как векторы $\overline{e}_i$ ${}_*{}^*D$-линейно независимы, мы имеем

$$v' - v = 0$$

□

**Определение 6.4.** Мы будем называть матрицу $v$ разложения (6.1) **координатной матрицей вектора $\overline{v}$ в ${}_*{}^*D$-базисе $\overline{\overline{e}}$** и её элементы **координатами вектора $\overline{v}$ в ${}_*{}^*D$-базисе $\overline{\overline{e}}$.** □

**Теорема 6.5.** *Множество координат $a$ вектора $\overline{a}$ в ${}_*{}^*D$-базисе $\overline{\overline{e}}$ порождают ${}_*{}^*D$-векторное пространство $D^n$, изоморфное ${}_*{}^*D$-векторному пространству $\overline{V}$. Это ${}_*{}^*D$-векторное пространство называется* **координатным ${}_*{}^*D$-векторным пространством**, *а изоморфизм* **координатным ${}_*{}^*D$-изоморфизмом**.

*Доказательство.* Допустим векторы $\overline{a}$ и $\overline{b} \in \overline{V}$ имеют разложение

$$\overline{a} = \overline{e}_*{}^*a$$
$$\overline{b} = \overline{e}_*{}^*b$$

в базисе $\overline{\overline{e}}$. Тогда

$$\overline{a} + \overline{b} = \overline{e}_*{}^*a + \overline{e}_*{}^*b = \overline{e}_*{}^*(a + b)$$
$$\overline{a}m = (\overline{e}_*{}^*a)m = \overline{e}_*{}^*(ma)$$

для любого $m \in D$. Таким образом, операции в векторном пространстве определены по координатно

$${}^i(a + b) = {}^ia + {}^ib$$
$${}^i(am) = {}^iam$$

Это доказывает теорему. □

**Пример 6.6.** Пусть $\overline{\overline{e}} = (\overline{e}_i, i \in I, |I| = n)$ - ${}_*{}^*D$-базис векторного пространства $D^n$. Координатная матрица

$$a = \begin{pmatrix} {}^1a \\ ... \\ {}^na \end{pmatrix} = ({}^ia, i \in I)$$

вектора $\overline{a}$ в ${}_*{}^*D$-базисе $\overline{\overline{e}}$ называется ${}_*{}^*D$**-вектором**.[4] Мы будем называть векторное пространство $D^n$ ${}_*{}^*D$**-векторным пространством** ${}^*$**-строк**.

Пусть ${}_*$-строка

$$\overline{A} = \begin{pmatrix} \overline{A}_1 & ... & \overline{A}_m \end{pmatrix} = (\overline{A}_j, j \in J) \tag{6.5}$$

задаёт множество векторов. Векторы $\overline{A}_j$ имеют разложение

$$\overline{A}_j = \overline{e}_*{}^* A_j$$

Если мы подставим координатные матрицы вектора $\overline{A}_j$ в матрицу (6.5), мы получим матрицу

$$A = \begin{pmatrix} \begin{pmatrix} {}^1A_1 \\ ... \\ {}^nA_1 \end{pmatrix} & ... & \begin{pmatrix} {}^1A_m \\ ... \\ {}^nA_m \end{pmatrix} \end{pmatrix} = \begin{pmatrix} {}^1A_1 & ... & {}^1A_m \\ ... & ... & ... \\ {}^nA_1 & ... & {}^nA_m \end{pmatrix} = ({}^iA_j)$$

Мы будем называть матрицу $A$ **координатной матрицей множества векторов** $(\overline{A}_j, j \in J)$ в базисе $\overline{\overline{e}}$ и её элементы **координатами множества векторов** $(\overline{A}_j, j \in J)$ в базисе $\overline{\overline{e}}$.

Мы будем представлять ${}_*{}^*D$**-базис** $\overline{\overline{f}}$ **векторного пространства** ${}^*$**-строк** $D^n$ в форме ${}_*$-строки

$$\overline{f} = \begin{pmatrix} \overline{f}_1 & ... & \overline{f}_n \end{pmatrix} = (\overline{f}_j, j \in I)$$

Мы будем говорить, что координатная матрица $f$ множества векторов $(\overline{f}_j, j \in I)$ определяет координаты ${}^if_j$ базиса $\overline{\overline{f}}$ относительно базиса $\overline{\overline{e}}$. □

**Пример 6.7.** Пусть $\overline{\overline{e}} = ({}^j\overline{e}, j \in J, |J| = m)$ - $D_*{}^*$-базис векторного пространства $D_n$. координатная матрица

$$a = \begin{pmatrix} a_1 & ... & a_m \end{pmatrix} = (a_j, j \in J)$$

вектора $\overline{a}$ в $D_*{}^*$-базисе $\overline{\overline{e}}$ называется $D_*{}^*$**-вектором**.[5] Мы будем называть векторное пространство $D_n$ $D_*{}^*$**-векторным пространством** ${}_*$**-строк**.

Пусть ${}^*$-строка

$$\overline{A} = \begin{pmatrix} {}^1\overline{A} \\ ... \\ {}^n\overline{A} \end{pmatrix} = ({}^i\overline{A}, i \in I) \tag{6.6}$$

задаёт множество векторов. Векторы ${}^i\overline{A}$ имеют разложение

$${}^i\overline{A} = {}^iA_*{}^*\overline{e}$$

[4] ${}_*{}^*D$-вектор является аналогом вектор-столбца. Мы можем называть его также $D\star$-вектор ${}^*$-строкой

[5] $D_*{}^*$-вектор является аналогом вектор-строки. Мы можем называть его также $D\star$-вектор ${}_*$-строкой

Если мы подставим координатные матрицы вектора ${}^i\overline{A}$ в матрицу (6.6), мы получим матрицу

$$A = \begin{pmatrix} \begin{pmatrix} {}^1A_1 & ... & {}^1A_m \end{pmatrix} \\ ... \\ \begin{pmatrix} {}^nA_1 & ... & {}^nA_m \end{pmatrix} \end{pmatrix} = \begin{pmatrix} {}^1A_1 & ... & {}^1A_m \\ ... & ... & ... \\ {}^nA_1 & ... & {}^nA_m \end{pmatrix} = ({}^iA_j)$$

Мы будем называть матрицу $A$ **координатной матрицей множества векторов** $({}^i\overline{A}, i \in M)$ в базисе $\overline{\overline{e}}$ и её элементы **координатами множества векторов** $({}^i\overline{A}, i \in M)$ в базисе $\overline{\overline{e}}$.

Мы будем представлять $D_*{}^*$**-базис** $\overline{\overline{f}}$ **векторного пространства** ${}_*$**-строк** $D_n$ в форме ${}^*$-строки

$$\overline{f} = \begin{pmatrix} {}^1\overline{f} & ... & {}^n\overline{f} \end{pmatrix} = ({}^i\overline{f}, i \in J)$$

Мы будем говорить, что координатная матрица $f$ множества векторов $({}^i\overline{f}, i \in J)$ определяет координаты ${}^if_j$ базиса $\overline{\overline{f}}$ относительно базиса $\overline{\overline{e}}$. □

Так как мы линейную комбинацию выражаем с помощью матриц, мы можем распространить принцип двойственности на теорию векторных пространств. Мы можем записать принцип двойственности в одной из следующих форм

**Теорема 6.8** (принцип двойственности)**.** *Пусть $\mathfrak{A}$ - истинное утверждение о векторных пространствах. Если мы заменим одновременно*

- *$D_*{}^*$-вектор и $D^*{}_*$-вектор*
- *${}_*{}^*D$-вектор и ${}^*{}_*D$-вектор*
- *${}_*{}^*$-произведение и ${}^*{}_*$-произведение*

*то мы снова получим истинное утверждение.*

**Теорема 6.9** (принцип двойственности)**.** *Пусть $\mathfrak{A}$ - истинное утверждение о векторных пространствах. Если мы одновременно заменим*

- *$D_*{}^*$-вектор и ${}_*{}^*D$-вектор или $D^*{}_*$-вектор и ${}^*{}_*D$-вектор*
- *${}_*{}^*$-квазидетерминант и ${}^*{}_*$-квазидетерминант*

*то мы снова получим истинное утверждение.*

## 7. ${}_*{}^*D$-линейное отображение векторных пространств

**Определение 7.1.** Пусть $\overline{V}$ - ${}_*{}^*S$-векторное пространство. Пусть $\overline{U}$ - ${}_*{}^*T$-векторное пространство. Мы будем называть морфизм

$$f : S \longrightarrow T \qquad \overline{A} : \overline{V} \longrightarrow \overline{U}$$

$T\star$-представлений тела в абелевой группе $({}_*{}^*S, {}_*{}^*T)$**-линейным отображением векторных пространств**. □

Согласно теореме 2.15 при изучении $({}_*{}^*S, {}_*{}^*T)$-линейного отображения мы можем ограничиться случаем $S = T$.

**Определение 7.2.** Пусть $\overline{V}$ и $\overline{W}$ - ${}_*{}^*D$-векторные пространства. Мы будем называть отображение

$$\overline{A} : \overline{V} \to \overline{W} \tag{7.1}$$

${}_*{}^*D$**-линейным отображением векторных пространств**, если[6]

$$\overline{A}(\overline{m}_*{}^*a)=\overline{A}(\overline{m})_*{}^*a \tag{7.2}$$

для любых ${}^ia\in D$, $\overline{m}_i\in\overline{V}$. □

**Теорема 7.3.** *Пусть*

$$\overline{\overline{f}}=(\overline{f}_i,i\in I)$$

${}_*{}^*D$*-базис в векторном пространстве* $\overline{V}$ *и*

$$\overline{\overline{e}}=(\overline{e}_j,j\in J)$$

${}_*{}^*D$*-базис в векторном пространстве* $\overline{U}$*. Тогда* ${}_*{}^*D$*-линейное отображение* (7.1) *векторных пространств имеет представление*

$$b=A_*{}^*a \tag{7.3}$$

*относительно заданных базисов. Здесь*

- $a$ *- координатная матрица вектора* $\overline{a}$ *относительно* ${}_*{}^*D$*-базиса* $\overline{\overline{f}}$
- $b$ *- координатная матрица вектора*
$$\overline{b}=\overline{A}(\overline{a})$$
*относительно* ${}_*{}^*D$*-базиса* $\overline{\overline{e}}$
- $A$ *- координатная матрица множества векторов* $(\overline{A}(\overline{f}_i))$ *в* ${}_*{}^*D$*-базисе* $\overline{\overline{e}}$*, которую мы будем называть* **матрицей** ${}_*{}^*D$**-линейного отображения** *относительно базисов* $\overline{\overline{f}}$ *и* $\overline{\overline{e}}$

*Доказательство.* Вектор $\overline{a}\in\overline{V}$ имеет разложение

$$\overline{a}=\overline{f}_*{}^*a$$

относительно ${}_*{}^*D$-базиса $\overline{\overline{f}}$. Вектор $\overline{b}=f(\overline{a})\in\overline{U}$ имеет разложение

$$\overline{b}=\overline{e}_*{}^*b \tag{7.4}$$

относительно ${}_*{}^*D$-базиса $\overline{\overline{e}}$.

Так как $\overline{A}$ - ${}_*{}^*D$-линейное отображение, то на основании (7.2) следует, что

$$\overline{b}=\overline{A}(\overline{a})=\overline{A}(\overline{f}_*{}^*a)=\overline{A}(\overline{f})_*{}^*a \tag{7.5}$$

$\overline{A}(\overline{f}_i)$ также вектор векторного пространства $\overline{U}$ и имеет разложение

$$\overline{A}(\overline{f}_i)=\overline{e}_*{}^*A_i=\overline{e}_j\,{}^jA_i \tag{7.6}$$

относительно базиса $\overline{\overline{e}}$. Комбинируя (7.5) и (7.6), мы получаем

$$\overline{b}=\overline{e}_*{}^*A_*{}^*a \tag{7.7}$$

(7.3) следует из сравнения (7.4) и (7.7) и теоремы 6.3. □

На основании теоремы 7.3 мы идентифицируем ${}_*{}^*D$-линейное отображение (7.1) векторных пространств и матрицу его представления (7.3).

[6]Выражение $\overline{A}(\overline{m})_*{}^*a$ означает выражение $\overline{A}(\overline{m}_i)\,{}^ia$

**Теорема 7.4.** *Пусть*

$$\overline{\overline{f}}=(\overline{f}_i, i\in I)$$

${}_*{}^*D$*-базис в векторном пространстве* $\overline{V}$*,*

$$\overline{\overline{e}}=(\overline{e}_j, j\in J)$$

${}_*{}^*D$*-базис в векторном пространстве* $\overline{U}$*, и*

$$\overline{\overline{g}}=(\overline{g}_l, l\in L)$$

${}_*{}^*D$*-базис в векторном пространстве* $\overline{W}$*. Предположим, что мы имеем коммутативную диаграмму отображений*

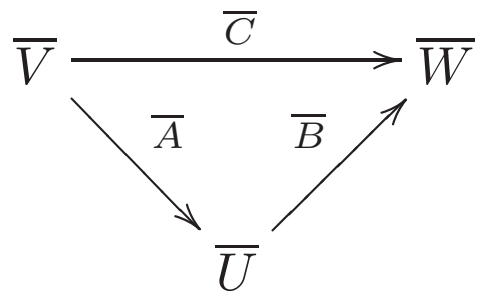

*где* ${}_*{}^*D$*-линейное отображение* $\overline{A}$ *имеет представление*

$$b=A_*{}^*a \tag{7.8}$$

*относительно заданных базисов и* ${}_*{}^*D$*-линейное отображение* $\overline{B}$ *имеет представление*

$$c=B_*{}^*b \tag{7.9}$$

*относительно заданных базисов. Тогда отображение* $\overline{C}$ *является* ${}_*{}^*D$*-линейным и имеет представление*

$$c=B_*{}^*A_*{}^*a \tag{7.10}$$

*относительно заданных базисов.*

*Доказательство.* Отображение $\overline{C}$ является ${}_*{}^*D$-линейным, так как

$$\begin{aligned}\overline{C}_*{}^*(\overline{f}_*{}^*a)&=(\overline{A}_*{}^*\overline{B})_*{}^*(\overline{f}_*{}^*a)\\&=\overline{B}_*{}^*(\overline{A}_*{}^*(\overline{f}_*{}^*a))\\&=\overline{B}_*{}^*(\overline{e}_*{}^*(A_*{}^*a))\\&=\overline{g}_*{}^*(B_*{}^*(A_*{}^*a))\\&=\overline{g}_*{}^*((B_*{}^*A)_*{}^*a)\\&=\overline{g}_*{}^*(C_*{}^*a)\end{aligned}$$

Равенство (7.10) следует из подстановки (7.8) в (7.9). □

Записывая ${}_*{}^*D$-линейное отображение в форме ${}_*{}^*$-произведения, мы можем переписать (7.2) в виде

$$\overline{A}_*{}^*(\overline{a}k)=(\overline{A}_*{}^*\overline{a})k \tag{7.11}$$

Утверждение теоремы 7.4 мы можем записать в виде

$$\overline{B}_*{}^*(\overline{A}_*{}^*\overline{a})=(\overline{B}_*{}^*\overline{A})_*{}^*\overline{a} \tag{7.12}$$

Равенства (7.11) и (7.12) представляют собой **закон ассоциативности для** ${}_*{}^*D$**-линейных отображений векторных пространств**. Это позволяет нам писать подобные выражения не пользуясь скобками.

Равенство (7.3) является координатной записью ${}_*{}^*D$-линейного отображения. На основе теоремы 7.3 бескоординатная запись также может быть представлена с помощью ${}_*{}^*$-произведения

$$\overline{b}=\overline{A}_*{}^*\overline{a}=\overline{A}_*{}^*\overline{f}_*{}^*a=\overline{e}_*{}^*A_*{}^*a \tag{7.13}$$

Если подставить равенство (7.13) в теорему 7.4, то мы получим цепочку равенств

$$\overline{c}=\overline{B}_*{}^*\overline{b}=\overline{B}_*{}^*\overline{e}_*{}^*b=\overline{g}_*{}^*B_*{}^*b$$
$$\overline{c}=\overline{B}_*{}^*\overline{A}_*{}^*\overline{a}=\overline{B}_*{}^*\overline{A}_*{}^*\overline{f}_*{}^*a=\overline{g}_*{}^*B_*{}^*A_*{}^*a$$

*Замечание* 7.5. На примере ${}_*{}^*D$-линейных отображений легко видеть насколько теорема 2.15 облегчает наши рассуждения при изучении морфизма $T\star$-представлений $\Omega$-алгебры. Договоримся в рамках этого замечания теорию ${}_*{}^*D$-линейных отображений называть сокращённой теорией, а теорию, излагаемую в этом замечании, называть расширенной теорией.

Пусть $\overline{V}$ - ${}_*{}^*S$-векторное пространство. Пусть $\overline{U}$ - ${}_*{}^*T$-векторное пространство. Пусть

$$r:S\longrightarrow T \qquad \overline{A}:\overline{V}\longrightarrow\overline{U}$$

$({}_*{}^*S,{}_*{}^*T)$-линейное отображение векторных пространств. Пусть

$$\overline{\overline{f}}=(\overline{f}_i,i\in I)$$

${}_*{}^*S$-базис в векторном пространстве $\overline{V}$ и

$$\overline{\overline{e}}=(\overline{e}_j,j\in J)$$

${}_*{}^*T$-базис в векторном пространстве $\overline{U}$.

Из определений 7.1 и 2.2 следует

$$\overline{b}=\overline{A}(\overline{a})=\overline{A}(\overline{f}_*{}^*a)=\overline{A}(\overline{f})_*{}^*r(a) \tag{7.14}$$

$\overline{A}(\overline{f}_i)$ также вектор векторного пространства $\overline{U}$ и имеет разложение

$$\overline{A}(\overline{f}_i)=\overline{e}_*{}^*A_i=\overline{e}_j{}^jA_i \tag{7.15}$$

относительно базиса $\overline{\overline{e}}$. Комбинируя (7.14) и (7.15), мы получаем

$$\overline{b}=\overline{e}_*{}^*A_*{}^*r(a) \tag{7.16}$$

Пусть $\overline{W}$ - ${}_*{}^*D$-векторное пространство. Пусть

$$p:T\longrightarrow D \qquad \overline{B}:\overline{U}\longrightarrow\overline{W}$$

$({}_*{}^*T,{}_*{}^*D)$-линейное отображение векторных пространств. Пусть

$$\overline{\overline{g}}=(\overline{g}_l,l\in L)$$

${}_*{}^*D$-базис в векторном пространстве $\overline{W}$. Тогда, согласно (7.16), произведение $({}_*{}^*S,{}_*{}^*T)$-линейного отображения $(r,\overline{A})$ и $({}_*{}^*T,{}_*{}^*D)$-линейного отображения $(p,\overline{B})$ имеет вид

$$\overline{c}=\overline{h}_*{}^*B_*{}^*p(A)_*{}^*pr(a) \tag{7.17}$$

Сопоставление равенств (7.10) и (7.17) показывает насколько расширеная теория линейных отображений сложнее сокращённой теории.

При необходимости мы можем пользоваться расширеной теорией, но мы не получим новых результатов по сравнению со случаем сокращённой теорией. В то же время обилие деталей делает картину менее ясной и требует постоянного внимания. □

## 8. Система ${}_*{}^*D$-линейных уравнений

**Определение 8.1.** Пусть $\overline{V}$ - ${}_*{}^*D$-векторное пространство и $\{\overline{A}_i \in \overline{V}, i \in I\}$ - множество векторов. **${}_*{}^*D$-линейная оболочка в векторном пространстве** - это множество $\mathrm{span}(\overline{A}_i, i \in I)$ векторов, ${}_*{}^*D$-линейно зависимых от векторов $\overline{A}_i$. □

**Теорема 8.2.** *Пусть $span(\overline{A}_i, i \in I)$ - ${}_*{}^*D$-линейная оболочка в векторном пространстве $\overline{V}$. Тогда $span(\overline{A}_i, i \in I)$ - подпространство векторного пространства $\overline{V}$.*

*Доказательство.* Предположим, что

$$\overline{b} \in \mathrm{span}(\overline{A}_i, i \in I)$$

$$\overline{c} \in \mathrm{span}(\overline{A}_i, i \in I)$$

Согласно определению 8.1

$$\overline{b} = \overline{A}_*{}^*b$$

$$\overline{c} = \overline{A}_*{}^*c$$

Тогда

$$\overline{b} + \overline{c} = \overline{A}_*{}^*b + \overline{A}_*{}^*c = \overline{A}_*{}^*(b+c) \in \mathrm{span}(\overline{A}_i, i \in I)$$

$$\overline{b}k = (\overline{A}_*{}^*b)k = \overline{A}_*{}^*(bk) \in \mathrm{span}(\overline{A}_i, i \in I)$$

Это доказывает утверждение. □

**Пример 8.3.** Пусть $\overline{V}$ - ${}_*{}^*D$-векторное пространство и ${}_*$-строка

$$\overline{\overline{A}} = \begin{pmatrix} \overline{A}_1 & ... & \overline{A}_n \end{pmatrix} = (\overline{A}_i, i \in I)$$

задаёт множество векторов. Чтобы ответить на вопрос, или вектор $\overline{b} \in \mathrm{span}(\overline{A}_i, i \in I)$, мы запишем линейное уравнение

$$\overline{b} = \overline{A}_*{}^*x \tag{8.1}$$

где

$$x = \begin{pmatrix} {}^1x \\ ... \\ {}^nx \end{pmatrix}$$

${}^*$-строка неизвестных коэффициентов разложения. $\overline{b} \in \mathrm{span}(\overline{A}_i, i \in I)$, если уравнение (8.1) имеет решение. Предположим, что $\overline{\overline{f}} = (\overline{f}_j, j \in J)$ - ${}_*{}^*D$-базис. Тогда векторы $\overline{b}$, $\overline{A}_i$ имеют разложение

$$\overline{b} = b_*{}^*\overline{f} \tag{8.2}$$

$$\overline{A}_i = \overline{f}_*{}^*A_i \tag{8.3}$$

Если мы подставим (8.2) и (8.3) в (8.1), мы получим

$$\overline{f}_*{}^*b = \overline{f}_*{}^*A_*{}^*x \tag{8.4}$$

Применяя теорему 6.3 к (8.4), мы получим **систему ${}_*{}^*D$-линейных уравнений**

$$A_*{}^*x = b \tag{8.5}$$

Мы можем записать систему ${}_*{}^*D$-линейных уравнений (8.5) в одной из следующих форм

$$\begin{pmatrix} {}^1A_1 & ... & {}^1A_n \\ ... & ... & ... \\ {}^mA_1 & ... & {}^mA_n \end{pmatrix} {}_*{}^* \begin{pmatrix} {}^1x \\ ... \\ {}^nx \end{pmatrix} = \begin{pmatrix} {}^1b \\ ... \\ {}^mb \end{pmatrix}$$

$${}^jA_i\ {}^ix = {}^jb \tag{8.6}$$

$$\begin{matrix} {}^1A_1\ {}^1x & +... & +{}^1A_n\ {}^nx & = {}^1b \\ ... & ... & ... & ... \\ {}^1A_m\ {}^1x & +... & +{}^mA_n\ {}^nx & = {}^mb \end{matrix}$$

□

**Пример 8.4.** Пусть $\overline{V}$ - $D_*{}^*$-векторное пространство и ${}^*$-строка

$$\overline{\overline{A}} = \begin{pmatrix} {}^1\overline{A} \\ ... \\ {}^m\overline{A} \end{pmatrix} = ({}^j\overline{A}, j \in J)$$

задаёт множество векторов. Чтобы ответить на вопрос, или вектор $\overline{b} \in \mathrm{span}({}^j\overline{A}, j \in J)$ , мы запишем линейное уравнение

$$\overline{b} = x_*{}^*\overline{A} \tag{8.7}$$

где

$$x = \begin{pmatrix} x_1 & ... & x_m \end{pmatrix}$$

${}_*$-строка неизвестных коэффициентов разложения. $\overline{b} \in \mathrm{span}({}^j\overline{A}, j \in J)$ , если уравнение (8.7) имеет решение. Предположим, что $\overline{\overline{f}} = ({}^i\overline{f}, i \in I)$ - $D_*{}^*$-базис. Тогда векторы $\overline{b}$, ${}^j\overline{A}$ имеют разложение

$$\overline{b} = b_*{}^*\overline{f} \tag{8.8}$$

$${}^j\overline{A} = {}^jA_*{}^*\overline{f} \tag{8.9}$$

Если мы подставим (8.8) и (8.9) в (8.7), мы получим

$$b_*{}^*\overline{f} = x_*{}^*A_*{}^*\overline{f} \tag{8.10}$$

Применяя теорему 6.3 к (8.10), мы получим **система $D_*{}^*$-линейных уравнений**[7]

$$x_*{}^*A = b \tag{8.11}$$

[7] Читая систему ${}_*{}^*D$-линейных уравнений (8.5) снизу вверх и слева направо, мы получим систему $D^*{}_*$-линейных уравнений (8.11).

Мы можем записать систему $D_*{}^*$-линейных уравнений (8.11) в одной из следующих форм

$$\begin{pmatrix} x_1 & ... & x_m \end{pmatrix} {}_*{}^* \begin{pmatrix} {}^1A_1 & ... & {}^1A_n \\ ... & ... & ... \\ {}^mA_1 & ... & {}^mA_n \end{pmatrix} = \begin{pmatrix} b_1 & ... & b_n \end{pmatrix}$$

$$x_j\ {}^jA_i = b_i \tag{8.12}$$

$$\begin{matrix} x_1\ {}^1A_1 & ... & x_1\ {}^1A_n \\ + & ... & + \\ ... & ... & ... \\ + & ... & + \\ x_m\ {}^mA_1 & ... & x_m\ {}^mA_n \\ = & ... & = \\ b_1 & ... & b_n \end{matrix}$$

□

**Определение 8.5.** Если $n \times n$ матрица $A$ имеет ${}_*{}^*$-обратную матрицу, мы будем называть такую матрицу **${}_*{}^*$-невырожденной матрицей**. В противном случае, мы будем называть такую матрицу **${}_*{}^*$-вырожденной матрицей**. □

**Определение 8.6.** Предположим, что $A$ - ${}_*{}^*$-невырожденная матрица. Мы будем называть соответствующую систему ${}_*{}^*D$-линейных уравнений

$$A_*{}^*x = b \tag{8.13}$$

**невырожденной системой ${}_*{}^*D$-линейных уравнений**. □

**Теорема 8.7.** *Решение невырожденной системы ${}_*{}^*D$-линейных уравнений (8.13) определено однозначно и может быть записано в любой из следующих форм*[8]

$$x = A^{-1_*{}^*}{}_*{}^*b \tag{8.14}$$

$$x = \overline{H}\det(A, {}_*{}^*)\,{}_*{}^*b \tag{8.15}$$

*Доказательство.* Умножая обе части равенства (8.13) слева на $A^{-1_*{}^*}$, мы получим (8.14). Пользуясь определением [6]-(3.12), мы получим (8.15). Решение системы единственно в силу теоремы [6]-(2.16). □

## 9. Ранг матрицы

**Определение 9.1.** Мы будем называть матрицу[9] ${}^SA_T$ минором порядка $k$. □

---

[8]Мы можем найти решение системы (8.13) в теореме [4]-1.6.1. Я повторяю это утверждение, так как я слегка изменил обозначения.

[9]Мы делаем следующие предположения в этом разделе

- $i \in M$, $|M| = m$, $j \in N$, $|N| = n$.
- $A = ({}^iA_j)$ - произвольная матрица.
- $k$, $s \in S \supseteq M$, $l$, $t \in T \supseteq N$, $k = |S| = |T|$.
- $p \in M \setminus S$, $r \in N \setminus T$.

**Определение 9.2.** Если минор ${}^{S}A_{T}$ - ${}_{*}{}^{*}$-невырожденная матрица, то мы будем говорить, что ${}_{*}{}^{*}$-ранг матрицы $A$ не меньше, чем $k$. ${}_{*}{}^{*}$**-ранг матрицы** $A$

$$\mathrm{rank}_{*}{}^{*}\,A$$

- это максимальное значение $k$. Мы будем называть соответствующий минор ${}_{*}{}^{*}$**-главным минором**. □

**Теорема 9.3.** *Пусть матрица $A$ - ${}_{*}{}^{*}$-вырожденная матрица и минор ${}^{S}A_{T}$ - главный минор, тогда*

$$(9.1)\qquad {}^{p}\det\left({}^{S\cup\{p\}}A_{T\cup\{r\}},{}_{*}{}^{*}\right)_{r}=0$$

*Доказательство.* Чтобы понять, почему минор ${}^{S\cup\{p\}}A_{T\cup\{r\}}$ не имеет ${}_{*}{}^{*}$-обратной матрицы,[10] мы предположим, что он имеет ${}_{*}{}^{*}$-обратную матрицу и запишем соответствующую систему [6]-(3.2), [6]-(3.3), полагая $i=r$, $j=p$ и попробуем решить эту систему. Положив

$$(9.2)\qquad B={}^{S\cup\{p\}}A_{T\cup\{r\}}$$

мы получим систему

$$(9.3)\qquad {}^{S}B_{T}{}_{*}{}^{*}{}^{T}B^{-1_{*}{}^{*}}{}_{p}+{}^{S}B_{r}\ {}^{r}B^{-1_{*}{}^{*}}{}_{p}=0$$

$$(9.4)\qquad {}^{p}B_{T}{}_{*}{}^{*}{}^{T}B^{-1_{*}{}^{*}}{}_{p}+{}^{p}B_{r}\ {}^{r}B^{-1_{*}{}^{*}}{}_{p}=1$$

Мы умножим (9.3) на $({}^{S}B_{T})^{-1_{*}{}^{*}}$

$$(9.5)\qquad {}^{T}B^{-1_{*}{}^{*}}{}_{p}+({}^{S}B_{T})^{-1_{*}{}^{*}}{}_{*}{}^{*}{}^{S}B_{r}\ {}^{r}B^{-1_{*}{}^{*}}{}_{p}=0$$

Теперь мы можем подставить (9.5) в (9.4)

$$(9.6)\qquad -\,{}^{p}B_{T}{}_{*}{}^{*}({}^{S}B_{T})^{-1_{*}{}^{*}}{}_{*}{}^{*}{}^{S}B_{r}\ {}^{r}B^{-1_{*}{}^{*}}{}_{p}+{}^{p}B_{r}\ {}^{r}B^{-1_{*}{}^{*}}{}_{p}=1$$

Из (9.6) следует

$$(9.7)\qquad ({}^{p}B_{r}-{}^{p}B_{T}{}_{*}{}^{*}({}^{S}B_{T})^{-1_{*}{}^{*}}{}_{*}{}^{*}{}^{S}B_{r})\ {}^{r}B^{-1_{*}{}^{*}}{}_{p}=1$$

Выражение в скобках является квазидетерминантом ${}^{p}\det(B,{}_{*}{}^{*})_{r}$. Подставляя это выражение в (9.7), мы получим

$$(9.8)\qquad {}^{p}\det(B,{}_{*}{}^{*})_{r}\ {}^{r}B^{-1_{*}{}^{*}}{}_{p}=1$$

Тем самым мы доказали, что квазидетерминант ${}^{p}\det(B,{}_{*}{}^{*})_{r}$ определён и условие его обращения в 0 необходимое и достаточное условие вырожденности матрицы $B$. Теорема доказана в силу соглашения (9.2). □

**Теорема 9.4.** *Предположим, что $A$ - матрица,*

$$\mathrm{rank}_{*}{}^{*}\,A=k<m$$

[10]Естественно ожидать связь между ${}_{*}{}^{*}$-вырожденностью матрицы и её ${}_{*}{}^{*}$-квазидетерминантом, подобную связи, известной в коммутативном случае. Однако ${}_{*}{}^{*}$-квазидетерминант определён не всегда. Например, если ${}_{*}{}^{*}$-обратная матрица имеет слишком много элементов, равных 0. Как следует из этой теоремы, не определён ${}_{*}{}^{*}$-квазидетерминант также в случае, когда ${}_{*}{}^{*}$-ранг матрицы меньше $n-1$.

*и* ${}^{S}A_T$ *-* ${}_*{}^*$*-главный минор. Тогда* ${}_*$*-строка* ${}^{p}A$ *является* $D_*{}^*$*-линейной комбинацией* ${}_*$*-строк* ${}^{S}A$*.*

$$ {}^{M\setminus S}A = R_*{}^*{}^{S}A \tag{9.9}$$

$$ {}^{p}A = {}^{p}R_*{}^*{}^{S}A \tag{9.10}$$

$$ {}^{p}A_b = {}^{p}R_s\ {}^{s}A_b \tag{9.11}$$

*Доказательство.* Если число ${}^*$-строк - $k$, то, полагая, что ${}_*$-строка ${}^{p}A$ - $D_*{}^*$-линейная комбинация (9.10) ${}_*$-строк ${}^{s}A$ с коэффициентами ${}^{p}R_s$, мы получим систему (9.11). Согласно теореме 8.7 эта система имеет единственное решение[11] и оно нетривиально потому, что все ${}_*{}^*$-квазидетерминанты отличны от 0.

Остаётся доказать утверждение в случае, когда число ${}^*$-строк больше чем $k$. Пусть нам даны ${}_*$-строка ${}^{p}A$ и ${}^*$-строка $A_r$. Согласно предположению минор ${}^{S\cup\{p\}}A_{T\cup\{r\}}$ - ${}_*{}^*$-вырожденная матрица и его ${}_*{}^*$-квазидетерминант

$$ {}^{p}\det\left({}^{S\cup\{p\}}A_{T\cup\{r\}},{}_*{}^*\right)_r = 0 \tag{9.12}$$

Согласно [6]-(3.14) равенство (9.12) имеет вид

$$ {}^{p}A_r - {}^{p}A_{T*}{}^*(({}^{S}A_T)^{-1_*{}^*})_*{}^*{}^{S}A_r = 0 $$

Матрица

$$ {}^{p}R = {}^{p}A_{T*}{}^*(({}^{S}A_T)^{-1_*{}^*}) \tag{9.13}$$

не зависят от $r$. Следовательно, для любых $r \in N \setminus T$

$$ {}^{p}A_r = {}^{p}R_*{}^*{}^{S}A_r \tag{9.14}$$

Из равенства

$$ (({}^{S}A_T)^{-1_*{}^*})_*{}^*{}^{S}A_l = {}^{T}\delta_l $$

следует, что

$$ {}^{p}A_l = {}^{p}A_{T*}{}^*{}^{T}\delta_l = {}^{p}A_{T*}{}^*(({}^{S}A_T)^{-1_*{}^*})_*{}^*{}^{S}A_l \tag{9.15}$$

Подставляя (9.13) в (9.15), мы получим

$$ {}^{p}A_l = {}^{p}R_*{}^*{}^{S}A_l \tag{9.16}$$

(9.14) и (9.16) завершают доказательство. □

**Следствие 9.5.** *Предположим, что* $A$ *- матрица,* $\mathrm{rank}_{*}{}^* A = k < m$*. Тогда* ${}_*$*-строки матрицы* $D_*{}^*$*-линейно зависимы.*

$$ \lambda_*{}^*A = 0 \tag{9.17}$$

*Доказательство.* Предположим, что ${}_*$-строка ${}^{p}A$ - $D_*{}^*$-линейная комбинация (9.10). Мы положим $\lambda_p = -1$, $\lambda_s = {}^{p}R_s$ и остальные $\lambda_c = 0$. □

**Теорема 9.6.** *Пусть* $({}^{i}\overline{A}, i \in M, |M| = m)$ *семейство* $D_*{}^*$*-линейно независимых векторов. Тогда* ${}_*{}^*$*-ранг их координатной матрицы равен* $m$*.*

[11] Мы положим, что неизвестные переменные здесь - это $x_s = {}^{p}R_s$

*Доказательство.* Согласно конструкции, изложенной в примере 6.7, координатная матрица семейства векторов $({}^i\overline{A})$ относительно базиса $\overline{\overline{e}}$ состоит из ${}_*$-строк, являющихся координатными матрицами векторов ${}^i\overline{A}$ относительно базиса $\overline{\overline{e}}$. Поэтому ${}_*{}^*$-ранг этой матрицы не может превышать $m$.

Допустим ${}_*{}^*$-ранг координатной матрицы меньше $m$. Согласно следствию 9.5 ${}_*$-строки матрицы $D_*{}^*$-линейно зависимы. Если мы умножим обе части равенства (9.17) на ${}^*$-строку $\overline{e}$ и положим $c = \lambda_*{}^*A$, то мы получим, что $D_*{}^*$-линейная комбинация

$$c_*{}^*\overline{e} = 0$$

векторов базиса равна 0. Это противоречит утверждению, что векторы $\overline{e}$ образуют базис. Утверждение теоремы доказано. □

**Теорема 9.7.** *Предположим, что $A$ - матрица,*

$$\operatorname{rank}_{*^*} A = k < n$$

*и ${}^SA_T$ - ${}_*{}^*$-главный минор. Тогда ${}^*$-строка $A_r$ является ${}_*{}^*D$-линейной композицией ${}^*$-строк $A_t$*

(9.18) $$A_{N\setminus T} = A_{T*}{}^*R$$

(9.19) $$A_r = A_{T*}{}^*R_r$$

(9.20) $${}^aA_r = {}^aA_t\ {}^tR_r$$

*Доказательство.* Если число ${}_*$-строк - $k$, то, полагая, что ${}^*$-строка $A_r$ - ${}_*{}^*D$-линейная комбинация (9.19) ${}^*$-строк $A_t$ с коэффициентами ${}^tR_r$, мы получим систему (9.20). Согласно теореме 8.7 эта система имеет единственное решение[12] и оно нетривиально потому, что все ${}_*{}^*$-квазидетерминанты отличны от 0.

Остаётся доказать утверждение в случае, когда число ${}_*$-строк больше, чем $k$. Пусть нам даны ${}^*$-строка $A_r$ и ${}_*$-строка ${}^pA$. Согласно предположению минор ${}^{S\cup\{p\}}A_{T\cup\{r\}}$ - ${}_*{}^*$-вырожденная матрица и его ${}_*{}^*$-квазидетерминант

(9.21) $${}^p\det\left({}^{S\cup\{p\}}A_{T\cup\{r\}}, {}_*{}^*\right)_r = 0$$

Согласно [6]-(3.14) (9.21) имеет вид

$${}^pA_r - {}^pA_{T*}{}^*(({}^SA_T)^{-1_*{}^*})_*{}^*{}^SA_r = 0$$

Матрица

(9.22) $$R_r = (({}^SA_T)^{-1_*{}^*})_*{}^*{}^SA_r$$

не зависит от $p$, Следовательно, для любых $p \in M \setminus S$

(9.23) $${}^pA_r = {}^pA_{T*}{}^*R_r$$

Из равенства

$${}^kA_{T*}{}^*(({}^SA_T)^{-1_*{}^*})_s = {}^k\delta_s$$

следует, что

(9.24) $${}^kA_r = {}^k\delta_{S*}{}^*{}^SA_r = {}^kA_{T*}{}^*(({}^SA_T)^{-1_*{}^*})^s{}_*{}^*{}^SA_r$$

Подставляя (9.22) в (9.24), мы получим

(9.25) $${}^kA_r = {}^kA_{T*}{}^*R_r$$

(9.23) и (9.25) завершают доказательство. □

[12] Мы положим, что неизвестные переменные здесь - это ${}^tx = {}^tR_r$

**Следствие 9.8.** *Предположим, что $A$ - матрица,* $\mathrm{rank}_{*}{}^{*} A = k < m$*. Тогда ${}^*$-строки матрицы ${}_*{}^*D$-линейно зависимы.*

$$A_*{}^*\lambda = 0$$

*Доказательство.* Предположим, что ${}^*$-строка $A_r$ - правая линейная комбинация (9.19). Мы положим ${}^r\lambda = -1$, ${}^t\lambda_= {}^tR_r$ и остальные ${}^c\lambda = 0$. □

Опираясь на теорему [6]-3.8, мы можем записать подобные утверждения для ${}^*{}_*$-ранга матрицы.

**Теорема 9.9.** *Предположим, что $A$ - матрица,*

$$\mathrm{rank}_{{}^*{}_*} A = k < m$$

*и ${}_TA^S$ - ${}^*{}_*$-главный минор. Тогда ${}_*$-строка $A^p$ является ${}^*{}_*D$-линейной композицией ${}_*$-строк $A^s$*

(9.26) $$A^{M\setminus S} = A^S{}^*{}_*R$$

(9.27) $$A^p = A^S{}^*{}_*R^p$$

(9.28) $${}_bA^p = {}_bA^s\ {}_sR^p$$

**Следствие 9.10.** *Предположим, что $A$ - матрица,* $\mathrm{rank}_{{}^*{}_*} A = k < m$*. Тогда ${}_*$-строки матрицы ${}^*{}_*D$-линейно зависимы.*

$$A^*{}_*\lambda = 0$$

**Теорема 9.11.** *Предположим, что $A$ - матрица,*

$$\mathrm{rank}_{{}^*{}_*} A = k < n$$

*и ${}_TA^S$ - ${}^*{}_*$-главный минор. Тогда ${}^*$-строка ${}_rA$ является $D^*{}_*$-линейной композицией ${}^*$-строк ${}_tA$*

(9.29) $${}_{N\setminus T}A = R^*{}_*{}_TA$$

(9.30) $${}_rA = {}_rR^*{}_*{}_TA$$

(9.31) $${}_rA^a = {}_rR^t\ {}_tA^a$$

**Следствие 9.12.** *Предположим, что $A$ - матрица,* $\mathrm{rank}_{{}^*{}_*} a = k < m$*. Тогда ${}^*$-строки матрицы $D^*{}_*$-линейно зависимы.*

$$\lambda^*{}_*A = 0$$

## 10. Система ${}_*{}^*D$-линейных уравнений

**Определение 10.1.** Предположим, что[9] $A$ - матрица системы $D_*{}^*$-линейных уравнений (8.12). Мы будем называть матрицу

(10.1) $$\begin{pmatrix} {}^jA_i \\ b_i \end{pmatrix} = \begin{pmatrix} {}^1A_1 & ... & {}^1A_n \\ ... & ... & ... \\ {}^mA_1 & ... & {}^mA_n \\ b_1 & ... & b_n \end{pmatrix}$$

**расширенной матрицей** этой системы. □

**Определение 10.2.** Предположим, что[9] $A$ - матрица системы ${}_*{}^*D$-линейных уравнений (8.6). Мы будем называть матрицу

$$
\begin{pmatrix} {}^jA_i & {}^jb \end{pmatrix} = \begin{pmatrix} {}^1A_1 & ... & {}^1A_n & {}^1b \\ ... & ... & ... & ... \\ {}^mA_1 & ... & {}^mA_n & {}^mb \end{pmatrix} \tag{10.2}
$$

**расширенной матрицей** этой системы. □

**Теорема 10.3.** *Система ${}_*{}^*D$-линейных уравнений (8.6) имеет решение тогда и только тогда, когда*

$$
\operatorname{rank}_{{}_*{}^*}({}^jA_i) = \operatorname{rank}_{{}_*{}^*}\begin{pmatrix} {}^jA_i & {}^jb \end{pmatrix} \tag{10.3}
$$

*Доказательство.* Пусть ${}^SA_T$ - ${}_*{}^*$-главный минор матрицы $A$.

Пусть система ${}_*{}^*D$-линейных уравнений (8.6) имеет решение ${}^ix = {}^id$. Тогда

$$
A_*{}^*d = b \tag{10.4}
$$

Уравнение (10.4) может быть записано в форме

$$
A_{T*}{}^{*T}d + A_{N\setminus T*}{}^{*N\setminus T}d = b \tag{10.5}
$$

Подставляя (9.18) в (10.5), мы получим

$$
A_{T*}{}^{*T}d + A_{T*}{}^*R_*{}^{*N\setminus T}d = b \tag{10.6}
$$

Из (10.6) следует, что ${}_*$-строка $b$ является ${}_*{}^*D$-линейной комбинацией ${}_*$-строк $A_T$

$$
A_{T*}{}^*({}^Td + R_*{}^{*N\setminus T}d) = b
$$

Это эквивалентно уравнению (10.3).

Нам осталось доказать, что существование решения ${}_*{}^*D$-системы линейных уравнений (8.6) следует из (10.3). Истинность (10.3) означает, что ${}^SA_T$ - так же ${}_*{}^*$-главный минор расширенной матрицы. Из теоремы 9.7 следует, что ${}_*$-строка $b$ является ${}_*{}^*D$-линейной композицией ${}_*$-строк $A_T$

$$
b = A_{T*}{}^{*T}R
$$

Полагая ${}^rR = 0$ , мы получим

$$
b = A_*{}^*R
$$

Следовательно, мы нашли, по крайней мере, одно решение системы ${}_*{}^*D$-линейных уравнений (8.6). □

**Теорема 10.4.** *Предположим, что (8.6) - система ${}_*{}^*D$-линейных уравнений, удовлетворяющих (10.3). Если $\operatorname{rank}_{{}_*{}^*} A = k \le m$, то решение системы зависит от произвольных значений $m - k$ переменных, не включённых в ${}_*{}^*$-главный минор.*

*Доказательство.* Пусть ${}^SA_T$ - ${}_*{}^*$-главный минор матрицы $a$. Предположим, что

$$
{}^pA_*{}^*x = {}^pb \tag{10.7}
$$

уравнение с номером $p$. Применяя теорему 9.4 к расширенной матрице (10.1), мы получим

$$ {}^pA = {}^pR_*{}^{*S}A \tag{10.8} $$

$$ {}^pb = {}^pR_*{}^{*S}b \tag{10.9} $$

Подставляя (10.8) и (10.9) в (10.7), мы получим

$$ {}^pR_*{}^{*S}A_*{}^*x = {}^pR_*{}^{*S}b \tag{10.10} $$

(10.10) означает, что мы можем исключить уравнение (10.7) из системы (8.6) и новая система эквивалентна старой. Следовательно, число уравнений может быть уменьшено до $k$.

В этом случае у нас есть два варианта. Если число переменных также равно $k$, то согласно теореме 8.7 система имеет единственное решение (8.15). Если число переменных $m > k$, то мы можем передвинуть $m-k$ переменных, которые не включены в ${}_*{}^*$-главный минор в правой части. Присваивая произвольные значения этим переменным, мы определяем значение правой части и для этого значения мы получим единственное решение согласно теореме 8.7. □

**Следствие 10.5.** *Система ${}_*{}^*D$-линейных уравнений (8.6) имеет единственное решение тогда и только тогда, когда её матрица невырожденная.* □

**Теорема 10.6.** *Решения однородной системы ${}_*{}^*D$-линейных уравнений*

$$ A_*{}^*x = 0 \tag{10.11} $$

*порождают ${}_*{}^*D$-векторное пространство.*

*Доказательство.* Пусть $\overline{X}$ - множество решений системы ${}_*{}^*D$-линейных уравнений (10.11). Предположим, что $x = ({}^ax) \in \overline{X}$ и $y = ({}^ay) \in \overline{X}$. Тогда

$$ x^a\ {}_aa^b = 0 $$

$$ y^a\ {}_aa^b = 0 $$

Следовательно,

$$ {}^iA_j\ ({}^jx + {}^jy) = {}^iA_j\ {}^jx + {}^iA_j\ {}^jy = 0 $$

$$ x + y = ({}^jx + {}^jy) \in \overline{X} $$

Таким же образом мы видим

$$ {}^iA_j\ ({}^jxb) = ({}^iA_j\ {}^jx)b = 0 $$

$$ xb = ({}^jxb) \in \overline{X} $$

Согласно определению 4.4 $\overline{X}$ является ${}_*{}^*D$-векторным пространством. □

## 11. Невырожденная матрица

Предположим, что нам дана $n \times n$ матрица $A$. Следствия 9.5 и 9.8 говорят нам, что если $\mathrm{rank}_{*}{}^{*} A < n$, то ${}_*$-строки $D_*{}^*$-линейно зависимы и ${}^*$-строки ${}_*{}^*D$-линейно зависимы.[13]

**Теорема 11.1.** *Пусть $A$ - $n \times n$ матрица и ${}^*$-строка $A_r$ - ${}_*{}^*D$-линейная комбинация других ${}^*$-строк. Тогда* $\mathrm{rank}_{*}{}^{*} A < n$.

[13] Это утверждение похоже на утверждение [5]-1.2.5.

*Доказательство.* Утверждение, что ${}^*$-строка $A_r$ является ${}_*{}^*D$-линейной комбинацией других ${}^*$-строк, означает, что система ${}_*{}^*D$-линейных уравнений

$$A_r = A_{[r]*}{}^*\lambda$$

имеет, по крайней мере, одно решение. Согласно теореме 10.3

$$\mathrm{rank}_{*^*} A = \mathrm{rank}_{*^*} A_{[r]}$$

Так как число ${}^*$-строк меньше, чем $n$, то $\mathrm{rank}_{*^*} A[r] < n$. □

**Теорема 11.2.** *Пусть $A$ - $n \times n$ матрица и ${}_*$-строка ${}^pA$ - $D_*{}^*$-линейная комбинация других ${}_*$-строк. Тогда* $\mathrm{rank}_{*^*} A < n$.

*Доказательство.* Доказательство утверждения похоже на доказательство теоремы 11.1 □

**Теорема 11.3.** *Предположим, что $A$ и $B$ - $n \times n$ матрицы и*

$$C = A_*{}^*B \tag{11.1}$$

*$C$ - ${}_*{}^*$-вырожденная матрица тогда и только тогда, когда либо матрица $A$, либо матрица $B$ - ${}_*{}^*$-вырожденная матрица.*

*Доказательство.* Предположим, что матрица $B$ - ${}_*{}^*$-вырожденная. Согласно теореме 9.7 ${}^*$-строки матрицы $B$ ${}_*{}^*D$-линейно зависимы. Следовательно,

$$0 = B_*{}^*\lambda \tag{11.2}$$

где $\lambda \neq 0$. Из (11.1) и (11.2) следует, что

$$C_*{}^*\lambda = A_*{}^*B_*{}^*\lambda = 0$$

Согласно теореме 11.1 матрица $C$ является ${}_*{}^*$-вырожденной.

Предположим, что матрица $B$ - не ${}_*{}^*$-вырожденная, но матрица $A$ - ${}_*{}^*$-вырожденная. Согласно теореме 9.4 ${}^*$-строки матрицы $A$ ${}_*{}^*D$-линейно зависимы. Следовательно,

$$0 = A_*{}^*\mu \tag{11.3}$$

где $\mu \neq 0$. Согласно теореме 8.7 система

$$B_*{}^*\lambda = \mu$$

имеет единственное решение, где $\lambda \neq 0$. Следовательно,

$$C_*{}^*\lambda = A_*{}^*B_*{}^*\lambda = A_*{}^*\mu = 0$$

Согласно теореме 11.1 матрица $C$ является ${}_*{}^*$-вырожденной.

Предположим, что матрица $C$ - ${}_*{}^*$-вырожденная матрица. Согласно теореме 9.4 ${}_*$-строки матрицы $C$ ${}_*{}^*D$-линейно зависимы. Следовательно,

$$0 = C_*{}^*\lambda \tag{11.4}$$

где $\lambda \neq 0$. Из (11.1) и (11.4) следует, что

$$0 = A_*{}^*B_*{}^*\lambda$$

Если

$$0 = B_*{}^*\lambda$$

выполнено, то матрица $B$ - ${}_*{}^*$-вырожденная. Предположим, что матрица $B$ не ${}_*{}^*$-вырожденная. Положим

$$\mu = B_*{}^*\lambda$$

где $\mu \neq 0$. Тогда

$$0 = A_*{}^*\mu \tag{11.5}$$

Из (11.5) следует, что матрица $A$ - ${}_*{}^*$-вырожденная. □

Опираясь на теорему [6]-3.8, мы можем записать подобные утверждения для $D^*{}_*$-линейной комбинации ${}_*$-строк или ${}^*{}_*D$-линейной комбинации ${}^*$-строк и ${}^*{}_*$-квазидетерминанта.

**Теорема 11.4.** *Пусть $A$ - $n\times n$ матрица и ${}^*$-строка ${}_rA$ - ${}^*{}_*D$-линейная комбинация других ${}^*$-строк. Тогда* $\mathrm{rank}_{{}^*{}_*} A < n$.

**Теорема 11.5.** *Пусть $A$ - $n\times n$ матрица и ${}_*$-строка $A^p$ - $D^*{}_*$-линейная комбинация других ${}_*$-строк. Тогда* $\mathrm{rank}_{{}^*{}_*} A < n$.

**Теорема 11.6.** *Предположим, что $A$ и $B$ - $n\times n$ матрицы и $C = A^*{}_*B$. $C$ - ${}^*{}_*$-вырожденная матрица тогда и только тогда, когда либо матрица $A$, либо матрица $B$ - ${}^*{}_*$-вырожденная матрица.*

**Определение 11.7.** **${}_*{}^*$-матричная группа** $GL_n^{{}_*{}^*D}$ - это группа ${}_*{}^*$-невырожденных матриц, где мы определяем ${}_*{}^*$-произведение матриц [6]-(2.1) и ${}_*{}^*$-обратную матрица $A^{-1^*{}_*}$. □

**Определение 11.8.** **${}^*{}_*$-матричная группа** $GL_n^{{}^*{}_*D}$ - это группа ${}^*{}_*$-невырожденных матриц где мы определяем ${}^*{}_*$-произведение матриц [6]-(2.2) и ${}^*{}_*$-обратную матрицу $A^{-1^*{}_*}$. □

**Теорема 11.9.**

$$GL_n^{{}_*{}^*D} \neq GL_n^{{}^*{}_*D}$$

*Замечание* 11.10. Из теоремы [6]-(3.10) следует, что существуют матрицы, которые ${}^*{}_*$-невырожденны и ${}_*{}^*$-невырожденны. Теорема 11.9 означает, что множества ${}^*{}_*$-невырожденных матриц и ${}_*{}^*$-невырожденных матриц не совпадают. Например, существует такая ${}_*{}^*$-невырожденная матрица, которая ${}^*{}_*$-вырожденная матрица. □

*Доказательство.* Это утверждение достаточно доказать для $n = 2$. Предположим, что каждая ${}^*{}_*$-вырожденная матрица

$$A = \begin{pmatrix} {}^1A_1 & {}^2A_1 \\ {}^1A_2 & {}^2A_2 \end{pmatrix} \tag{11.6}$$

является ${}_*{}^*$-вырожденной матрицей. Из теоремы 9.7 и теоремы 9.7 следует, что ${}_*{}^*$-вырожденная матрица удовлетворяет условию

$${}^1A_2 = b\ {}^1A_1 \tag{11.7}$$

$${}^2A_2 = b\ {}^2A_1 \tag{11.8}$$

$${}^2A_1 = {}^1A_1c \tag{11.9}$$

$${}^2A_2 = {}^1A_2c \tag{11.10}$$

Если мы подставим (11.9) в (11.8), мы получим

$${}^2A_2 = b\ {}^1A_1\ c$$

$b$ и $c$ - произвольные элементы тела $D$ и ${}_*{}^*$-вырожденная матрица (11.6) имеет вид ($d={}^1A_1$)

$$A=\begin{pmatrix}d&dc\\bd&bdc\end{pmatrix}\tag{11.11}$$

Подобным образом мы можем показать, что ${}^*{}_*$-вырожденная матрица имеет вид

$$A=\begin{pmatrix}d&c'd\\db'&c'db'\end{pmatrix}\tag{11.12}$$

Из предположения следует, что (11.12) и (11.11) представляют одну и ту же матрицу. Сравнивая (11.12) и (11.11), мы получим, что для любых $d,c\in D$ существует такое $c'\in D$, которое не зависит от $d$ и удовлетворяет уравнению

$$dc=c'd$$

Это противоречит утверждению, что $D$ - тело. □

**Пример 11.11.** В случае тела кватернионов мы положим $b=1+k$, $c=j$, $d=k$. Тогда мы имеем

$$A=\begin{pmatrix}k&kj\\(1+k)k&(1+k)kj\end{pmatrix}=\begin{pmatrix}k&-i\\k-1&-i-j\end{pmatrix}$$

$$\begin{aligned}{}^2\det(A,{}_*{}^*)_2&={}^2A_2-{}^2A_1({}^1A_1)^{-1}\,{}^1A_2\\&=-i-j-(k-1)(k)^{-1}(-i)=-i-j-(k-1)(-k)(-i)\\&=-i-j-kki+ki=-i-j+i+j\\&=0\end{aligned}$$

$$\begin{aligned}{}_1\det(A,{}_*{}^*)^1&={}_1A^1-{}_1A^2({}_2A^2)^{-1}\,{}_2A^1\\&=k-(k-1)(-i-j)^{-1}(-i)=k-(k-1)\frac{1}{2}(i+j)(-i)\\&=k+\frac{1}{2}((k-1)i+(k-1)j)i=k+\frac{1}{2}(ki-i+kj-j)i\\&=k+\frac{1}{2}(j-i-i-j)i=k-ii\\&=k+1\end{aligned}$$

$$\begin{aligned}{}_1\det(A,{}^*{}_*)^2&={}_1A^2-{}_1A^1({}_2A^1)^{-1}\,{}_2A^2\\&=k-1-k(-i)^{-1}(-i-j)=k-1+ki(i+j)\\&=k-1+j(i+j)=k-1+ji+jj\\&=k-1-k-1\\&=-2\end{aligned}$$

$$\begin{aligned}
{}_2\det(A,{}^*{}_*)^1 &= {}_2A^1 - {}_2A^2({}_1A^2)^{-1}\ {}_1A^1\\
&= (-i) - (-i-j)(k-1)^{-1}k = -i + (i+j)\frac{1}{2}(-k-1)k\\
&= -i - \frac{1}{2}(i+j)(k+1)k = -i - \frac{1}{2}(ik+i+jk+j)k\\
&= -i - \frac{1}{2}(-j+i+i+j)k = -i - ik\\
&= -i + j
\end{aligned}$$

$$\begin{aligned}
{}_2\det(A,{}^*{}_*)^2 &= {}_2A^2 - {}_2A^1({}_1A^1)^{-1}\ {}_1A^2\\
&= -i-j-(-i)(k)^{-1}(k-1) = -i-j+i(-k)(k-1)\\
&= -i-j+j(k-1) = -i-j+jk-j = -i-j+i-j\\
&= -2j
\end{aligned}$$

Система ${}_*{}^*D$-линейных уравнений

$$\begin{pmatrix} k & -i\\ k-1 & -i-j\end{pmatrix} {}_*{}^* \begin{pmatrix}{}^1x\\ {}^2x\end{pmatrix} = \begin{pmatrix}{}^1b\\ {}^2b\end{pmatrix} \tag{11.13}$$

имеет ${}_*{}^*$-вырожденную матрицу. Мы можем записать систему ${}_*{}^*D$-линейных уравнений (11.13) в виде

$$\left\{\begin{array}{rl} k\ {}^1x - & i\ {}^2x = {}^1b\\ (k-1)\ {}^1x - & (i+j)\ {}^2x = {}^2b\end{array}\right.$$

Система ${}^*{}_*D$-линейных уравнений

$$\begin{pmatrix} k & -i\\ k-1 & -i-j\end{pmatrix} {}^*{}_* \begin{pmatrix}{}_1x & {}_2x\end{pmatrix} = \begin{pmatrix}{}_1b & {}_2b\end{pmatrix} \tag{11.14}$$

имеет ${}^*{}_*$-невырожденную матрицу. Мы можем записать систему ${}^*{}_*D$-линейных уравнений (11.14) в виде

$$\left\{\begin{array}{rr} k\ {}_1x & - \quad i\ {}_1x\\ +\ (k-1)\ {}_2x & -\ (i+j)\ {}_2x\\ = \quad {}_1b & = \quad {}_2b\end{array}\right. \qquad \left\{\begin{array}{l} k\ {}_1x + (k-1)\ {}_2x = {}_1b\\ -i\ {}_1x - (i+j)\ {}_2x = {}_2b\end{array}\right.$$

Система $D_*{}^*$-линейных уравнений

$$\begin{pmatrix} x_1 & x_2\end{pmatrix} {}_*{}^* \begin{pmatrix} k & -i\\ k-1 & -i-j\end{pmatrix} = \begin{pmatrix} b_1 & b_2\end{pmatrix} \tag{11.15}$$

имеет ${}_*{}^*$-вырожденную матрицу. Мы можем записать систему $D_*{}^*$-линейных уравнений (11.15) в виде

$$\left\{\begin{array}{ll} x_1k & -x_1i\\ +\ x_2(k-1) & -x_2(i+j)\\ =b_1 & =b_2\end{array}\right. \qquad \left\{\begin{array}{l} x_1k \ +x_2(k-1) = b_1\\ -x_1i\ -x_2(i+j) = b_2\end{array}\right.$$

Система $D^*{}_*$-линейных уравнений

$$\begin{pmatrix} x^1\\ x^2\end{pmatrix} {}^*{}_* \begin{pmatrix} k & -i\\ k-1 & -i-j\end{pmatrix} = \begin{pmatrix}{}^1b\\ {}^2b\end{pmatrix} \tag{11.16}$$

имеет ${}^*{}_*$-невырожденную матрицу. Мы можем записать систему $D^*{}_*$-линейных уравнений (11.16) в виде

$$\begin{cases} k\ {}^1x - \quad i\ {}^2x = {}^1b \\ (k-1)\ {}^1x - (i+j)\ {}^2x = {}^2b \end{cases}$$

□

## 12. Размерность ${}_*{}^*D$-векторного пространства

**Теорема 12.1.** *Пусть $\overline{V}$ - ${}_*{}^*D$-векторное пространство. Предположим, что $\overline{V}$ имеет ${}_*{}^*D$-базисы $\overline{\overline{e}} = (\overline{e}_i, i \in I)$ и $\overline{\overline{g}} = (\overline{g}_j, j \in J)$. Если $|I|$ и $|J|$ - конечные числа, то $|I| = |J|$.*

*Доказательство.* Предположим, что $|I| = m$ и $|J| = n$. Предположим, что

$$m < n \tag{12.1}$$

Так как $\overline{\overline{e}}$ - ${}_*{}^*D$-базис, любой вектор $\overline{g}_j$, $j \in J$ имеет разложение

$$\overline{g}_j = \overline{e}_*{}^*A_j$$

Так как $\overline{\overline{g}}$ - ${}_*{}^*D$-базис

$$\lambda = 0 \tag{12.2}$$

следует из

$$\overline{g}_*{}^*\lambda = \overline{e}_*{}^*A_*{}^*\lambda = 0$$

Так как $\overline{\overline{e}}$ - ${}_*{}^*D$-базис, мы получим

$$A_*{}^*\lambda = 0 \tag{12.3}$$

Согласно (12.1) $\mathrm{rank}_{{}_*{}^*}\, A \le m$ и система (12.3) имеет больше переменных, чем уравнений. Согласно теорема 10.4 $\lambda \ne 0$. Это противоречит утверждению (12.2). Следовательно, утверждение $m < n$ неверно.

Таким же образом мы можем доказать, что утверждение $n < m$ неверно. Это завершает доказательство теоремы. □

**Определение 12.2.** Мы будем называть **размерностью ${}_*{}^*D$-векторного пространства** число векторов в базисе □

**Теорема 12.3.** *Координатная матрица ${}_*{}^*D$-базиса $\overline{\overline{g}}$ относительно ${}_*{}^*D$-базиса $\overline{\overline{e}}$ векторного пространства $\overline{V}$ является ${}_*{}^*$-невырожденной матрицей.*

*Доказательство.* Согласно лемме 9.6 ${}_*{}^*D$-ранг координатной матрицы базиса $\overline{\overline{g}}$ относительно базиса $\overline{\overline{e}}$ равен размерности векторного пространства, откуда следует утверждение теоремы. □

**Определение 12.4.** Мы будем называть взаимно однозначное отображение $\overline{A} : \overline{V} \to \overline{W}$ **${}_*{}^*D$-изоморфизмом векторных пространств**, если это отображение - ${}_*{}^*D$-линейное отображение векторных пространств. □

**Определение 12.5.** **${}_*{}^*D$-автоморфизм векторного пространства** $\overline{V}$ - это ${}_*{}^*D$-изоморфизм $\overline{A} : \overline{V} \to \overline{V}$. □

**Теорема 12.6.** *Предположим, что* $\overline{\overline{f}}$ *-* ${}_*{}^*D$*-базис в векторном пространстве* $\overline{V}$*. Тогда любой* ${}_*{}^*D$*-автоморфизм* $\overline{A}$ *векторного пространства* $\overline{V}$ *имеет вид*

$$v' = A_*{}^*v \tag{12.4}$$

*где* $A$ *-* ${}_*{}^*$*-невырожденная матрица.*

*Доказательство.* (12.4) следует из теоремы 7.3. Так как $\overline{A}$ - изоморфизм, то для каждого вектора $\overline{v}'$ существует единственный вектор $\overline{v}$ такой, что $\overline{v}' = \overline{v}_*{}^*\overline{A}$. Следовательно, система ${}_*{}^*D$-линейных уравнений (12.4) имеет единственное решение. Согласно следствию 10.5 матрица $A$ невырожденна. □

**Теорема 12.7.** *Автоморфизмы* ${}_*{}^*D$*-векторного пространства порождают группу* $GL_n^{{}_*{}^*D}$*.*

*Доказательство.* Если даны два автоморфизма $\overline{A}$ и $\overline{B}$, то мы можем записать

$$v' = A_*{}^*v$$

$$v'' = B_*{}^*v' = B_*{}^*A_*{}^*v$$

Следовательно, результирующий автоморфизм имеет матрицу $A_*{}^*B$. □

## 13. Список литературы

[1] S. Burris, H.P. Sankappanavar, A Course in Universal Algebra, Springer-Verlag (March, 1982),
eprint http://www.math.uwaterloo.ca/ snburris/htdocs/ualg.html
(The Millennium Edition)

[2] А. Г. Курош, Общая алгебра, (лекции 1969 - 70 учебного года), М., МГУ, 1970

[3] Lev V. Sabinin, Smooth Quasigroups and Loops, Kluwer Academic Publisher, 1999

[4] I. Gelfand, S. Gelfand, V. Retakh, R. Wilson, Quasideterminants,
eprint arXiv:math.QA/0208146 (2002)

[5] I.Gelfand, V.Retakh, Quasideterminants, I,
eprint arXiv:q-alg/9705026 (1997)

[6] Александр Клейн, Бикольцо матриц,
eprint arXiv:math.OA/0612111 (2007)

[7] П. Кон, Универсальная алгебра, М., Мир, 1968

# 14. Предметный указатель

# 15. Специальные символы и обозначения